\documentclass[10pt,a4paper]{article}
\usepackage{amsmath}
\usepackage{amsfonts}
\usepackage{amssymb}
\usepackage{graphicx}
\usepackage{bbm}
\usepackage{hyperref}

\usepackage{amsmath}
\usepackage{amsfonts}
\usepackage{amssymb}
\usepackage{amsrefs}
\usepackage{mathtools}
\usepackage[draft]{todonotes}   

\usepackage{cases}
\usepackage{xcolor}

\DeclareMathOperator{\spn}{Span}

\newcommand{\dt}{\partial_t}

\newcommand{\dvh}{\mathrm{div}_h\,}
\newcommand{\nablah}{\nabla_h}
\newcommand{\deltah}{\Delta_h}
\newcommand{\dz}{\partial_z}

\newcommand{\idx}{\,d\vec{x}}
\newcommand{\idxh}{\,dxdy}
\newcommand{\subeqref}[2]{$\eqref{#1}_{#2}$}
\newcommand{\abs}[2]{\bigl| #1 \bigr|^{#2}}
\newcommand{\norm}[2]{\bigl\Arrowvert #1 \bigr\Arrowvert_{#2}}
\newcommand{\normh}[2]{\bigl| #1 \bigr|_{#2}}
\newcommand{\spacetime}{(\Omega\times(0,T))}

\newtheorem{df}{Definition}
\newtheorem{thm}{Theorem}[section]
\newtheorem{lm}{Lemma}

\newtheorem{prop}{Proposition}
\newenvironment{pf}{\paragraph{Proof}}{\par\hfill$\square$\par}

\numberwithin{equation}{section}
\allowdisplaybreaks
\title{Global Existence of Weak Solutions to the Compressible Primitive Equations of Atmospheric Dynamics with Degenerate Viscosities}
\date{12 August 2018}

\author{Xin Liu\footnote{Department of Mathematics, Texas A{\&}M University, College Station, Texas 77843, USA. Email: stleonliu@gmail.com} \,\, and \,  Edriss S. Titi\footnote{Department of Mathematics, Texas A{\&}M University, College Station, Texas 77843. USA. Also: Department of Computer Science and Applied Mathematics, The Weizmann Institute of Science, Rehovot 76100, Israel. Email: titi@math.tamu.edu \, and \, edriss.titi@weizmann.ac.il}}

\newcommand\blfootnote[1]{%
  \begingroup
  \renewcommand\thefootnote{}\footnote{#1}%
  \addtocounter{footnote}{-1}%
  \endgroup
}

\begin{document}
	\allowdisplaybreaks
	\maketitle	
	\blfootnote{\textup{2010} \textit{Mathematics Subject Classification}: \textup{35Q30}, \textup{35Q35}, \textup{76N10}}

\begin{abstract}
	We show the existence of global weak solutions to the three-dimensional compressible primitive equations of atmospheric dynamics with degenerate viscosities. In analogy with the case of the compressible Navier-Stokes equations, the weak solutions satisfy the basic energy inequality, the Bresh-Desjardins entropy inequality and the Mellet-Vasseur estimate. These estimates play an important role in establishing the compactness of  the vertical velocity of the approximating solutions, and therefore are essential to recover the vertical velocity in the weak solutions.
\end{abstract}

\tableofcontents

\section{Introduction}
\subsection{The compressible primitive equations}
In this work, we aim at studying the following compressible primitive equations of atmospheric dynamics without buoyancy (no gravity).
\begin{equation*}\tag{CPE} \label{CPE}
	\begin{cases}
		\dt \rho + \dvh (\rho v) + \partial_z (\rho w) = 0 & \text{in} ~ \Omega , \\
		\dt (\rho v) + \dvh( \rho v\otimes v) + \partial_z (\rho w v) + \nablah P(\rho) \\
		~~~~ ~~~~ = \dvh ( h(\rho) \mathcal D(v)) + \partial_{z} ( h(\rho) \partial_z v ) + \nablah ( g(\rho) \dvh v)& \text{in} ~ \Omega , \\
		\partial_z P(\rho) = 0& \text{in} ~ \Omega ,
	\end{cases}
\end{equation*}
where $ P(\rho) = \rho^\gamma $, with $ \gamma > 1 $.
Here, $ \rho $ is a scaler function denoting the density profile, and $ u = (v, w) \in \mathbb{R}^2 \times \mathbb{R} $ denotes the velocity filed. The domain is assumed to be $ \Omega := \Omega_h \times (0,1) $, with $ \Omega_h = \mathbb{T}^2 $,  being the fundamental periodic domain in $ \mathbb{R}^2 $.
$ \nablah= (\partial_x, \partial_y)^\top $ and $ \dvh= \nablah \cdot $ denote the horizontal gradient and divergence, respectively. Also, $ \deltah = \partial_{xx} + \partial_{yy} $ is the Laplacian operator in the horizontal direction.
The viscosity coefficients $ h, g $ are functions of the density. The horizontal viscosity tensor $ \mathcal D(v) $ is assumed to be symmetric with respect to the gradient, that is
\begin{equation*}
	\mathcal{D}(v) := \dfrac{1}{2} (\nablah v + (\nablah v)^\top).
\end{equation*}
%
The following boundary condition is imposed on system \eqref{CPE},
\begin{equation}\label{CPE:bc}
(\dz v, w) \bigr|_{z=0,1} = (0,0).
\end{equation}
The isentropic compressible primitive equations \eqref{CPE} are obtained by replacing the equation of the vertical velocity (vertical momentum) in the isentropic compressible Navier-Stokes equations with the hydrostatic balance equation, i.e., \subeqref{CPE}{3}. The compressible primitive equations are used to model the motion of flows in the atmosphere, which is an extremely complicated dynamical system. Indeed, as described in \cite{Lions1996}, the full hydrodynamic and thermodynamic equations with Coriolis force and gravity describe the motion and state of the atmosphere. However, such equations are too complicated to study and compute. By making use of the fact that the planetary horizontal scales are significantly larger than the vertical one in the atmosphere, one can replace the equation of the vertical momentum with the hydrostatic balance equation. Indeed, since the vertical scale of the atmosphere is significantly smaller than the planetary horizontal  scales, the authors in \cite{Ersoy2011a} take advantage, as it is commonly done in planetary scale geophysical models,  of the smallness of this aspect ratio  between these two orthogonal directions to formally   derive the compressible primitive equations of planetary atmospheric dynamics from the compressible Navier-Stokes equations. It turns out that the hydrostatic approximation balance equation is  accurate enough for practical applications and  is used as a fundamental equation in  atmospheric science models. It is the starting point of many large scale models in the theoretical investigations and practical weather and climate predictions (see, e.g., \cite{Lions1992}). This has also been observed by meteorologists (see, e.g., \cite{Richardson1965,Washington2005}). In fact, such an approximation is reliable and useful in the sense that the balance of gravity and pressure dominates the dynamic in the vertical direction and that the vertical velocity is usually hard to observe in the atmosphere. On the other hand, instead of the molecular viscosity,  eddy viscosity is used to model the statistical effect of turbulent motion in the atmosphere. The above observations and more perceptions from the meteorological point of view can be found in \cite[Chapter 4]{Richardson1965}.

In the simplified model \eqref{CPE}, we assume that the atmosphere is under adiabatic process. Moreover, we neglect here the Coriolis force and the gravity. However, it is worth stressing that adding the Coriolis force will not cause any difficulties and our results will be equally valid.
In addition, for the compressible Navier-Stokes equations, the kinetic theory implies that the viscosity coefficients $ h, g $ depend on the temperature, and for an isentropic flow, this dependence is translated into the dependence of the viscosity on the density (see, e.g., \cite{Liu1998a}).
Moreover, the relation 
\begin{equation}\label{19July2018}
	g(\rho) = h'(\rho) \rho - h (\rho).
\end{equation}
of the viscosity coefficients is introduced in \cite{Bresch2006} by Bresch and Desjardins for the compressible Navier-Stokes equations, which will induce a mathematical entropy, i.e., the Bresch-Desjardins entropy (see, also, \cite{Bresch2007}).
This motivates us to assume 
$ h (\rho) = \rho, g(\rho) = 0 $.

Consequently, system \eqref{CPE} can be written as,
\begin{equation*}\tag{CPE'}\label{CPE'}
	\begin{cases}
	\dt \rho + \dvh (\rho v) + \partial_z (\rho w) = 0& \text{in} ~ \Omega , \\
	\dt (\rho v) + \dvh( \rho v\otimes v) + \partial_z (\rho w v) + \nablah P(\rho) \\
	~~~~ ~~~~ = \dvh (\rho \mathcal D(v)) + \partial_{z} ( \rho \partial_z v )& \text{in} ~ \Omega,  \\
		\dz P(\rho) = 0 & \text{in} ~ \Omega,
	\end{cases}
\end{equation*}
with initial data
\begin{equation}\label{CPE-ini} (\rho, v)|_{t=0} = (\rho_0, v_0), \end{equation}
and boundary conditions \eqref{CPE:bc}.
Our goal is to establish the existence of global weak solutions to \eqref{CPE'}. Notice, in the vacuum set $ \lbrace \rho = 0 \rbrace $, the viscosity vanishes in system \eqref{CPE'}, and hence the system has a degenerate viscosity.

The first mathematical treatment of the compressible primitive equations (referred to as CPE) reformulated the system in the pressure coordinates ($p$-coordinates) and showed that in the new coordinate system, the equations are in the form of classical primitive equations (called the primitive equations, or PE hereafter) with the incompressibility condition in the new coordinates. This was done by Lions, Temam and Wang in \cite{Lions1992}. The authors also showed the existence of global weak solutions to the PE and therefore indirectly studied the CPE. We point out that the pressures on the upper and lower boundaries were assumed to be constant in \cite{Lions1992}, which is not the case in reality. In yet another work \cite{JLLions1992}, the authors modeled the nearly incompressible ocean by the PE. It is formulated as the hydrostatic approximation of the Boussinesq equations. The authors show the existence of global weak solutions and therefore indirectly study the CPE (see, e.g., \cite{JLLions1994,Lions2000} for additional work  by the authors).
Since then, the PE have been the subject of intensive mathematical research.  For instance, Guill{\'e}n-Gonz{\'a}lez, Masmoudi and Rodr{\'\i}guez-Bellido in \cite{GuillenGonzalez2001} study the local existence of strong solutions and global existence of strong solutions for small data to the PE. In \cite{HuTemamZiane2003} the authors address the global existence of strong solutions to the PE in a domain with small depth for large, but restricted initial data that depends on the depth. In \cite{Petcu2005}, the authors study the Sobolev and Gevrey regularity of the solutions to the PE. The first breakthrough concerning the global well-posedness of the PE is obtained by Cao and Titi in \cite{Cao2007}, in which the authors show the existence of unique global strong solutions (see, also, \cite{Cao2003,Kobelkov2006,Kukavica2007a,Kukavica2007,Kukavica2014,Zelati2015,Hieber2016,Li2017a,Ignatova2012,LiTiti2015global,LiTiti2016moisture,LiTiti2016book,hittmeir2017} and the references therein for related study). On the other hand, with partial anisotropic diffusion and viscosity, Cao, Li and Titi in \cite{Cao2012,Cao2014b,Cao2014,Cao2017,Cao2016,Cao2016a} establish the global well-posedness of strong solutions to PE. For the inviscid primitive equations, or hydrostatic incompressible Euler equations, in \cite{Brenier1999,Wong2012,Kukavica2011}, the authors show the existence of solutions in the analytic function space and in $H^s $ space. More recently, the authors in \cite{Wong2014,Cao2015} construct finite-time blowup for the inviscid PE in the absence of rotation. Also, in \cite{GerardVaret2018} , the authors establish the Gevrey regularity of hydrostatic Navier-Stokes equations with only vertical viscosity.

In this work, we focus on showing the global existence of weak solutions to \eqref{CPE'}. It should be emphasized that prior to this work, even though it has been studied indirectly, the global existence of weak solutions to the CPE has been open. This is due to the fact that the former study has neglected the effect of the physical boundary condition on the pressure. In particular, and as we stated above, in \cite{Lions1992} the pressures of the CPE on the boundaries are taken to be constant. In fact, the physical pressures on the boundaries should be changing in time cooperating for the evolving flows in the atmosphere. In fact, this is the very reason and motivation of our study of the existence of  weak solution to \eqref{CPE'}, which is our main goal in this work.

Recently, Gatapov, Kazhikhov, Ersoy, Ngom construct a global weak solution to some two dimensional compressible primitive equations in \cite{Ersoy2012,Gatapov2005}. Meanwhile, Ersoy, Ngom, Sy, Tang, Gao study the stability of weak solutions to CPE in \cite{Ersoy2011a,Tang2015}, in the sense that a sequence of weak solutions satisfying some entropy conditions contains a subsequence converging to another weak solution. In this work, we address the problem of existence of such global weak solutions to \eqref{CPE'}.

The existence of global weak solutions to the isentropic compressible Navier-Stokes equations with constant viscosity coefficients was established by P.-L.Lions in \cite{Lions1998}. It was further developed by Feireisl, Novotn\'y and Petzeltov\'a in \cite{Feireisl2001} (see also \cite{Feireisl2004}). Very recently, independently by Vasseur and Yu in \cite{Vasseur2016}, and by Li and Xin in \cite{Li2015a}, the problem of existence of global weak solutions to the compressible Navier-Stokes equations with degenerate viscosity is addressed.

In order to establish the existence of global weak solutions to \eqref{CPE'}, our main difficulty is due to the lack of evolutionary equation for the vertical velocity (momentum). Indeed, the vertical velocity is recovered from the density and the horizontal velocity by making use of the stratifying property of the density (induced by the hydrostatic balance equation \subeqref{CPE'}{3}) and the continuity equation \subeqref{CPE'}{1} (see, e.g., \eqref{GR-011}, below). Furthermore, the vertical velocity is less regular than the horizontal velocity in general. To overcome this  difficulty, in \cite{Ersoy2011a,Tang2015}, the authors make use of the Bresch-Desjardins entropy inequality by taking the viscosity coefficients as in \eqref{CPE'}. This will give some extra regularity to the density and the vertical velocity. That is, one will have $ \rho^{1/2} \in L^\infty(0,T;H^1(\Omega)) $ and $ \rho^{1/2}\dz w \in L^2(\Omega\times(0,T)) $ (see \eqref{BD-007} and \eqref{ene:BD-entropy}, below). Moreover, one has, as in the compressible Navier-Stokes equations, the Mellet-Vasseur estimate for the horizontal velocity (see \eqref{ene:MV-estimate}, below). Such estimates will help us to deal with the required compactness argument in the nonlinearities involving the vertical velocity.

We follow here the same strategy of \cite{Li2015a}. Indeed, we introduce the approximating problem in section \ref{sec:approximation-problem}. In particular, such an approximation captures the bound of the basic energy, the Bresch-Desjardins entropy inequality and the Mellet-Vasseur estimate (see Proposition \ref{prop:uniform-est}, below). The first main part of this work is to construct global strong solutions (the approximating solutions) to the approximating problem. This is done in section \ref{sec:aprx-sol-exist}. In fact, our approximating solutions are less regular than those for the compressible Navier-Stokes equations in \cite{Li2015a} (see section \ref{sec:aprx-propri-est}). We employ a Galerkin type approximating scheme, modified from that in \cite{Feireisl2001},  to show the existence of global strong solutions to the approximating problem \eqref{eq:approximating-CPE}, below (see section \ref{sec:aprx-exist}). Then we show the compactness result which will be the final ingredient in this work. We emphasize that the main convergences to restore the compactness in the nonlinearities, involving the vertical velocity, are the strong convergences in \eqref{CM-density}, \eqref{CM-energy} and the weak convergence in \eqref{CV-002}, below. In particular, the vertical velocity is obtained in \eqref{CV-def-vertical-w}, below.

In this work, we will assume the initial data $ (\rho,v)\bigr|_{t=0} = (\rho_0, v_0) $ satisfies the following bounds. For $ \gamma > 1 $ and some $ \varpi > 0 $,
\begin{equation}\label{initial-data}
	\begin{gathered}
		\rho_0 = \rho_0(x,y) \geq 0, ~ \text{almost everywhere in} ~ \Omega, ~ \rho_0 \in L^1(\Omega) \cap L^\gamma(\Omega),\\
		 \nabla_h \rho_0^{1/2} \in L^2(\Omega), ~
		m_0 := \rho_0 v_0 \in L^{2\gamma/(\gamma+1)}(\Omega), \\
		m_0=0, ~ \text{almost everywhere in the set} ~ \lbrace \rho_0 = 0 \rbrace \subset \Omega, \\
		\rho_{0}^{-1- \varpi}\abs{m_0}{2+\varpi} \in L^1(\Omega), ~
		\text{and} ~ \rho_0 v_0^2 = \rho_0^{-1} m_0^2 \in L^1(\Omega),
	\end{gathered}
\end{equation}
where we set $ \rho_0^{-1} m_0^2 = \rho_{0}^{-1- \varpi}\abs{m_0}{2+\varpi} = 0 $, almost everywhere in the set $ \lbrace \rho_0 = 0 \rbrace $, i.e., in the vacuum set of $ \rho_0 $.

This work is organized as follows. In the next section, we will introduce some notations, that we will be using in this work, and we will also state our main theorem. In section \ref{sec:approximation-problem}, we will introduce our approximating problem of \eqref{CPE'}. With the approximating solutions defined there, we establish the uniform estimates in section \ref{sec:uniform-estimates}, which will capture the basic energy estimate, the Bresch-Desjardins entropy and the Mellet-Vasseur estimates. Next, in section \ref{sec:aprx-sol-exist}, we focus on the construction of the approximating solutions. This is done by presenting the a priori estimates in section \ref{sec:aprx-propri-est} and the constructing the solutions via a modified Galerkin type approximating scheme in section \ref{sec:aprx-exist}. Finally, in section \ref{sec:compactness}, we summarize the compactness argument and show the existence of weak solutions to \eqref{CPE'}. Moreover, we show the energy inequality and the entropy inequalities. This in turn will complete the proof of our main theorem.

\subsection{Preliminaries and the main theorem}
We will use the notations
\begin{equation*}
	\int \cdot \idx = \int_{\Omega} \cdot \idx = \int_{\Omega} \cdot \,dxdydz, ~ \int_{\Omega_h} \cdot \idxh
\end{equation*}
to denote the integrals in the whole domain and in the horizontal domain, respectively. Furthermore, we use
\begin{equation*}
	\norm{\cdot}{L^p} = \norm{\cdot}{L^p(\Omega)}, ~  \normh{\cdot}{L^p} = \norm{\cdot}{L^p(\Omega_h)},
\end{equation*}
to denote the $ L^p $, $ p > 1 $, norms in the whole spatial domain and in the horizontal domain, respectively. $ W^{k,p}(\Omega), W^{k,p}(\Omega_h), k \in \mathbb N^+$, $p > 1 $,  are the corresponding Sobolev space in the whole domain and in the horizontal domain, respectively, and when $ p = 2 $, they are denoted by $ H^k(\Omega), H^k (\Omega_h) $. $ W^{-1,1}(\Omega) = \bigl(W^{1,\infty}(\Omega)\bigr)^* $ is the dual space of $ W^{1,\infty}(\Omega) $. $ C_c^\infty (\overline\Omega\times[0,T)) $ consists of the smooth function with compact support in $ \overline{\Omega}\times[0,T) $ for any given positive constant $ T \in (0,\infty) $ (e.g., trigonometric polynomials of variables $ (x,y) $ such that the coefficients are $ C^\infty $ functions of $ (z, t) $ with support in any compact subset of $ [0,1] \times [0,T) $).
$ \mathcal D'(0,T) $ is the space of distribution functions in the temporal variable. Also, we have the following  De Giorgi-type lemma.
\begin{lm}[De Giorgi]\label{lm:De-Giorgi}
	Consider a function on the positive real line $$ g: \mathbb R^+ \mapsto [0,\infty). $$ Suppose that $ g $ is bounded and non-increasing, and that there exist constants $ \alpha , \beta , C  \in (0, \infty) $ such that, for any $ l, k \in \mathbb R^+ $, with $ l > k $,
	\begin{equation}\label{lm:dg-001}
		g(l) \leq C (l-k)^{-\beta} g(k)^{1+\alpha}.
	\end{equation}
	Then there exists a constant $ L \in (0,\infty) $ such that $ g \equiv 0 $ on $ [L,\infty) $.
\end{lm}
\begin{pf}
	Take $ k = 1 $ in \eqref{lm:dg-001}. It follows that there is a constant $ C' $ such that for any $ l > 1 $,
	\begin{equation}\label{lm:dg-002}
		g(l) \leq C(l-1)^{-\beta}g(1)^{1+\alpha} \leq \dfrac{C'}{l^{\beta}}.
	\end{equation}
	Define the monotonically increasing sequence $ \lbrace h_i \rbrace_{i=0,1,2,\cdots} $, by $ h_i =  h_{i-1} + i^{-2} $, for $ i \geq 1 $, and $ h_0 = \kappa $, with $ \kappa\in (1,\infty) $ to be determined later. Then from \eqref{lm:dg-001}, for $ i \geq 2 $,
	\begin{align*}
		& g(h_i) \leq C i^{2\beta} g(h_{i-1})^{1+\alpha}, ~~~~ \text{or equivalently} ~~\\
		& i^{\frac{4\beta}{\alpha}} g(h_i) \leq C \biggl( \dfrac{i}{i-1} \biggr)^{\frac{2\alpha\beta + 4\beta}{\alpha}} g(h_{i-1})^{\alpha/2} \bigl( (i-1)^{\frac{4\beta}{\alpha}} g(h_{i-1}) \bigr)^{1+\alpha/2} \\
		& ~~~~ \leq C 2^{\frac{2\alpha\beta + 4\beta}{\alpha}} g(h_0)^{\alpha/2} \bigl( (i-1)^{\frac{4\beta}{\alpha}} g(h_{i-1}) \bigr)^{1+\alpha/2} \\
		& ~~~~ \leq C(C')^{\alpha/2} 2^{\frac{2\alpha\beta + 4\beta}{\alpha}} \kappa^{-\frac{\alpha\beta}{2}} \bigl( (i-1)^{\frac{4\beta}{\alpha}} g(h_{i-1}) \bigr)^{1+\alpha/2},
	\end{align*}	
	where we have applied inequality \eqref{lm:dg-002} for $ l = h_0 $. Then for $ \kappa $ large enough, we have for $ i \geq 2 $,
	\begin{align*}
		i^{\frac{4\beta}{\alpha}} g(h_i) \leq \bigl( (i-1)^{\frac{4\beta}{\alpha}} g(h_{i-1}) \bigr)^{1+\alpha/2}. 
	\end{align*}
	We iterate this inequality $ i $ times to obtain
	\begin{equation*}
		i^{\frac{4\beta}{\alpha}} g(h_i) \leq \bigl(g(h_0)\bigr)^{(1+\alpha/2)i} = \bigl(\dfrac{C'}{\kappa}\bigr)^{(1+\alpha/2)i} \leq 1,
	\end{equation*}
	for $ \kappa $ large enough, where we have applied inequality \eqref{lm:dg-002} for $ l = h_0 = \kappa $. Thus
	\begin{equation*}
		g(h_i) \leq i^{-\frac{4\beta}{\alpha}}.
	\end{equation*}
	Since $ h_i $ is increasing, and $ g $ is non-increasing and non-negative, we have
	\begin{equation*}
		\lim\limits_{i\rightarrow \infty} g(h_i) = g(h_\infty) = 0,
	\end{equation*}
	where $ h_\infty = \kappa + \sum_{i=1}^\infty i^{-2} \in (0,\infty) $.
	Consequently, since $ g $ is non-increasing and non-negative, then  $ g \equiv 0 $ on $ [L,\infty) $, for any $ L \geq h_\infty $.
\end{pf}
The following is the definition of weak solution to CPE in this work:
\begin{df}[Global weak solutions to CPE]\label{def:weak-sol-CPE}
	Given any positive constant time $ T \in (0,\infty) $, a triple $ (\rho, v, w) $ is called a weak solution to \eqref{CPE'}, with \eqref{CPE:bc} and \eqref{initial-data}, in $ \Omega \times (0,T) $, if $\rho  = \rho(x,y,t) $ and it satisfies the following weak formulation,
	\begin{gather*}
		\int_\Omega \rho_0 \psi|_{t=0} \idx + \int_0^T \int_\Omega \biggl\lbrack \rho \dt \psi + \rho v \cdot\nablah \psi + \rho w \dz \psi \biggr\rbrack \idx\,dt = 0,\\
		\int_\Omega m_0 \cdot \phi|_{t=0} \idx + \int_0^T \int_\Omega  \biggr\lbrack \rho v \cdot \dt \phi + \rho v \otimes v : \nablah \phi + \rho w v \cdot \dz \phi \\
		 + \rho^\gamma \dvh \phi \biggr\rbrack  \idx\,dt - \int_0^T \int_\Omega \biggl\lbrack \rho \mathcal D(v) : \nablah \phi + \rho \dz v \cdot \dz \phi \biggl\rbrack \idx\,dt =0,
	\end{gather*}
	for any smooth functions $ (\psi, \phi) \in C_c^{\infty}(\overline\Omega\times [0,T);\mathbb R) \times C_c^{\infty}(\overline\Omega\times [0,T);\mathbb R^2) $, where the following regularity holds
	\begin{gather*}
		0 \leq \rho = \rho(x,y,t) \in L^\infty(0,T;L^1(\Omega)\cap L^\gamma(\Omega)), \\
		\nablah \rho^{1/2} \in L^\infty(0,T;L^2(\Omega)), \rho^{1/2} v \in L^\infty(0,T;L^2(\Omega)), \\
		\nablah \rho^{\gamma/2} \in L^2(\Omega\times(0,T)), ~ \rho^{1/2} \nabla v \in L^2(\Omega\times(0,T)), \\
		 \rho^{1/2} \dz w \in L^2(\Omega\times(0,T)).
	\end{gather*}
	In addition, the following energy inequality holds,
	\begin{align*}
		& \dfrac{d}{dt} \bigl\lbrace \dfrac 1 2 \int \rho \abs{v}{2} \idx + \dfrac{1}{\gamma - 1} \int \rho^\gamma \idx \bigr\rbrace + \int \rho \abs{\mathcal{D}(v)}{2} \idx \\
		& ~~~~ + \int \rho \abs{\dz v}{2} \idx \leq 0,
	\end{align*}
in $ \mathcal D'(0,T) $. Also, we have the following entropy inequality,
	\begin{align*}
		& \sup_{0\leq t\leq T} \bigl\lbrace \int \abs{\nablah \rho^{1/2}}{2} \idx + \int \rho (e + v^2) \log(e+v^2) \idx \bigr\rbrace \\
		& ~~~~  ~~~~ + \int_0^T \int_\Omega \rho \abs{\nablah v}{2} + \rho \abs{\dz w}{2} \idx\,dt
		 + \int_0^T \int_\Omega \rho^{\gamma-2} \abs{\nablah \rho}{2} \idx\,dt \\
		 & ~~~~ \leq C_T \mathfrak E_0 ,
	\end{align*}
	where $ \mathfrak E_0 $ is defined as
	\begin{equation}\label{bound-initial-data}
		\begin{aligned}
		& \mathfrak E_0 : =  \norm{\rho_0^{-1}m_0^2}{L^1(\Omega)} + \norm{\rho_0^{-1-\varpi}\abs{m_0}{2+\varpi}}{L^1(\Omega)} + \norm{\nablah \rho_0^{1/2}}{L^2(\Omega)} \\
		 & ~~~~ ~~~~ + \norm{\rho_0}{L^1(\Omega)} + \norm{\rho_0}{L^\gamma(\Omega)}^\gamma.
		 \end{aligned}
	\end{equation}
\end{df}
Next, we introduce the main theorem in this work.
\begin{thm}\label{main-theorem}
	Considering the initial data given in \eqref{initial-data}, then there is a global weak solutions to the problem \eqref{CPE'}, \eqref{CPE:bc}  \eqref{CPE-ini}, satisfying Definition  \ref{def:weak-sol-CPE}.
\end{thm}

\section{The approximating problem}\label{sec:approximation-problem}
In this section, we will study the following approximating problem of \eqref{CPE'}. Consider the following parabolic regularizing system,
\begin{equation}\label{eq:approximating-CPE}
	\begin{cases}
		\dt \rho_\varepsilon + \dvh (\rho_\varepsilon v_\varepsilon) + \partial_z (\rho_\varepsilon w_\varepsilon) = G(\rho_\varepsilon)& \text{in} ~ \Omega , \\
		\rho_\varepsilon \dt v_\varepsilon +  \rho_\varepsilon v_\varepsilon \cdot\nablah v_\varepsilon + \rho_\varepsilon w_\varepsilon \dz v_\varepsilon + \nablah P(\rho_\varepsilon) \\
		~~~~ ~~~~ = \dvh (\rho_\varepsilon \mathcal D(v_\varepsilon)) + \partial_{z} ( \rho_\varepsilon \partial_z v_\varepsilon ) + F(v_\varepsilon,\rho_\varepsilon)& \text{in} ~ \Omega , \\
		\partial_z P(\rho_\varepsilon) = 0& \text{in} ~ \Omega ,
	\end{cases}
\end{equation}
where
\begin{align*}
	G(\rho_\varepsilon) & : = \varepsilon \rho_\varepsilon^{1/2} \deltah \rho_\varepsilon^{1/2} + \varepsilon \rho_\varepsilon^{1/2} \dvh ( \abs{\nablah \rho_\varepsilon^{1/2}}{2} \nablah \rho_\varepsilon^{1/2} ) + \varepsilon \rho_\varepsilon^{-p_0} ,\\
	F(v_\varepsilon,\rho_\varepsilon) & : =  \sqrt{\varepsilon} \dvh(\rho_\varepsilon \nablah v_\varepsilon) + \widetilde F (v_\varepsilon,\rho_\varepsilon) ,\\
	\widetilde F (v_\varepsilon,\rho_\varepsilon) & := \varepsilon \rho_\varepsilon^{1/2} \abs{\nablah \rho_\varepsilon^{1/2} }{2} \nablah \rho_\varepsilon^{1/2} \cdot \nablah v_\varepsilon - \varepsilon \rho_\varepsilon^{-p_0} v_\varepsilon - \varepsilon \rho_\varepsilon \abs{v_\varepsilon}{3}v_\varepsilon,
\end{align*}
with the positive constants $ \varepsilon, p_0 $ satisfying,
\begin{align*}
	\varepsilon > 0 , p_0 > \max\lbrace 24 , \gamma - 1 \rbrace > 0.
\end{align*}
Similarly, we complement the above system with the following initial and boundary conditions,
\begin{equation}\label{approx:ini-bc}
	(\rho_\varepsilon,v_{\varepsilon})|_{t=0} = (\rho_{\varepsilon,0}, v_{\varepsilon,0}), ~ (\dz v, w)|_{z=0,1} = (0,0).
\end{equation}
Next, we will define the approximating initial data and solutions of the approximating compressible primitive equation \eqref{eq:approximating-CPE}.
\begin{df}[Approximating initial data]\label{def:aprxm-initial}
Let $ \varepsilon_0 \in (0,1) $ be a finite positive constant. A sequence of approximating initial data to $ (\rho_0, v_0) $, with $ m_0 = \rho_0 v_0 $, as in \eqref{initial-data}, is $ \lbrace (\rho_{\varepsilon,0}, v_{\varepsilon,0}) \rbrace_{\varepsilon\in (0,\varepsilon_0)} $ satisfying the following conditions:
\begin{itemize}
	\item As $ \varepsilon \rightarrow 0^+ $,
	\begin{gather*}
		\norm{\rho_{\varepsilon,0} - \rho_0}{L^1(\Omega)} + \norm{\rho_{\varepsilon,0} - \rho_0}{L^\gamma(\Omega)} + \norm{\nablah(\rho^{1/2}_{\varepsilon,0} - \rho_0^{1/2})}{L^2(\Omega)} \\
		 + \norm{\rho_{\varepsilon,0} v_{\varepsilon,0}^2 - \rho_0^{-1} m_0^2 }{L^1(\Omega)} 
		 \rightarrow 0;
	\end{gather*}
	
	\item for any fixed $ \varepsilon \in (0,\varepsilon_0) $,
\begin{gather*}
	\varepsilon^{1/p_0+1} < \rho_{\varepsilon,0}< \varepsilon^{-1/p_0-1}, ~ \dz \rho_{\varepsilon,0} = 0, ~ \dz v_{\varepsilon,0}|_{z=0,1} = 0, \\
	\rho_{\varepsilon,0} \in W^{1,4}(\Omega), ~
	v_{\varepsilon,0} \in H^1(\Omega),
\end{gather*}
with
\begin{align*}
	& \norm{\rho_{\varepsilon,0}^{1/2} v_{\varepsilon,0}}{L^2(\Omega)}^2 + \norm{\rho_{\varepsilon,0}}{L^\gamma(\Omega)}^\gamma + \norm{\rho_{\varepsilon,0}}{L^1(\Omega)} + \norm{\nablah \rho^{1/2}_{\varepsilon,0}}{L^2(\Omega)}^2 \\
	& ~~~~~ + \norm{\rho_{\varepsilon,0} v_{\varepsilon,0}^{2+\varpi}}{L^1(\Omega)}+ \norm{\varepsilon^{1/p_0} \rho_{\varepsilon,0}^{-1}}{L^{p_0}(\Omega)}^{p_0} + \norm{\varepsilon^{1/4} \nablah \rho^{1/2}_{\varepsilon,0}}{L^4(\Omega)}^4 \\
	  & < 2 \norm{\rho_0^{-1}m_0^2}{L^1(\Omega)}+ 2 \norm{\rho_0}{L^\gamma(\Omega)}^\gamma + 2 \norm{\rho_0}{L^1(\Omega)} + 2 \norm{\nablah \rho_0^{1/2}}{L^2(\Omega)}  \\
	  & ~~~~ + 2 \norm{\rho_0^{-1-\varpi}\abs{m_0}{2+\varpi}}{L^1(\Omega)} =: C_0.
\end{align*}
\end{itemize}
\end{df}

\begin{df}[Approximating solutions]\label{def:aprxm-sols}
	Let $ \varepsilon_0 \in (0, 1) $ be a finite positive constant.
	Then a sequence of approximating solutions $ \lbrace (\rho_\varepsilon, v_\varepsilon) \rbrace_{\varepsilon \in (0,\varepsilon_0) } $ consists of the global strong solutions $(\rho_\varepsilon, v_\varepsilon) $  to the approximating equations \eqref{eq:approximating-CPE} with the approximating initial data $ (\rho_{\varepsilon, 0}, v_{\varepsilon, 0})$, as defined in Definition \ref{def:aprxm-initial}, for all $ \varepsilon \in (0,\varepsilon_0) $. In particular, for any positive time $ T \in (0, \infty) $, there is a positive constant $ C= C_{\varepsilon,T} < \infty $ such that $ C^{-1} < \rho_\varepsilon < C $ and,
	\begin{gather*}
		\rho_{\varepsilon}^{1/2} \in L^\infty(0,T;W^{1,4}(\Omega)) \cap L^2(0,T;H^2(\Omega)), \nablah \rho_\varepsilon \in L^{30}(\Omega\times(0,T)), \\
		 \abs{\nablah \rho_\varepsilon^{1/2}}{2} \nablah^2 \rho_\varepsilon^{1/2} \in L^2(\Omega\times(0,T)), \dt \rho_\varepsilon\in L^2(\Omega\times(0,T)), \\
		v_\varepsilon \in L^\infty(0,T;H^1(\Omega))\cap L^2(0,T;H^2(\Omega)) \cap L^5(\Omega\times(0,T)) \cap L^{10}(\Omega\times(0,T)), \\
		\nabla v_\varepsilon \in L^{10/3}(\Omega\times(0,T)), \dt v_\varepsilon\in L^2(\Omega\times(0,T)), \dz w_\varepsilon \in L^2(\Omega\times(0,T)) .
	\end{gather*}
	Here $ w_\varepsilon $ is given by
	\begin{align*}
		& w_\varepsilon(x,y,z,t) = - \rho_\varepsilon^{-1}(x,y,t) \int_0^z \dvh \bigl(\rho_\varepsilon(x,y,t) ( v_\varepsilon(x,y,z',t) \\
		& ~~~~ ~~~~ ~~~~ ~~~~ - \int_0^1 v_\varepsilon(x,y,z'',t) \,dz'' )\bigr)\,dz'.
	\end{align*}
\end{df}

We will establish the global existence of approximating solutions for $ \varepsilon \in (0,\varepsilon_0) $, for some $ \varepsilon_0 $ small enough. In addition, we will establish some estimates for the approximating solutions, that are independent of $ \varepsilon $, capturing the basic energy, the Bresch-Desjardins entropy and the Mellet-Vasseur estimates of \eqref{CPE'}. In fact, we will show the following:

\begin{prop}[Global existence of approximating solutions]\label{prop:aprx-sol-exist}
	For any fixed $ \varepsilon > 0 $, there is an approximating solution $ (\rho_\varepsilon, v_\varepsilon) $ as defined in Definition \ref{def:aprxm-sols}.
\end{prop}
	The proof will be presented in section \ref{sec:aprx-sol-exist}. In particular, see the arguments in section \ref{sec:aprx-exist}.

%

\section{$ \varepsilon $-independent uniform estimates}\label{sec:uniform-estimates}
In this section we will establish some uniform a priori estimates of the approximating solutions defined in Definition \ref{def:aprxm-sols}
independent of $ \varepsilon $. These estimates will capture the basic energy, the Bresch-Desjardins entropy and the Mellet-Vasseur estimates of the solutions to \eqref{CPE'}. In fact, we are going to show the following:
\begin{prop}\label{prop:uniform-est}
There is a constant $ \varepsilon_0 > 0 $, small enough, such that for every $ T \in (0,\infty) $ the approximating solution $ (\rho_\varepsilon, v_\varepsilon) $ on $ [0,T) $, with $ \varepsilon \in (0,\varepsilon_0) $, as defined 	
	in Definition \ref{def:aprxm-sols} will satisfy the following estimates: for some constant $ C_T > 0 $, which is independent of $ \varepsilon $, but depends on $ T $ and the initial data \eqref{initial-data} and \eqref{bound-initial-data},
	\begin{itemize}
		\item the bound of the basic energy bound:
			\begin{align*}
				& \sup_{0\leq t\leq T} \biggl\lbrace \int \rho_\varepsilon \abs{v_\varepsilon}{2} \idx + \int \rho^\gamma_\varepsilon \idx + \int \rho_\varepsilon \idx + \varepsilon \int \rho_\varepsilon^{-p_0} \idx \biggr\rbrace \\
				& ~~~~ + \int_0^T \int \biggl( \rho_\varepsilon \abs{\mathcal D (v_\varepsilon)}{2} + \rho_\varepsilon \abs{\dz v_\varepsilon}{2} \biggr) \idx\,dt
				+ \varepsilon \int_0^T\int \biggl(\rho_\varepsilon^{-p_0} \abs{v_\varepsilon}{2} + \rho_\varepsilon \abs{v_\varepsilon}{5} \\
				& ~~~~ + \abs{\nablah \rho_\varepsilon^{1/2}}{4}\abs{v_\varepsilon}{2} +  \abs{\nablah \rho_\varepsilon^{1/2}}{2} \abs{v_\varepsilon}{2} + \abs{\nablah \rho_\varepsilon^{1/2}}{4} \biggr)\idx \,dt \\
				& ~~~~ + \varepsilon^2 \int_0^T \int \rho_\varepsilon^{-2p_0-1}\idx \,dt \leq C_T,
			\end{align*}
		\item the Bresch-Desjardins entropy inequality:
			\begin{align*}
				& \sup_{0\leq t\leq T} \biggl\lbrace \int  \abs{\nablah \rho_\varepsilon^{1/2}}{2} \idx  + \varepsilon \int  \abs{\nablah \rho_\varepsilon^{1/2}}{4}\idx \biggr\rbrace + \int_0^T \int \biggl( \rho_\varepsilon\abs{\nablah v_\varepsilon}{2} \\
				& ~~~~ + \rho_\varepsilon \abs{\dz w_\varepsilon}{2} \biggr) \idx\,dt + \int_0^T \int  \rho_\varepsilon^{\gamma-2}\abs{\nablah \rho_\varepsilon}{2} \idx \,dt  + \varepsilon \int_0^T \int \biggl( \abs{\deltah \rho_\varepsilon^{1/2}}{2} \\
				& ~~~~ + \abs{\nablah \rho_\varepsilon^{1/2}}{2} \abs{\nablah^2 \rho_\varepsilon^{1/2}}{2} \biggr) \idx \,dt + \varepsilon^2 \int_0^T \int  \abs{\nablah \rho_\varepsilon^{1/2}}{4} \abs{\nablah^2 \rho_\varepsilon^{1/2}}{2} \idx \,dt\\
				& ~~~~ \leq C_T,
			\end{align*}
		\item and the Mellet-Vasseur estimate:
			\begin{equation*}
				\sup_{0 \leq t\leq T} \int \rho_\varepsilon(e + v_\varepsilon^2) \log(e+ v_\varepsilon^2) \idx \leq C_T,
			\end{equation*}
	\end{itemize}
	where $ w_\varepsilon $ is given by
	\begin{align*}
		& w_\varepsilon(x,y,z,t) = - \rho_\varepsilon^{-1}(x,y,t) \int_0^z \dvh \bigl(\rho_\varepsilon(x,y,t) ( v_\varepsilon(x,y,z',t) \\
		& ~~~~ ~~~~ ~~~~ ~~~~ - \int_0^1 v_\varepsilon(x,y,z'',t) \,dz'' )\bigr)\,dz'.
	\end{align*}
\end{prop}
\begin{pf} This follows directly from
	\eqref{ene:basic-energy},  \eqref{ene:BD-entropy} and \eqref{ene:MV-estimate} in the following subsections.
\end{pf}

 For the sake of convenience, we shall denote, in the following subsections, the approximating solution $ (\rho_\varepsilon, v_\varepsilon) $ as $ (\rho, v) $ and $ w_\varepsilon $ as $ w $. It should be noted that the following calculations are rigorous in view of the regularity given in Definition \ref{def:aprxm-sols}.

\subsection{The basic energy bound}\label{sec:basic-energy}
Multiplying \subeqref{eq:approximating-CPE}{2} with $ v $ and integrating the resultant over $ \Omega $ to obtain
\begin{equation}\label{BE-003}
	\begin{aligned}
		& \dfrac{d}{dt} \biggl\lbrace \dfrac 1 2 \int \rho \abs{v}{2} \idx \biggr\rbrace - \int \dfrac{1}{2} G(\rho) \abs{v}{2} \idx + \int \nablah P(\rho) \cdot v \idx \\
		& ~~~~  + \int \biggl( \rho \abs{\mathcal{D}(v)}{2} + \rho \abs{\dz v}{2} \biggr) \idx  + \sqrt{\varepsilon} \int \rho \abs{\nablah v}{2} \idx  = \int \widetilde F (v,\rho) \cdot v \idx.
	\end{aligned}
\end{equation}
In the meantime, by making use of \subeqref{eq:approximating-CPE}{1}, after integration by parts, the following identity holds, for any  $ q \neq 1 $,
\begin{equation}\label{BE-002}
	\begin{aligned}
		& \int \nablah \rho^{q} \cdot v\idx = \int \dfrac{q}{q - 1} \nablah \rho^{q - 1} \cdot (\rho v) \idx = \int \biggl( \dfrac{q}{q - 1} \rho^{q - 1}  \\
		& ~~~~ ~~~~ ~~~~ \times (\dt \rho - \dz(\rho w) - G(\rho)) \biggr) \idx = \dfrac{d}{dt} \bigl\lbrace \dfrac{1}{q - 1} \int \rho^q\idx \bigr\rbrace \\
		& ~~~~ + \dfrac{\varepsilon (2q-1) q}{16(q-1)} \int \rho^{q-3}(4 \rho + \abs{\nablah \rho}{2}) \abs{\nablah \rho}{2} \idx - \dfrac{\varepsilon q}{q-1} \int \rho^{q-p_0-1} \idx, 
	\end{aligned}
\end{equation}
where we have also used the fact that $ \rho $ is independent of the $ z $ variable and the following identity,
\begin{align*}
	& \int \rho^{q-1} G(\rho) \idx = \int \rho^{q-1} \bigl(\varepsilon \rho^{1/2} \deltah \rho^{1/2} + \varepsilon \rho^{1/2} \dvh ( \abs{\nablah \rho^{1/2}}{2} \nablah \rho^{1/2} ) \\
	& ~~~~ ~~~~ + \varepsilon \rho^{-p_0} \bigr)\idx  = - \varepsilon (2q-1) \int \biggl(  \rho^{q-1} \abs{\nablah \rho^{1/2}}{2} + \rho^{q-1} \abs{\nablah \rho^{1/2}}{4} \biggr) \idx \\
	& ~~~~ ~~~~  + \varepsilon \int \rho^{q-p_0-1} \idx.
\end{align*}
Therefore, applying \eqref{BE-002}, with $ q = \gamma $ in \eqref{BE-003}, one has
\begin{equation}\label{BE-001}
	\begin{aligned}
		& \dfrac{d}{dt} \biggl\lbrace \dfrac 1 2 \int \rho \abs{v}{2} \idx + \dfrac{1}{\gamma - 1} \int \rho^\gamma \idx \biggr\rbrace + \int \rho \abs{\mathcal{D}(v)}{2} \idx + \int \rho \abs{\dz v}{2} \idx \\
		& ~~~~  + (\sqrt{\varepsilon} - \varepsilon) \int \rho \abs{\nablah v}{2} \idx  + \dfrac{\varepsilon}{2} \int \rho^{-p_0} \abs{v}{2} \idx + \varepsilon \int \rho \abs{v}{5} \idx \\
		& ~~~~  + \dfrac{\varepsilon}{32} \int \rho^{-2} \abs{\nablah \rho}{4} \abs{v}{2} \idx  + \dfrac{\varepsilon}{16} \int \rho^{-1} \abs{\nablah \rho}{2} \abs{v}{2} \idx \\
		& ~~~~  + \dfrac{\varepsilon(2\gamma - 1)\gamma}{16(\gamma - 1)} \int \rho^{\gamma-3}(4\rho + \abs{\nablah \rho}{2})\abs{\nablah \rho}{2} \idx \leq \dfrac{\varepsilon \gamma}{\gamma - 1} \int \rho^{\gamma - p_0 - 1 }\idx,
	\end{aligned}
\end{equation}
where we have substituted above the following calculation,
\begin{align*}
	& \int \widetilde F (v,\rho) \cdot v\idx + \int \dfrac{1}{2} G(\rho) \abs{v}{2}\idx
	 = - \varepsilon \int \rho^{-p_0} \abs{v}{2} \idx - \varepsilon \int \rho \abs{v}{5} \idx \\
	& ~~~~ - \dfrac{\varepsilon}{2} \int \dvh (\rho^{1/2} \abs{\nablah \rho^{1/2}}{2} \nablah \rho^{1/2}) \abs{v}{2} \idx + \int \biggl( \dfrac{\varepsilon}{2} \rho^{1/2}\deltah \rho^{1/2} \abs{v}{2} \\
	& ~~~~ + \dfrac{\varepsilon}{2} \dvh(\rho^{1/2} \abs{\nablah \rho^{1/2}}{2} \nablah \rho^{1/2}) \abs{v}{2} \biggr) \idx  - \dfrac{\varepsilon}{2} \int \abs{\nablah \rho^{1/2}}{4} \abs{v}{2} \idx \\
	& ~~~~ + \dfrac{\varepsilon}{2} \int \rho^{-p_0} \abs{v}{2} \idx
	= - \dfrac{\varepsilon}{2} \int \rho^{-p_0} \abs{v}{2} \idx - \varepsilon \int \rho \abs{v}{5} \idx \\
	& ~~~~ - \dfrac{\varepsilon}{32} \int \rho^{-2} \abs{\nablah \rho}{4} \abs{v}{2} \idx - \dfrac{\varepsilon}{8} \int \rho^{-1} \abs{\nablah \rho}{2} \abs{v}{2} \idx \\
	& ~~~~ - \dfrac{\varepsilon}{2} \int  ( \nablah\rho \cdot \nablah ) v  \cdot v \idx,
\end{align*}
and the estimate
\begin{equation*}
	- \dfrac{\varepsilon}{2} \int  ( \nablah\rho \cdot \nablah ) v \cdot v \idx \leq \dfrac{\varepsilon}{16} \int \rho^{-1}\abs{\nablah \rho}{2} \abs{v}{2} \idx + \varepsilon \int \rho \abs{\nablah v}{2} \idx.
\end{equation*}
On the other hand, direct integration of \subeqref{eq:approximating-CPE}{1} over $ \Omega $  gives
\begin{equation}\label{BE-004}
	\dfrac{d}{dt} \int \rho \idx + \dfrac{\varepsilon}{16} \int \rho^{-2} (4\rho + \abs{\nablah \rho}{2}) \abs{\nablah \rho}{2}
	 = \varepsilon \int \rho^{-p_0} \idx.
\end{equation}
Furthermore, using identity \eqref{BE-002}, with  $ q = -p_0 $, yields
\begin{equation}\label{BE-005}
	\begin{aligned}
		& \dfrac{d}{dt} \biggl\lbrace \dfrac{1}{1+p_0} \int \rho^{-p_0} \idx \biggr\rbrace + \dfrac{\varepsilon p_0 (2p_0 + 1)}{16(p_0 + 1)} \int \rho^{-p_0 - 3} (4\rho + \abs{\nablah \rho}{2}) \abs{\nablah \rho}{2} \idx \\
		& ~~~~ + \dfrac{\varepsilon p_0}{p_0+1} \int \rho^{-2p_0 - 1} \idx = \int \rho^{-p_0} \dvh v\idx \leq \dfrac{\varepsilon p_0}{2(p_0+1)} \int \rho^{-2p_0 - 1} \idx \\
		& ~~~~  + \dfrac{p_0 + 1}{ \varepsilon p_0} \int \rho \abs{\mathcal D (v)}{2}\idx,
	\end{aligned}
\end{equation}
where we have used the fact that
\begin{equation*}
	\abs{\dvh v}{2} = \abs{\partial_1 v_1 + \partial_2 v_2}{2} \leq 2 (\abs{\partial_1 v_1}{2} + \abs{\partial_2 v_2}{2} ) \leq  2 \abs{\mathcal{D}(v)}{2}.
\end{equation*}
Therefore, multiplying \eqref{BE-005} by $ \dfrac{\varepsilon}{4} $ and adding the result to \eqref{BE-001} and \eqref{BE-004} 
yield
\begin{equation}\label{BE-006}
	\begin{aligned}
		& \dfrac{d}{dt} \biggl\lbrace \dfrac 1 2 \int \rho \abs{v}{2} \idx + \dfrac{1}{\gamma - 1} \int \rho^\gamma \idx + \int \rho \idx + \dfrac{\varepsilon}{4(1+p_0)} \int \rho^{-p_0} \idx \biggr\rbrace\\
		& ~~~~ + (1 -\dfrac{p_0 + 1}{4p_0}  ) \int \rho \abs{\mathcal{D}(v)}{2} \idx + \int \rho \abs{\dz v}{2} \idx   + (\sqrt{\varepsilon} - \varepsilon) \int \rho \abs{\nablah v}{2} \idx \\
		& ~~~~ + \dfrac{\varepsilon}{2} \int \rho^{-p_0} \abs{v}{2} \idx + \varepsilon \int \rho \abs{v}{5} \idx  + \dfrac{\varepsilon}{32} \int \rho^{-2} \abs{\nablah \rho}{4} \abs{v}{2} \idx\\
		& ~~~~ + \dfrac{\varepsilon}{16} \int \rho^{-1} \abs{\nablah \rho}{2} \abs{v}{2} \idx + \dfrac{\varepsilon(2\gamma - 1)\gamma}{16(\gamma - 1)} \int \rho^{\gamma-3}(4\rho + \abs{\nablah \rho}{2})\abs{\nablah \rho}{2} \idx  \\
		& ~~~~ + \dfrac{\varepsilon}{16} \int \rho^{-2}(4\rho + \abs{\nablah \rho}{2})\abs{\nablah \rho}{2} \idx + \dfrac{\varepsilon^2 p_0}{8(p_0+1)} \int \rho^{-2p_0-1}\idx \\
		& ~~~~ + \dfrac{\varepsilon^2 p_0 (2p_0 + 1)}{64(p_0+1)} \int \rho^{-p_0- 3}(4\rho + \abs{\nablah \rho}{2}) \abs{\nablah \rho}{2} \idx  \\
		& ~~ \leq \dfrac{\varepsilon \gamma}{\gamma - 1} \int \rho^{\gamma - p_0 - 1 }\idx + \varepsilon \int \rho^{-p_0}\idx \leq 2 \varepsilon \int \rho^{-p_0} \idx \\
		& ~~~~ + \varepsilon C_\gamma \int \rho \idx,
	\end{aligned}
\end{equation}
where $ C_\gamma > 0 $ is a constant depending only on $\gamma $ and we have employed the fact
\begin{equation*}
	p_0 \geq \max\lbrace1 , \gamma - 2 \rbrace.
\end{equation*}
Then for sufficiently small $ \varepsilon_0 $, with $ \varepsilon \in (0,\varepsilon_0) $, by applying the Gr\"onwall's inequality,  \eqref{BE-006} implies,
\begin{equation}\label{ene:basic-energy}
	\begin{aligned}
		& \sup_{0\leq t\leq T} \biggl\lbrace \int \rho \abs{v}{2} \idx + \int \rho^\gamma \idx + \int \rho \idx + \varepsilon \int \rho^{-p_0} \idx \biggr\rbrace \\
		& ~~~~ + \int_0^T \int \biggl( \rho \abs{\mathcal D (v)}{2} + \rho \abs{\dz v}{2} \biggr) \idx\,dt
		+ \varepsilon \int_0^T\int \biggl(\rho^{-p_0} \abs{v}{2} + \rho \abs{v}{5} \\
		& ~~~~ + \rho^{-2}\abs{\nablah \rho}{4}\abs{v}{2} + \rho^{-1} \abs{\nablah \rho}{2} \abs{v}{2} + \rho^{-2} \abs{\nablah \rho}{4} \biggr)\idx \,dt \\
		& ~~~~ + \varepsilon^2 \int_0^T \int \rho^{-2p_0-1}\idx \,dt \leq C_T,
	\end{aligned}
\end{equation}
where $ C_T $ depends only on $ T $ and the bounds of initial data given in \eqref{bound-initial-data}.

\subsection{The Bresh-Desjardins entropy estimate}

Multiplying \subeqref{eq:approximating-CPE}{1} with $ \rho^{-1} $ will lead to
\begin{equation*}
	\dt \log \rho + v \cdot\nablah \log \rho + \dvh v + \dz w = G(\rho)/\rho~~~~ \text{in} ~ L^2(0,T;L^{2}(\Omega)).
\end{equation*}
Now, we take the horizontal derivative $ \nablah $ to the above equation which gives 
\begin{align*}
	& \dt \nablah \log \rho + ( v\cdot\nablah) \nablah \log \rho + \nablah v \cdot \nablah \log \rho + \nablah \dvh v + \nablah \dz w \\
	& ~~~~ ~~~~ ~~~~ = \nablah (\dfrac{G(\rho)}{\rho}), ~~~~ \text{in} ~ L^2(0,T;H^{-1}(\Omega)).
\end{align*}
Taking the $L^2$-inner produce (or the duality action) of this equation with $ \rho \nablah\log \rho = 2 \rho^{1/2} \nablah \rho^{1/2} \in L^2(0,T;H^1(\Omega)) $ 
and using \subeqref{eq:approximating-CPE}{1} imply
\begin{equation}\label{BD-001}
	\begin{aligned}
		& \dfrac{d}{dt} \biggl\lbrace \dfrac{1}{2} \int \rho \abs{\nablah \log \rho}{2} \idx  \biggr\rbrace
		= \int \biggl( \dfrac 1 2  (G(\rho) - \dz(\rho w) ) \abs{\nablah\log \rho}{2} \\
		& ~~~~ + \nablah (\dfrac{G(\rho)}{\rho}) \cdot \rho \nablah\log\rho \biggr) \idx  - \int \biggl( \rho \nablah\log \rho \cdot (\nablah v \cdot \nablah \log \rho \\
		& ~~~~ + \nablah \dvh v  + \nablah \dz w) \biggr) \idx = - \int  G(\rho) ( \deltah \log \rho + \dfrac 1 2 \abs{\nablah \log\rho}{2} ) \idx \\
		& ~~~~ - \int \biggl( \rho^{-1} \nablah \rho \cdot\nablah v\cdot \nablah \rho - (\dvh v) \, \deltah \rho \biggr) \idx.  
	\end{aligned}
\end{equation}
On the other hand, taking the $ L^2 $-inner produce of \subeqref{eq:approximating-CPE}{2} with $ \nablah \log \rho = \nablah \rho / \rho $ will lead to
\begin{equation}\label{BD-002}
	\begin{aligned}
		 & \int \dt v \cdot \nablah \rho \idx + \int \rho^{-1}P'(\rho) \abs{\nablah\rho}{2} \idx + \int \biggl( (v \cdot\nablah) v \cdot \nablah \rho \\
		 & ~~~~ ~~~~ + w \dz v \cdot \nablah \rho \biggr) \idx  + (1+\sqrt{\varepsilon}) \int \dvh v \, \deltah \rho \idx \\
		 & ~~~~  - (1+\sqrt{\varepsilon}) \int \rho^{-1} (\nablah\rho \cdot \nablah) v \cdot \nablah \rho \idx = \int \widetilde F(v,\rho) \cdot\nablah \log \rho \idx.
	\end{aligned}
\end{equation}
As before, applying the horizontal derivative $ \nablah $ to \subeqref{eq:approximating-CPE}{1} yields,
\begin{equation*}
	\dt \nablah\rho + \nablah \dvh (\rho v) + \nablah \dz (\rho w) = \nablah G(\rho), ~~~~ \text{in} ~ L^2(0,T;H^{-1}(\Omega)).
\end{equation*}
Therefore, thanks to the above, we have the following equality, since $ v \in L^2(0,T; H^1(\Omega)) $,
\begin{align*}
	& \int \dt v \cdot \nablah \rho \idx = \dfrac{d}{dt} \int v \cdot\nablah \rho \idx - \int v \cdot \dt \nablah \rho \idx = \dfrac{d}{dt} \int v\cdot \nablah \rho \idx \\
	& ~~  + \int v \cdot (\nablah \dvh (\rho v) + \nablah \dz (\rho w) - \nablah G(\rho) ) \idx = \dfrac{d}{dt} \int v\cdot \nablah \rho \idx \\
	& ~~ + \int \biggl( v \cdot\nablah \dvh (\rho v) - \dz(\rho w) \dvh v + G(\rho) \dvh v \biggr) \idx  .
\end{align*}
Together with \eqref{BD-002}, we have
\begin{equation}\label{BD-003}
	\begin{aligned}
		& \dfrac{d}{dt} \int v \cdot\nablah \rho \idx + \gamma \int \rho^{\gamma-2} \abs{\nablah \rho}{2} \idx + \int \biggl( (v \cdot\nablah) v\cdot\nablah \rho \\
		& ~~~~ + v \cdot\nablah \dvh (\rho v) \biggr) \idx + \int \biggl( w \dz v \cdot\nablah \rho - \rho \dz w\, \dvh v \biggr) \idx \\
		& ~~~~ + \int G(\rho) \dvh v\idx + (1+\sqrt \varepsilon) \int \biggl( (\dvh v) \,\deltah \rho \\
		& ~~~~ ~~~~ - \rho^{-1} (\nablah \rho \cdot \nablah) v \cdot \nablah \rho \biggr) \idx  = \int \widetilde F (v,\rho) \cdot\nablah \log \rho\idx.
	\end{aligned}
\end{equation}
Notice that applying integration  by parts yields
\begin{align}
	& \int \biggl( (v \cdot\nablah) v\cdot\nablah \rho + v \cdot\nablah \dvh (\rho v) \biggr)  \idx = - \int \rho \nablah v : \nablah v^\top \idx {\nonumber} 
	\\
	& ~~~~ ~~~ = - 2 \int \rho \mathcal{D}(v) :\nablah v \idx  + \int \rho \abs{\nablah v}{2} \idx, {\nonumber} \\
	& \int \biggl( w \dz v \cdot\nablah \rho - \rho \dz w \dvh v \biggr) \idx = \int - \dz w \, \dvh (\rho v) \idx {\nonumber} \\
	& ~~~~ ~~~~ = \int \dz w (\dt \rho + \dz (\rho w) - G(\rho)) \idx  = \int \rho \abs{\dz w}{2}\idx, \label{BD-007}
\end{align}
where we have used the fact $ \rho $ is independent of $ z $ from \subeqref{eq:approximating-CPE}{3}. Therefore, multiplying \eqref{BD-001} by $ (1+\sqrt \varepsilon ) $ and adding the result to \eqref{BD-003} yield
\begin{equation}\label{BD-004}
	\begin{aligned}
		& \dfrac{d}{dt} \biggl\lbrace \dfrac{1+\sqrt \varepsilon}{2} \int \rho^{-1} \abs{\nabla \rho }{2} \idx + \int v \cdot \nablah \rho \idx  \biggr\rbrace + \int \rho \abs{\nablah v}{2} \,dx \\
		& ~~~~ + \int \rho \abs{\dz w}{2}\idx + \gamma\int \rho^{\gamma-2} \abs{\nablah \rho}{2} \idx + (1+\sqrt\varepsilon)\int  G(\rho) ( \deltah \log \rho \\
		& ~~~~ + \dfrac 1 2 \abs{\nablah \log\rho}{2} ) \idx  = 2 \int \rho \mathcal{D} v :\nablah v \idx  - \int G(\rho) \dvh v\idx \\
		& ~~~~ + \int \widetilde F (v,\rho) \cdot\nablah \log \rho \idx =: I_1 + I_2 + I_3.
	\end{aligned}
\end{equation}
By virtue of the identity $$ \deltah \log \rho + \dfrac{1}{2} \abs{\nablah \log\rho}{2} = 2 \rho^{-1/2} \deltah \rho^{1/2} = \rho^{-1}\deltah \rho - \dfrac{1}{2} \rho^{-2} \abs{\nablah \rho}{2},  $$
we have
\begin{align*}
	& \int  G(\rho) ( \deltah \log \rho + \dfrac 1 2 \abs{\nablah \log\rho}{2} ) \idx 
	= 2 \varepsilon \int \abs{\deltah \rho^{1/2}}{2}\idx \\
	& ~~~~ + 2 \varepsilon \int \dvh(\abs{\nablah \rho^{1/2}}{2} \nablah \rho^{1/2}) \deltah \rho^{1/2} \idx + 2\varepsilon \int \rho^{-p_0 - 1/2} \deltah \rho^{1/2} \idx \\
	& ~~~~ = 2 \varepsilon\int \abs{\deltah \rho^{1/2}}{2} \idx + 2 \varepsilon \int \abs{\nablah \rho^{1/2}}{2} \abs{\nablah^2 \rho^{1/2}}{2} \idx \\
	& ~~~~ + \varepsilon \int \abs{\nablah \abs{\nablah \rho^{1/2}}{2}}{2} \idx + \varepsilon(4p_0 + 2) \int \rho^{-p_0-1} \abs{\nablah \rho^{1/2}}{2} \idx.
\end{align*}
Moreover, by applying integration by parts and the Young's inequality, the right-hand side of \eqref{BD-004} satisfies the following estimates:
\begin{align*}
	& I_1 \leq \dfrac{1}{2} \int \rho \abs{\nablah v}{2}\idx + 2 \int \rho \abs{\mathcal D (v)}{2} \idx,\\
	& I_2 = - \int \bigl(\varepsilon \rho^{1/2} \deltah \rho^{1/2} + \varepsilon \rho^{1/2} \dvh ( \abs{\nablah \rho^{1/2}}{2} \nablah \rho^{1/2} ) + \varepsilon \rho^{-p_0}\bigr) \dvh v \idx \\
	& ~~~~ \leq \dfrac{\varepsilon^2}{4} \int \abs{\deltah \rho^{1/2} + \dvh(\abs{\nablah \rho^{1/2}}{2}\nablah \rho^{1/2} )}{2} \idx + \dfrac{\varepsilon^2}{4} \int \rho^{-2p_0-1} \idx \\
	& ~~~~  + 2 \int \rho \abs{\dvh v}{2} \idx,\\
	& I_3  
	= \varepsilon \int \biggl( 2 (\abs{\nablah \rho^{1/2}}{2} \nablah \rho^{1/2} \cdot \nablah) v \cdot \nablah \rho^{1/2} - \rho^{-p_0-1} v \cdot \nablah \rho - \abs{v}{3} v \cdot\nablah \rho \biggr) \idx \\
	& ~~~~ = - \varepsilon \int 2 \sum_{i,j\in\lbrace 1,2\rbrace} v_j \partial_i (\abs{\nablah \rho^{1/2}}{2} \partial_i \rho^{1/2} \partial_j \rho^{1/2}) \idx - \varepsilon \int \dfrac{1}{p_0} \rho^{-p_0} \dvh v \idx \\
	& ~~~~ ~~~~ - \varepsilon \int \abs{v}{3} v \cdot\nablah \rho \idx =: I_3' + I_3'' + I_3''', ~~ \text{with} \\
	& I_3' \leq \varepsilon \int \abs{\nablah \rho^{1/2}}{2} \abs{\nablah^2\rho^{1/2}}{2} \idx + \varepsilon C \int \rho^{-2}\abs{\nablah \rho}{4} v^2 \idx, \\
	& I_3'' \leq \varepsilon^2 \int \rho^{-2p_0-1} \idx + C \int \rho \abs{\dvh v}{2}\idx \leq \varepsilon^2 \int \rho^{-2p_0-1} \idx \\
	& ~~~~ + C \int \rho \abs{\mathcal D v}{2}\idx, \\
	& I_3''' \leq \varepsilon \int \underbrace{\rho^{-1/5}}_{L^{20}} \cdot \underbrace{\rho^{7/10} |v|^{7/2}}_{L^{10/7}} \cdot\underbrace{\rho^{-1/2} \abs{\nablah \rho}{} |v|^{1/2}}_{L^4} \idx \leq  \varepsilon \int \rho^{-2} \abs{\nablah \rho}{4} v^2 \idx \\
	& ~~~~ + \varepsilon \int \rho v^5 \idx + \varepsilon C \int \rho^{-4} \idx  \leq  \varepsilon \int \rho^{-2} \abs{\nablah \rho}{4} v^2 \idx + \varepsilon \int \rho v^5 \idx \\
	& ~~~~ + \varepsilon C \int \rho^{-p_0} \idx + \varepsilon C \abs{\Omega}{} .
\end{align*}
where we have used the fact
\begin{equation*}
	p_0 \geq 4,
\end{equation*}
and
\begin{align*}
	& \int \rho^{-4} \idx = \int_{\lbrace \rho > 1\rbrace} \rho^{-4} \idx  + \int_{\lbrace \rho \leq 1 \rbrace} \rho^{-4} \idx \\
	& ~~~~ \leq \int 1 \idx + \int_{\lbrace \rho \leq 1 \rbrace} \rho^{-p_0} \idx \leq \abs{\Omega}{} + \int \rho^{-p_0}\idx.
\end{align*}
%
Therefore, after summing up the above inequalities,  \eqref{BD-004} implies
\begin{equation}\label{BD-005}
	\begin{aligned}
		& \dfrac{d}{dt} \biggl\lbrace \dfrac{1+\sqrt \varepsilon}{2} \int \rho^{-1} \abs{\nablah \rho }{2} \idx + \int v \cdot \nablah \rho \idx  \biggr\rbrace + \dfrac{1}{2} \int \rho \abs{\nablah v}{2} \,dx \\
		& ~~ + \int \rho \abs{\dz w}{2}\idx  + \gamma\int \rho^{\gamma-2} \abs{\nablah \rho}{2} \idx + \varepsilon \int \abs{\deltah \rho^{1/2}}{2} \idx \\
		& ~~ + \varepsilon \int \abs{\nablah \rho^{1/2}}{2} \abs{\nablah^2 \rho^{1/2}}{2} \idx  \leq C \int \rho \abs{\mathcal D (v)}{2} \idx \\
		& ~~ + \varepsilon^2 C \int \rho^{-2p_0-1}\idx + \varepsilon C \int \rho^{-2}\abs{\nablah \rho}{4} v^2 \idx + \varepsilon C \int \rho v^5 \idx \\
		& ~~ + \varepsilon C \int \rho^{-p_0}  \idx + \dfrac{\varepsilon^2}{4} \int \abs{\deltah \rho^{1/2} + \dvh(\abs{\nablah \rho^{1/2}}{2}\nablah \rho^{1/2} )}{2} \idx + \varepsilon C .
	\end{aligned}
\end{equation}
On the other hand, \subeqref{eq:approximating-CPE}{1} can be written as
\begin{equation}\label{BD-006}
\begin{aligned}
	& 2 \dt \rho^{1/2} - \varepsilon \deltah \rho^{1/2} - \varepsilon  \dvh ( \abs{\nablah \rho^{1/2}}{2} \nablah \rho^{1/2} ) \\
	& ~~~~ = - 2 v \cdot\nablah \rho^{1/2} - \rho^{1/2} (\dvh v + \dz w) + \varepsilon \rho^{-p_0 - 1/2}.
\end{aligned}
\end{equation}
After multiplying \eqref{BD-006} with $ - \varepsilon \deltah \rho^{1/2} - \varepsilon \dvh (\abs{\nablah \rho^{1/2}}{2} \nablah \rho^{1/2}) $ and integrating the resultant over $ \Omega $, we will have
\begin{equation}\label{BD-1001}
\begin{aligned}
	& \dfrac{d}{dt} \biggl\lbrace  \varepsilon \int \abs{\nablah \rho^{1/2}}{2} \idx + \dfrac{\varepsilon}{2} \int \abs{\nablah \rho^{1/2}}{4} \idx \biggr\rbrace  \\
	& ~~~~ ~~~~ + \varepsilon^2 \int \abs{\deltah \rho^{1/2}  + \dvh(\abs{\nablah \rho^{1/2}}{2}\nablah \rho^{1/2} )}{2} \idx \\
	& ~~~~ ~~~~ + \varepsilon^2 (2p_0 + 1) \int \rho^{-p_0-1} (\abs{\nablah \rho^{1/2}}{2}  + \abs{\nablah \rho^{1/2}}{4}) \idx \\
	& ~~~~  = \int \rho^{1/2}( \dvh v + \dz w ) (\varepsilon \deltah \rho^{1/2}  + \varepsilon  \dvh ( \abs{\nablah \rho^{1/2}}{2} \nablah \rho^{1/2} )) \idx \\
	& ~~~~ ~~~~ + \int 2 (v \cdot\nablah) \rho^{1/2} (\varepsilon \deltah \rho^{1/2} + \varepsilon  \dvh ( \abs{\nablah \rho^{1/2}}{2} \nablah \rho^{1/2} )) \idx  \\
	& ~~~~ \leq \dfrac{\varepsilon^2}{4} \int \abs{\deltah \rho^{1/2} + \dvh(\abs{\nablah \rho^{1/2}}{2}\nablah \rho^{1/2} )}{2} \idx + C \int \rho \abs{\dvh v}{2} \idx \\
	& ~~~~ ~~~~ + \dfrac{\varepsilon}{2} \int \biggl( \abs{\deltah \rho^{1/2}}{2} + \abs{\nablah \rho^{1/2}}{2} \abs{\nablah^2\rho^{1/2}}{2} \biggr) \idx \\
	& ~~~~ ~~~~ + \varepsilon C  \int \biggl( \rho^{-1} \abs{\nablah \rho}{2} v^2 + \rho^{-2}\abs{\nablah \rho}{4} v^2 \biggr) \idx.
\end{aligned}
\end{equation}
After adding this estimate to \eqref{BD-005}, one will end up with
\begin{equation}\label{BD-007}
	\begin{aligned}
		& \dfrac{d}{dt} \biggl\lbrace \dfrac{1+\sqrt \varepsilon}{2} \int \rho^{-1} \abs{\nablah \rho }{2} \idx + \int v \cdot \nablah \rho \idx + \varepsilon \int \abs{\nablah \rho^{1/2}}{2} \idx \\
		& ~~~~ + \dfrac{\varepsilon}{2} \int \abs{\nablah \rho^{1/2}}{4} \idx  \biggr\rbrace + \dfrac{1}{2} \int \rho \abs{\nablah v}{2} \,dx + \int \rho \abs{\dz w}{2}\idx \\
		& ~~~~ + \gamma\int \rho^{\gamma-2} \abs{\nablah \rho}{2} \idx + \dfrac{\varepsilon}{2} \int \abs{\deltah \rho^{1/2}}{2} \idx + \dfrac{\varepsilon}{2} \int \abs{\nablah \rho^{1/2}}{2} \abs{\nablah^2 \rho^{1/2}}{2} \idx \\
		& ~~~~ + \dfrac{\varepsilon^2}{2} \int \abs{\deltah \rho^{1/2} + \dvh(\abs{\nablah \rho^{1/2}}{2}\nablah \rho^{1/2} )}{2} \idx \\
		& ~~~~ + \varepsilon^2 (2p_0 + 1) \int \rho^{-p_0-1} (\abs{\nablah \rho^{1/2}}{2} + \abs{\nablah \rho^{1/2}}{4}) \idx \\
		& ~~~~ \leq C \int \rho \abs{\mathcal D (v)}{2} \idx + \varepsilon C \int \biggl( \rho^{-2}\abs{\nablah \rho}{4} v^2 + \rho^{-1} \abs{\nablah\rho}{2}v^2 \biggr) \idx \\
		& ~~~~ + \varepsilon C \int \rho v^5 \idx  + \varepsilon^2 C \int \rho^{-2p_0-1}\idx + \varepsilon C \int \rho^{-p_0}  \idx + \varepsilon C.
	\end{aligned}
\end{equation}
Moreover, for $ \varepsilon_0 $ sufficiently small and $ \varepsilon \in (0,\varepsilon_0) $,
\begin{align*}
	& \dfrac{\varepsilon^2}{2} \int \abs{\deltah \rho^{1/2} + \dvh(\abs{\nablah \rho^{1/2}}{2}\nablah \rho^{1/2} )}{2} \idx \\
	& ~~~~ \geq \dfrac{\varepsilon^2}{4} \int \abs{\dvh(\abs{\nablah \rho^{1/2}}{2} \nablah \rho^{1/2})}{2} \idx - \dfrac{\varepsilon}{4} \int \abs{\deltah \rho^{1/2}}{2} \idx.
	\end{align*}
	{Applying integration by parts twice yields}
	\begin{align*}
	& \int \abs{\dvh(\abs{\nablah \rho^{1/2}}{2} \nablah \rho^{1/2})}{2} \idx = \int \biggl( \sum_{i,j\in \lbrace 1, 2 \rbrace} \partial_j(\abs{\nablah \rho^{1/2}}{2} \partial_{i}\rho^{1/2} ) \\
	& ~~~~ ~~~~ \times \partial_i ( \abs{\nablah \rho^{1/2}}{2} \partial_j \rho^{1/2} ) \biggr) \idx = \int \biggl( \abs{\nablah \rho^{1/2}}{4} \abs{\nablah^2 \rho^{1/2}}{2} \\
	& ~~~~ + \abs{\nablah \abs{\nablah \rho^{1/2}}{2} \cdot \nablah \rho^{1/2}}{2} +  \dfrac{1}{2} \nablah \abs{\nablah \rho^{1/2}}{4} \cdot\nablah  \abs{\nablah\rho^{1/2}}{2} \biggr) \idx \\
	& ~~~~ = \int \biggl( \abs{\nablah \rho^{1/2}}{4} \abs{\nablah^2 \rho^{1/2}}{2} + \abs{\nablah \abs{\nablah \rho^{1/2}}{2} \cdot \nablah \rho^{1/2}}{2} \\
	& ~~~~ +  \abs{\nablah \rho^{1/2}}{2} \abs{\nablah  \abs{\nablah\rho^{1/2}}{2}}{2} \biggr) \idx.
\end{align*}
Therefore, after plugging the above into \eqref{BD-007} and integrating the result with respect to the temporal variable, we obtain the following estimate,
\begin{equation}\label{ene:BD-entropy}
	\begin{aligned}
		& \sup_{0\leq t\leq T} \biggl\lbrace \int \rho^{-1} \abs{\nablah \rho}{2} \idx  + \varepsilon \int \rho^{-2} \abs{\nablah \rho}{4}\idx \biggr\rbrace + \int_0^T \int \biggl( \rho\abs{\nablah v}{2} \\
		& ~~~~ + \rho \abs{\dz w}{2} \biggr) \idx\,dt + \int_0^T \int  \rho^{\gamma-2}\abs{\nablah \rho}{2} \idx \,dt  + \varepsilon \int_0^T \int \biggl( \abs{\deltah \rho^{1/2}}{2} \\
		& ~~~~ + \abs{\nablah \rho^{1/2}}{2} \abs{\nablah^2 \rho^{1/2}}{2} \biggr) \idx \,dt + \varepsilon^2 \int_0^T \int  \abs{\nablah \rho^{1/2}}{4} \abs{\nablah^2 \rho^{1/2}}{2} \idx \,dt\\
		& ~~~~ \leq C_T,
	\end{aligned}
\end{equation}
where we have also employed the energy estimate \eqref{ene:basic-energy}. Here $ C_T $ depends only on $ T $ and the bounds of initial data given in \eqref{bound-initial-data}.

\subsection{The Mellet-Vasseur estimate}
We take the inner product of \subeqref{eq:approximating-CPE}{2} with $ (1+\log(e+v^2)) v $ and integrate the resultant to get
	\begin{align}\label{MV-001}
		& {\nonumber} \dfrac{d}{dt} \bigl\lbrace \dfrac{1}{2} \int \rho (e + v^2) \log (e + v^2) \idx \bigr\rbrace  + \int \rho (1 + \log(e + v^2)) (\abs{\mathcal{D} (v)}{2} + \abs{\dz v}{2} \\
		& {\nonumber} ~~~~ ~~~~ + \sqrt{\varepsilon} \abs{\nablah v}{2} ) \idx
		 = - \int \nablah \rho^\gamma \cdot (1+ \log (e + v^2)) v \idx \\
		 &{\nonumber} ~~~~ - \sum_{i,j,k \in \lbrace 1,2 \rbrace} \int \biggl( \rho \dfrac{2 v_j v_k }{e+v^2}\partial_i v_k \mathcal{D}(v)_{ij} + \rho \dfrac{2 v_i v_j}{e+ v^2}\dz v_i \dz v_j \\
		 &{\nonumber} ~~~~ ~~~~ + \sqrt\varepsilon \rho \dfrac{2 v_j v_k}{e + v^2} \partial_i v_j \partial_i v_k \biggr) \idx + \int \widetilde F (v,\rho) \cdot (1+\log(e+v^2)) v \idx \\
		&{\nonumber} ~~~~ + \dfrac{1}{2} \int (e + v^2) \log(e + v^2) G(\rho) \idx = - \int \nablah \rho^\gamma \cdot (1+ \log (e + v^2)) v \idx \\
		&{\nonumber} ~~~~ - \int \rho \sum_{i,j,k \in \lbrace 1,2 \rbrace} \biggl( \dfrac{2 v_j v_k }{e+v^2}\partial_i v_k \mathcal{D}(v)_{i,j}  + \dfrac{2 v_i v_j}{e+ v^2}\dz v_i \dz v_j \\
		&{\nonumber} ~~~~ ~~~~ +  \sqrt\varepsilon \dfrac{2 v_j v_k}{e + v^2} \partial_i v_j \partial_i v_k  \biggr) \idx + \dfrac{\varepsilon}{2} \int \rho^{1/2} \deltah \rho^{1/2} (e+v^2)\log(e+v^2) \idx \\
		&{\nonumber} ~~~~ - \varepsilon \int \rho \abs{v}{3}(1+\log(e+v^2))v^2 \idx + \varepsilon\int \biggl( \dfrac{1}{2} \rho^{-p_0}(e+v^2) \log(e+v^2) \\
		&{\nonumber} ~~~~ ~~~~ - \rho^{-p_0} (1+\log(e+v^2)) v^2 \biggr) \idx  + \dfrac{\varepsilon}{2} \int \biggl( ( \rho^{1/2}\dvh(\abs{\nablah \rho^{1/2}}{2} \nablah \rho^{1/2}) \\
		&{\nonumber} ~~~~ ~~~~ - \dvh(\rho^{1/2} \abs{\nablah \rho^{1/2}}{2}\nablah \rho^{1/2} ) ) (e+v^2) \log(e +v^2) \biggr) \idx \\
		&{\nonumber} ~~~~ \leq - \int \nablah \rho^\gamma \cdot (1+\log(e+v^2)) v \idx + C \int \rho (\abs{\nablah v}{2} + \abs{\dz v}{2}) \idx \\
		& ~~~~ + \dfrac{\varepsilon^2}{2} \int \rho^{1/2} \deltah \rho^{1/2} (e+ v^2) \log(e+ v^2) \idx  + \varepsilon C \int \rho^{-p_0} \idx .
	\end{align}
The first two integrals on the right-hand side of the last two lines of \eqref{MV-001} can be calculated after applying integration by parts and Young's inequality as follows:
\begin{align*}
	& - \int \nablah \rho^\gamma \cdot (1+\log(e+v^2)) v \idx 
	= \int \rho^\gamma (1+\log(e+v^2)) \dvh v \idx \\
	& ~~~~ + \int \rho^\gamma \dfrac{2 (v \cdot \nablah) v \cdot v}{e+v^2} \idx \leq C \int \rho \abs{\nablah v}{2} \idx + C \int \rho^{2\gamma - 1} \log^2(e+v^2) \idx,\\
	& \dfrac{\varepsilon^2}{2} \int \rho^{1/2} \deltah \rho^{1/2} (e+ v^2) \log(e+ v^2) \idx = - \dfrac{\varepsilon^2}{2} \int \biggl( (e+v^2) \log(e+v^2) \abs{\nablah \rho^{1/2}}{2} \\
	& ~~~~ ~~~~ + (1+ \log (e+v^2)) (\nablah \rho \cdot \nablah) v \cdot v \biggr) \idx \\
	& \leq - \dfrac{\varepsilon^2}{2} \int (e+v^2) \log(e+v^2) \abs{\nablah \rho^{1/2}}{2} \idx + \varepsilon^2 \int \rho (1 + \log(e+v^2)) \abs{\nablah v}{2} \idx \\
	& ~~~~ + \dfrac{\varepsilon^2}{16} \int \underbrace{ \rho^{-1} \abs{\nablah \rho}{2} v^2 (1+ \log (e +v^2)) }_{\leq 8 v^2 \log(e+v^2)|\nablah \rho^{1/2}|^{2}  }\idx\\
	& \leq - \dfrac{\varepsilon^2}{2} \int e \log (e+v^2) \abs{\nablah \rho^{1/2}}{2} \idx + \varepsilon^2 \int \rho (1+ \log(e+v^2)) \abs{\nablah v}{2} \idx.
\end{align*}
Therefore, after substituting the above estimates in \eqref{MV-001}
and integrating with respect to the temporal variable, one obtains,
\begin{equation}\label{MV-002}
	\begin{aligned}
		& \sup_{0\leq t\leq T} \int \rho(e+v^2) \log(e+v^2) \idx + \int_0^T \int \biggl\lbrack \rho (1+\log(e+v^2)) (\abs{\mathcal{D} (v)}{2} \\
		& ~~~~ ~~~~ + \abs{\dz v}{2}  + (\sqrt{\varepsilon} - \varepsilon^2)\abs{\nablah v}{2} ) \biggr\rbrack \idx\,dt  \leq C + \\
		& ~~~~ ~~~~ +  C \int_0^T \int \rho^{2\gamma- 1} \log^2 (e+v^2) \idx \,dt,
	\end{aligned}
\end{equation}
where we have also used estimates \eqref{ene:basic-energy}, \eqref{ene:BD-entropy}. Meanwhile, the last integral on the right-hand side of \eqref{MV-002} can be estimated by interpolation and embedding inequalities. Indeed, for any $ \delta \in (0,1) $, small enough satisfying
\begin{equation*}
	\dfrac{4\gamma - 2 - 2\delta}{1-\delta} > 1,
\end{equation*}
since $ \rho = \rho(x,y,t) $ is independent of the $ z $ variable, applying the H\"older and the Young's inequalities yields,
\begin{align*}
	& \int_0^T \int \rho^{2\gamma- 1} \log^2 (e+v^2) \idx \,dt \leq \int_0^T \biggl( \bigl( \int \rho \log^{2/\delta}(e+v^2) \idx \bigr)^{\delta} \\
	& ~~~~ \times \bigl( \int  \rho^{(2\gamma-1-\delta)/(1-\delta)}\idx  \bigr)^{1-\delta} \biggr) \,dt  \leq \delta \int_0^T \int \rho (e + v^2) \idx \,dt \\
	& ~~~~ + (1-\delta) \int_0^T \normh{\rho^{1/2}}{L^{(4\gamma - 2 - 2\delta)/(1-\delta)}}^{(4\gamma - 2 - 2\delta)/(1-\delta)} 
	\,dt  \leq \delta C \\
	& ~~~~ + (1- \delta) \int_0^T \biggl( \normh{\rho^{1/2}}{L^2}^{2} \normh{\nablah\rho^{1/2}}{L^2}^{(4\gamma-4)/(1-\delta)} + \normh{\rho^{1/2}}{L^2}^{(4\gamma-2-2\delta)/(1-\delta)}\biggr) \,dt \\
	& ~~~~ \leq C_T,
\end{align*}
where we have used the two-dimensional Gagliardo-Nirenberg inequalities in the third inequality and the facts that
\begin{equation*}
	\normh{\rho^{1/2}}{L^2} = \bigl(\int \rho \idx\bigr)^{1/2} < C , ~ \normh{\nablah\rho^{1/2}}{L^2} = \dfrac{1}{2}\bigl( \int \rho^{-1} \abs{\nablah \rho}{2} \idx \bigr)^{1/2} < C_T,
\end{equation*}
as consequences of \eqref{ene:basic-energy} and \eqref{ene:BD-entropy}.
Therefore, from \eqref{MV-002}, with $ \varepsilon_0 $ small enough, and $ \varepsilon \in (0, \varepsilon_0) $, we have established the Mellet-Vasseur type estimate
\begin{equation}\label{ene:MV-estimate}
		\sup_{0\leq  t\leq T} \int \rho(e + v^2) \log(e+ v^2) \idx \leq C_T.
\end{equation}
Here $ C_T $ depends only on $ T $ and the bounds of initial data given in \eqref{bound-initial-data}, but is independent of $ \varepsilon $.

\section{The existence of the approximating solutions}\label{sec:aprx-sol-exist}
In this section, we aim at showing the global existence of the strong solutions $ (\rho_\varepsilon, v_\varepsilon) $ to the approximating system \eqref{eq:approximating-CPE}. That is, we will show that there are approximating solutions $ (\rho_\varepsilon, v_\varepsilon) $ satisfying Definition \ref{def:aprxm-sols}, for any given $ \varepsilon \in (0,\varepsilon_0) $, with $ \varepsilon_0 $ small enough. This is done in two steps. We first present some a priori estimates, which show the regularity of the global strong solutions, provided they exist.
Then we employ a modified Galerkin approximating scheme to show the existence of the approximating solutions.


\subsection{The global a priori estimates}\label{sec:aprx-propri-est}
 We shall first present the global a priori estimates. In fact, we will show the following:
\begin{prop}\label{prop:aprx-priori-est}
Let $ T \in (0,\infty) $, and suppose that there is a smooth enough solution $(\rho_\varepsilon, v_\varepsilon)$ to \eqref{eq:approximating-CPE} on the interval $[0,T] $, with the approximating initial data $ ( \rho_{\varepsilon,0}, v_{\varepsilon,0}) $ defined in Definition \ref{def:aprxm-initial}. Then there exists a positive constant $ C_{\varepsilon, T} $ such that the following estimates hold,
	\begin{align*}
		& C_{\varepsilon, T}^{-1} < \rho_\varepsilon < C_{\varepsilon, T}, \\
		& \sup_{0\leq t\leq T} \bigl\lbrace \norm{\nabla \rho_\varepsilon(t)}{L^2} + \norm{\nablah \rho_\varepsilon(t)}{L^4} + \norm{v_\varepsilon(t)}{L^2} + \norm{\nabla v_\varepsilon(t)}{L^2}  \bigr\rbrace \\
		& ~~~~ + \norm{\nabla^2 \rho_\varepsilon^{1/2}}{L^2(\Omega\times(0,T))} + \norm{\nabla \rho_\varepsilon}{L^{30}(\Omega\times(0,T))} + \norm{\dt \rho_\varepsilon}{L^2(\Omega\times(0,T))} \\
		& ~~~~ + \norm{\abs{\nablah \rho_\varepsilon^{1/2}}{2} \nablah^2 \rho_\varepsilon^{1/2}}{L^2(\Omega\times(0,T))}  + \norm{v_{\varepsilon}}{L^5(\Omega\times (0,T))} + \norm{\nabla v_\varepsilon}{L^2(\Omega\times(0,T))}\\
		& ~~~~ + \norm{\nabla^2 v_\varepsilon}{L^2(\Omega\times (0,T))} + \norm{\dt v_{\varepsilon}}{L^2(\Omega\times(0,T))} + \norm{v_\varepsilon}{L^{10}(\Omega\times(0,T))} \\
		& ~~~~ + \norm{\nabla v_\varepsilon}{L^{10/3}(\Omega\times(0,T))} + \norm{\dz w_{\varepsilon}}{L^2(\Omega\times(0,T))} < C_{\varepsilon, T}.
	\end{align*}
	Here $ w_\varepsilon $ is given by
	\begin{align*}
		& w_\varepsilon(x,y,z,t) = - \rho_\varepsilon^{-1}(x,y,t) \int_0^z \dvh \bigl(\rho_\varepsilon(x,y,t) ( v_\varepsilon(x,y,z',t) \\
		& ~~~~ ~~~~ ~~~~ ~~~~ - \int_0^1 v_\varepsilon(x,y,z'',t) \,dz'' )\bigr)\,dz'.
	\end{align*}
	In particular, we have
	\begin{gather*}
		\rho_{\varepsilon}^{1/2} \in L^\infty(0,T;W^{1,4}(\Omega)) \cap L^2(0,T;H^2(\Omega)), \nablah \rho_\varepsilon \in L^{30}(\Omega\times(0,T)), \\
		 \abs{\nablah \rho_\varepsilon^{1/2}}{2} \nablah^2 \rho_\varepsilon^{1/2} \in L^2(\Omega\times(0,T)), \dt \rho_\varepsilon\in L^2(\Omega\times(0,T)), \\
		v_\varepsilon \in L^\infty(0,T;H^1(\Omega))\cap L^2(0,T;H^2(\Omega)) \cap L^5(\Omega\times(0,T)) \cap L^{10}(\Omega\times(0,T)), \\
		\nabla v_\varepsilon \in L^{10/3}(\Omega\times(0,T)), \dt v_\varepsilon\in L^2(\Omega\times(0,T)), \dz w_\varepsilon \in L^2(\Omega\times(0,T)),
	\end{gather*}
	and therefore $ (\rho_\varepsilon, v_\varepsilon) $ is the approximating solution satisfying Definition \ref{def:aprxm-sols}.
\end{prop}
\begin{pf}
	This is a direct consequence of  \eqref{ene:density-UL}, \eqref{GR-001}, \eqref{ene:H^1-of-density}, \eqref{ene:H^1-of-v}, below.
\end{pf}

As before, we shall denote $ (\rho_\varepsilon, v_\varepsilon, w_\varepsilon) $ as $ (\rho, v, w) $, for the sake of convenience. Moreover, all the estimates in this section may depend on $ \varepsilon > 0 $. In particular, through out this section, the generic constant $ C $ may depend on $ \varepsilon \in (0, \varepsilon_0) $, $ T >  0 $ and the bounds of initial data given in Definition \ref{def:aprxm-initial}, unless it is stated differently.

\subsubsection{The upper and lower bounds of the density}\label{sec:ul-bounds}
We first establish the upper and lower bound of the density $ \rho $. That is, we will establish the following:
\begin{prop}\label{lm:density-ul-bounds}
	Let $ T \in (0,\infty) $, and
	consider a solution $ (\rho, v) $ on $ [0,T] $ of equation \eqref{eq:approximating-CPE}, satisfying the estimates in Proposition \ref{prop:uniform-est}. Then there is a constant $ C \in (0, \infty) $, which may depend on $ \varepsilon $ and $ T $, such that
	\begin{equation}\label{ene:density-UL}
		C^{-1} \leq \rho \leq C .
	\end{equation}
\end{prop}
\begin{pf}
It suffices to show the upper and lower bounds of $ \rho^{1/2} $. Indeed,
define
\begin{equation}\label{UL-001}
	\eta = \eta(x,y,t):= \rho^{1/2}(x,y,t),~ \sigma = \sigma(x,y,t): = \rho^{-1/2}(x,y,t).
\end{equation}
Then the two-dimensional Sobolev embedding inequality yields
\begin{align*}
	& \sup_{0\leq t\leq T} \normh{\eta}{L^\infty} \leq C \sup_{0\leq t\leq T} \bigl\lbrace \normh{\eta}{L^2} + \normh{\nablah\eta}{L^4} \bigr\rbrace = C \sup_{0\leq t\leq T} \bigl\lbrace \bigl( \int \rho \idx \bigr)^{1/2} \\
	& ~~~~ ~~~~ + \dfrac{1}{2} \bigl(\int \rho^{-2} \abs{\nablah \rho}{4} \idx\bigr)^{1/4} \bigr\rbrace < C,
\end{align*}
as a consequence of \eqref{ene:basic-energy} and \eqref{ene:BD-entropy}.
In particular,
\begin{equation}\label{UL-002}
	\sup_{(x,y,z,t)\in \Omega\times [0,T]} \rho(x,y,z,t) < C.
\end{equation}
On the other hand, recall $ \sigma = \rho^{-1/2} $. Thanks to \subeqref{eq:approximating-CPE}{1}, $ \sigma $ satisfies
\begin{equation}\label{UL-003}
	\begin{aligned}
		& 2 \dt \sigma + 2 v \cdot\nablah \sigma - \sigma (\dvh v + \dz w) + \varepsilon \sigma^{3+2p_0} + 2 \varepsilon \sigma^{-1} \abs{\nablah \sigma}{2}\\
		& ~~~~ + 2 \varepsilon \sigma^{-5} \abs{\nablah \sigma}{4} = \varepsilon (\deltah \sigma + \dvh (\sigma^{-4} \abs{\nablah \sigma}{2} \nablah \sigma)).
	\end{aligned}
\end{equation}
We will employ a De Giorgi type procedure to obtain the upper bound of $ \sigma $. Multiple \eqref{UL-003} with $ (\sigma - k)_+ = (\sigma - k) \mathbbm 1_{\lbrace \sigma > k \rbrace} $ with $ k\geq \norm{\sigma(\cdot,0)}{L^\infty} = \norm{\rho_0^{-1/2}}{L^\infty} $, where $ \mathbbm{1}_{\lbrace \sigma > k \rbrace} $ is the characteristic function of the set $\lbrace \sigma > k \rbrace \subset \Omega\times (0,T)$, and integrate the resultant over $ \Omega $. It follows, after applying the Young's inequality,
\begin{equation}\label{UL-007}
	\begin{aligned}
		& \dfrac{d}{dt} \int \abs{(\sigma-k)_+}{2} \idx + \varepsilon \int \biggl( \abs{\nablah(\sigma-k)_+}{2} + \sigma^{-4}\abs{\nablah (\sigma-k)_+}{4} \biggr) \idx \\
		& ~~~~ = - 3 \int (\sigma-k)_+  v \cdot\nablah \sigma \idx - \int \sigma v \cdot \nablah (\sigma - k)_+ \idx   \\
		& ~~~~ ~~~~ -  \varepsilon \int ( \sigma^{3+2p_0} + 2\sigma^{-1} \abs{\nablah \sigma}{2} +2 \sigma^{-5} \abs{\nablah \sigma}{4}) (\sigma-k)_+\idx \\
		& ~~~~ \leq C \int  \mathbbm{1}_{\lbrace\sigma > k\rbrace} (\abs{(\sigma-k)_+}{} + \abs{\sigma}{}) \abs{v}{} \abs{\nablah (\sigma-k)_+}{} \idx \\
		& ~~~~ \leq  \dfrac{\varepsilon}{4} \int \sigma^{-4} \abs{\nablah(\sigma-k)_+}{4} \idx + C \varepsilon^{-1/3} \int \mathbbm{1}_{\lbrace\sigma > k\rbrace} \sigma^{8/3} \abs{v}{4/3} \idx .
	\end{aligned}
\end{equation}
Integrating the above expression with respect to the temporal variable yields
\begin{equation}\label{UL-004}
	\begin{aligned}
		& \sup_{0\leq t\leq T} \int \abs{(\sigma-k)_+}{2} \idx + \varepsilon \int_0^T \int \abs{\nablah(\sigma-k)_+}{2} \\
		& ~~~~ + \sigma^{-4}\abs{\nablah (\sigma-k)_+}{4} \idx \,dt \leq C \varepsilon^{-1/3} \int_0^T \int \mathbbm{1}_{\lbrace \sigma > k \rbrace } \rho^{-4/3} \abs{v}{4/3}\idx \,dt.
	\end{aligned}
\end{equation}
On the other hand,  since
\begin{equation*}
	p_0 > 24/11,
\end{equation*}
applying H\"older's inequality yields
\begin{equation}\label{UL-005}
	\begin{aligned}
		& \int_0^T \int \mathbbm{1}_{\lbrace \sigma > k \rbrace } \rho^{-4/3} \abs{v}{4/3}\idx \,dt \leq \bigl(\int_0^T \int \rho v^5\idx \,dt \bigr)^{4/15} \\
		& ~~~~ \times \bigl( \int_0^T \int \mathbbm{1}_{\lbrace \sigma > k \rbrace} \rho^{-24/11} \idx \,dt \bigr)^{11/15}
		\leq C \bigl( \int_0^T \int \rho^{-p_0} \idx \,dt \bigr)^{\frac{24}{15p_0}} \\
		& ~~~~ \times \abs{\lbrace \sigma > k \rbrace }{\frac{11}{15}-\frac{24}{15p_0}} \leq C \abs{\lbrace \sigma > k \rbrace }{\frac{11}{15}-\frac{24}{15p_0}},
	\end{aligned}
\end{equation}
where we have applied \eqref{ene:basic-energy}.
Here $ \abs{\lbrace \sigma > k\rbrace }{} $ denotes the measure of the set $ \lbrace \sigma > k\rbrace \subset \Omega \times (0,T) $.
Since $ \sigma = \rho^{-1/2} $ is independent of the $ z $ variable, so is $ (\sigma - k)_+ $. We will apply the two-dimensional Sobolev embedding inequality together with \eqref{UL-004} and \eqref{UL-005},
\begin{align*}
	& \normh{(\sigma-k)_+}{L^{4}(\Omega_h \times(0,T))}^4 \leq C \int_0^T \bigl( \normh{(\sigma-k)_+}{L^2}^{2} \normh{\nablah(\sigma-k)_+}{L^2}^2 + \normh{(\sigma-k)_+}{L^2}^4 \bigr) \,dt \\
	& ~~~~ \leq C \sup_{0\leq t\leq T} \int \abs{(\sigma-k)_+}{2} \idx \times \int_0^T \int \abs{\nablah(\sigma-k)_+}{2} \idx \,dt \\
	& ~~~~ +  C \bigl( \sup_{0\leq t\leq T} \int \abs{(\sigma-k)_+}{2} \idx \bigr)^2  \leq C \abs{\lbrace \sigma > k \rbrace }{\frac{22}{15}-\frac{48}{15p_0}}.
\end{align*}
Then notice for for any positive number $ l > k $ we have
$$ \abs{\lbrace \sigma > l \rbrace}{}(l-k)^4 \leq \norm{(\sigma - k)_+}{L^4(\Omega\times(0,T))}^4 = \normh{(\sigma-k)_+}{L^4(\Omega_h\times(0,T))}^4. $$
Then by denoting $ a_k : = \abs{\lbrace \sigma > k \rbrace}{} $, we arrive at
\begin{equation}\label{UL-006}
	a_l \leq C (l-k)^{-4} a_k^{\frac{22}{15}-\frac{48}{15p_0}},
\end{equation}
with $ \frac{22}{15} - \frac{48}{15p_0} > 0 $.
Notice $ a_k \in (0,\abs{\Omega \times (0,T)}{}) $ is bounded and non-increasing with respect to $ k $. In fact, from \eqref{UL-006} it can be easily shown that, for every $ k \geq  k_0 + 1 $, one has
\begin{equation}\label{UL-008}
	0 \leq a_k \leq a_{k_0+1} \leq C(k_0 +1 - 1) ^{-4} a_1^{\frac{22}{15}-\frac{48}{15p_0}} \leq C k_0^{-4} < C.
\end{equation}
Similarly from \eqref{UL-006}, for $ l > k > 1 $, and
since
$
	p_0 \geq 24,
$
it follows that
\begin{equation}\label{UL-010}
	a_l \leq C(l-k)^{-4} a_k^{4/3}.
\end{equation}
Then the De Giorgi-type Lemma will imply the upper bound of $ \sigma $.
Indeed, apply Lemma \ref{lm:De-Giorgi} with $ g(l) = a_l, \alpha = \frac{1}{3}, \beta = 4 $, which implies that there exists $ L \in (0,\infty) $ such that $ a_L = 0 $.
This will imply that $ \sigma \leq L < \infty $. Therefore, it yields
\begin{equation}\label{UL-009}
	\rho \geq c,
\end{equation}
for some $  c \in (0, \infty) $ which may depend on $ \varepsilon $ and $ T $. This together with \eqref{UL-002} shows the upper and lower bounds of $ \rho $. This finishes the proof.
\end{pf}

\subsubsection{The global-in-time regularity estimates: for $ \eta $}
From the previous estimates \eqref{ene:basic-energy}, \eqref{ene:BD-entropy}, \eqref{ene:density-UL}, we have
\begin{equation}\label{GR-001}
	\begin{aligned}
		& \sup_{0\leq t\leq T} \bigl\lbrace \norm{v(t)}{L^2} + \norm{\nablah \eta(t)}{L^2} + \norm{\nablah \eta(t)}{L^4} \bigr\rbrace  + \norm{v}{L^5(\Omega\times (0,T))} \\
		& ~~~~ + \norm{\nabla v}{L^2(\Omega\times (0,T))} + \norm{\dz w}{L^2(\Omega\times(0,T))} + \norm{\nablah^2 \eta}{L^2(\Omega\times(0,T))} \\
		& ~~~~ + \norm{\abs{\nablah \eta}{2}\nablah^2\eta}{L^2(\Omega\times(0,T))} \leq C.
	\end{aligned}
\end{equation}
Also, from \subeqref{eq:approximating-CPE}{1}, $ \eta = \rho^{1/2}
$ satisfies
\begin{equation}\label{GR-002}
	\begin{aligned}
	& 2 \dt \eta - \varepsilon \dvh ( (1+\abs{\nablah \eta}{2}) \nablah \eta) = - \dvh (\eta v) - \eta \dz w \\
	& ~~~~ ~~~~ ~~~~ - ( v \cdot \nablah \eta  - \varepsilon \eta ^{-2p_0 - 1}).
	\end{aligned}
\end{equation}
%
%
%
%
%
%
We have the following proposition concerning the parabolic estimate of equation \eqref{GR-002}.
\begin{prop}\label{lm:density-parabolic-est}
Let $ T \in (0,\infty) $.
Consider $ \eta = \rho^{1/2}, v, w $ satisfying the estimates \eqref{ene:density-UL}, \eqref{GR-001} and $ \eta $ is the solution to \eqref{GR-002} over $[0,T] $ with $ \dz \eta \equiv 0 $. We have the following bounds,
\begin{equation}\label{ene:H^1-of-density}
	\norm{\dt \eta}{L^2(\Omega\times(0,T))} + \norm{\nablah^2 \eta}{L^2(\Omega\times(0,T))} < C.
\end{equation}
\end{prop}

\begin{pf} We only need to derive the bound of $ \dt \eta $.
Notice $ \eta $ is independent of the $ z $ variable. Integrating \eqref{GR-002} in the vertical direction yields
\begin{equation}\label{GR-003}
	2 \dt \eta - \varepsilon \dvh ( (1+\abs{\nablah \eta}{2}) \nablah \eta) = - \dvh (\eta \overline{v})- ( \overline{v} \cdot \nablah \eta - \varepsilon \eta ^{-2p_0 - 1}),
\end{equation}
where the average of $ v $ over the vertical variable is defined as,
\begin{equation}\label{GR-004}
	\overline v = \overline v (x,y,t) : = \int_0^1 v(x,y,z',t) \,dz'.
\end{equation}
\noindent Then it is straightforward to check from \eqref{ene:density-UL} and \eqref{GR-001} that
\begin{align*}
& \norm{\dvh ((1+\abs{\nablah \eta}{2})\nablah \eta)}{L^2(\Omega\times(0,T))} \leq C\norm{\deltah \eta}{L^2(\Omega\times(0,T))} \\
& ~~~~ ~~~~ + C\norm{\abs{\nablah \eta}{2}\nablah^2\eta}{L^2(\Omega\times(0,T))} \leq C,\\
& \norm{\dvh(\eta \overline v) + \overline v\cdot\nablah \eta}{L^2(\Omega\times(0,T))} \leq C \norm{\nablah v}{L^2(\Omega\times(0,T))} \\
& ~~~~ ~~~~ + C \norm{ v}{L^4(\Omega\times(0,T))}\norm{\nablah \eta }{L^4(\Omega\times(0,T))} \leq C \norm{\nablah v}{L^2(\Omega\times(0,T))}\\
& ~~~~ ~~~~ + C  \sup_{0\leq t\leq T} \norm{\nablah \eta }{L^4} \cdot \norm{ v}{L^5(\Omega\times(0,T))} \leq C,
\end{align*}
where we have applied H\"older's inequality and noticed the fact that for any $ \varphi \in L^p(\Omega)$, $ p>1 $, the vertical average $ \overline \varphi : = \int_0^1 \varphi(\cdot,z') \,dz' $ satisfies, 
\begin{equation}\label{GR-007}
\normh{\overline \varphi}{L^p} \leq \int_0^1 \normh{\varphi}{L^p} (z',t)\,dz' \leq \bigl(\int_0^1 \normh{\varphi}{L^p}^p (z',t)\,dz' \bigr)^{1/p} = \norm{\varphi}{L^p},
\end{equation}
by employing the Minkowski and H\"older's inequalities. Then \eqref{ene:H^1-of-density} follows directly from the equation \eqref{GR-002}.
\end{pf}

Before we move on to the estimates for $ v $, we will derive the a $ L^p $ estimate of $ \nablah \eta $ first. This will be useful in deriving the estimates on $ v $.
\begin{prop}\label{lm:density-Lp-est}
	Under the same assumptions as in Proposition \ref{lm:density-parabolic-est}, we have the following inequality,
	\begin{equation}\label{GR-009}
		\norm{\nablah \rho}{L^{3p}(\Omega \times (0,T))}^{3p} \leq C_p \norm{\nablah \eta}{L^{3p}(\Omega \times (0,T))}^{3p} \leq C_p + C_p \norm{v}{L^p(\Omega\times(0,T))}^{2p},
	\end{equation}
	for any $ p > 4/3 $.
\end{prop}
\begin{pf}
To establish this inequality,
we first write \eqref{GR-002} in the divergence form. In order to do so, define $ \xi $ to be the solution to the following elliptic equation in $ \Omega_h $.
\begin{equation}\label{GR-005}
	\begin{cases}
		\deltah \xi = \overline{v} \cdot\nablah \eta - \varepsilon \eta^{-2p_0 - 1} - \int (\overline{v} \cdot\nablah \eta - \varepsilon \eta^{-2p_0 - 1}) \idxh, \\
		\int \xi \idxh = 0.
	\end{cases}
\end{equation}
Here
\begin{align*}
	& \abs{\int (\overline{v} \cdot\nablah \eta - \varepsilon \eta^{-2p_0 - 1}) \idxh }{} \leq C \normh{\overline v}{L^2} \normh{\nablah \eta}{L^2} + C \\
	& ~~~~ ~~~~ \leq C \norm{v}{L^2} \norm{\nablah \eta}{L^2} + C < C,
\end{align*}
where we have applied \eqref{GR-007}.
Then,  thanks to the Sobolev embedding inequality, and since $ \int \nablah \xi \idxh = 0 $,
$ \xi $ will satisfy, for any $ p > 2 $, 
\begin{equation}\label{GR-008}
	\begin{aligned}
		& \normh{\nablah \xi}{L^p} \leq C\normh{\nablah^2 \xi}{L^{2p/(2+p)}} \leq C \normh{\deltah \xi}{L^{2p/(2+p)}} \leq C \normh{\bar v \cdot\nablah \eta }{L^{2p/(2+p)}}\\
		& ~~~~ + C \leq C \normh{\overline v}{L^p} \normh{\nablah \eta}{L^2}+ C \leq C \norm{\overline v}{L^p} \norm{\nablah \eta}{L^2}+ C \\
		& ~~~~ \leq C \norm{v}{L^p} + C.
	\end{aligned}
\end{equation}
Then \eqref{GR-003} can be written as
\begin{align*}
	& 2\dt \eta - \varepsilon \dvh((1+\abs{\nablah \eta}{2}) \nablah \eta) = - \dvh(\eta \overline v + \nablah \xi) \\
	& ~~~~ ~~~~ - \int (\overline{v} \cdot\nablah \eta - \varepsilon\eta^{-2p_0-1}) \idxh.
\end{align*}
Or, by denoting $$ \hat{\eta} = \eta + \dfrac{1}{2} \int_0^t \int_{\Omega_h} (\overline{v} \cdot\nablah \eta - \varepsilon\eta^{-2p_0-1}) \idxh\,dt, $$
the following equation holds
\begin{equation}\label{GR-006}
	2 \dt \hat{\eta} - \varepsilon\dvh (\abs{\nablah \hat\eta}{2} \nablah \hat{\eta}) = \dvh (\varepsilon \nablah \eta - \eta \overline v - \nablah \xi) =: \dvh\hat{f},
\end{equation}
with $ \hat{\eta}(x,y,0) = \eta(x,y,0) = \rho_{\varepsilon,0}^{1/2}(x,y)$. Then by applying the $ L^p $ estimate 
 in \cite[Theorem 1]{Acerbi2007} for the $p$-Laplacian equation \eqref{GR-006},
for $ p > 4/3 $, we have
\begin{align*}
	\normh{\nablah \hat{\eta}}{L^{3p}(\Omega_h\times(0,T))}^4 \leq C_p\bigl( \normh{\nablah \hat{\eta}}{L^4(\Omega_h\times(0,T))}^4 + \normh{\hat f}{L^{p}(\Omega_h\times(0,T))}^{4/3} + 1 \bigr)^2.
\end{align*}
Therefore, by employing \eqref{ene:density-UL}, \eqref{GR-001}, \eqref{GR-007}, \eqref{GR-008},
\begin{align*}
	& \norm{\nablah \eta}{L^{3p}(\Omega \times (0,T))}^{3p} = \normh{\nablah \hat \eta}{L^{3p}(\Omega_h\times(0,T))}^{3p} \leq C_p \bigl( \norm{\nablah \eta}{L^4(\Omega\times(0,T))}^{3p} \\
	& ~~~~ ~~~~ + \normh{\hat{f}}{L^p(\Omega_h\times(0,T))}^{p} + 1 \bigr)^2  \leq C_p \bigl( 1 + \normh{\nablah \eta}{L^p(\Omega_h\times(0,T))}^p \\
	& ~~~~ ~~~~ + \normh{\overline v}{L^p(\Omega_h\times(0,T))}^p + \normh{\nablah \xi}{L^p(\Omega_h\times(0,T))}^p \bigr)^{2}\\
	& ~~~~ \leq C_p + C_p \norm{v}{L^p(\Omega\times(0,T))}^{2p} + \dfrac{1}{2} \norm{\nablah \eta}{L^{3p}(\Omega\times(0,T))}^{3p},
\end{align*}
where in the last inequality, H\"older's and Young's inequalities have been applied as follows,
\begin{align*}
	& \normh{\nablah \eta}{L^p(\Omega_h\times(0,T))}^{2p} \leq \norm{\nablah \eta}{L^p(\Omega\times(0,T))}^{2p} \leq \norm{\nablah \eta}{L^{3p}(\Omega\times(0,T))}^{2p} \\
	& ~~~~ ~~~~ \leq C_p + \dfrac{1}{2}\norm{\nablah \eta}{L^{3p}(\Omega\times(0,T))}^{3p}.
\end{align*}
This together with \eqref{ene:density-UL} completes the proof.
\end{pf}
%
%

\subsubsection{The global-in-time regularity estimates: for $ v $}\label{sec:global-approx-priori-v}

Next, we will derive the estimates of $ v $. To do so, \subeqref{eq:approximating-CPE}{2} can be written as
\begin{equation}\label{GR-010}
	\dt v - \bigl( \dfrac{1}{2} + \sqrt{\varepsilon} \bigr) \deltah v - \partial_{zz} v - \dfrac{1}{2} \nablah \dvh v = g,
\end{equation}
where
\begin{equation}\label{approximate:v-source}
	\begin{aligned}
		& g : = ( \dfrac{1}{2} + \sqrt{\varepsilon}) \nablah \log \rho \cdot\nablah v + \dfrac{1}{2} \nablah v \cdot\nablah \log  \rho - v \cdot\nablah v - w \dz v \\
		& ~~~~ - \rho^{-1} \nablah \rho^{\gamma}   + \varepsilon \rho^{-1/2} \abs{\nablah \rho^{1/2}}{2} \nablah \rho^{1/2} \cdot\nablah v - \varepsilon \rho^{-p_0-1} v - \varepsilon\abs{v}{3}v.
	\end{aligned}
\end{equation}
Let $ \overline v $ be the vertical average defined in \eqref{GR-004} and
\begin{equation} \label{GR-017}
	\widetilde{v} = v - \overline v.
\end{equation}
Then we have, for any $ p > 1 $, the following inequalities as in \eqref{GR-007},
\begin{equation}\label{GR-014}
	\begin{gathered}
		\norm{\overline v}{L^p} = \normh{\overline v}{L^p} \leq \norm{v}{L^p}, ~~
		\norm{\widetilde{v}}{L^p} \leq \norm{v}{L^p} + \norm{\overline v}{L^p} \leq 2 \norm{v}{L^p}, \\
		\norm{\nablah \overline v}{L^p} = \normh{\nablah \overline v}{L^p} \leq \norm{\nablah v}{L^p} , ~~ \norm{\nablah^2 \overline v}{L^p} = \normh{\nablah^2 \overline v}{L^p} \leq \norm{\nablah^2 v}{L^p}, \\
		 \norm{\nablah \widetilde{v}}{L^p} \leq 2 \norm{\nablah v}{L^p},~ \norm{\nablah^2 \widetilde{v}}{L^p} \leq 2 \norm{\nablah^2 v}{L^p}.
	\end{gathered}
\end{equation}
Moreover, after integrating \subeqref{eq:approximating-CPE}{1} in the vertical direction and taking the difference, we have the following equations
\begin{equation*}
	\begin{gathered}
		\dt \rho + \dvh(\rho \overline v) = G(\rho), \\
		\dvh(\rho \widetilde v) + \dz(\rho w) = 0.
	\end{gathered}
\end{equation*}
Then we have
\begin{equation}\label{GR-011}
	w = - \rho^{-1} \int_0^z \dvh (\rho \widetilde v) \,dz'.
\end{equation}
In the rest of this section, we will show the following:
\begin{prop}\label{lm:velocity-parabolic-est}
	Under the same assumptions as in Proposition \ref{lm:density-parabolic-est}, let $ v $ be the solution to the parabolic equation \eqref{GR-010} and $ w $ be given by \eqref{GR-011}. Then the following estimate holds
	\begin{equation}\label{ene:H^1-of-v}
		\begin{gathered}
		\sup_{0\leq t\leq T} \norm{\nabla v(t)}{L^2} + \norm{\nabla^2 v}{L^2(\Omega\times(0,T))} + \norm{v_t}{L^2(\Omega\times(0,T))} < C, \\
		\norm{v}{L^{10}(\Omega\times(0,T))} + \norm{\nabla v}{L^{10/3}(\Omega\times(0,T))} + \norm{\nablah \rho}{L^{30}(\Omega\times(0,T))} < C.
		\end{gathered}
	\end{equation}
\end{prop}
\begin{pf}
We will establish the vertical derivative estimate for $ v $ first. Take the inner product of \eqref{GR-010} with $ - 2 \partial_{zz} v $ and integrate the resultant in the spatial variable.
After integration by parts, one obtains,
\begin{equation}\label{GR-012}
	\begin{aligned}
		& \dfrac{d}{dt} \int \abs{\dz v}{2} \idx + \int \biggl( (1+2\sqrt{\varepsilon}) \abs{\nablah\dz v}{2} + 2 \abs{\partial_{zz} v}{2} + \abs{\dz \dvh v}{2} \biggr) \idx \\
		& ~~~~ = \int \biggl( (1+2\sqrt{\varepsilon}) (\nablah \log \rho \cdot \nablah) \dz v \cdot\dz v \\
		& ~~~~ ~~~~ + (\nablah \dz v \cdot\nablah) \log \rho \cdot \dz v \biggr) \idx - 2 \int \dz(v \cdot\nablah v + w \dz v) \cdot \dz v\idx \\
		& ~~~~ ~~~~ +  2 \varepsilon\int ( \rho^{-1/2} \abs{\nablah \rho^{1/2}}{2} \nablah \rho^{1/2} \cdot\nablah ) \dz v \cdot\dz v \idx \\
		& ~~~~ ~~~~	- 2 \varepsilon \int \biggl( \rho^{-p_0-1} \abs{\dz v}{2} + \abs{v}{3} \abs{\dz v}{2} + \dfrac{3}{4} \abs{v}{} \abs{\dz \abs{v}{2}}{2} \biggr) \idx \\
		& ~~~~ ~~~~ =: I_1 + I_2 +I_3 + I_4.
	\end{aligned}
\end{equation}
As a consequence of \eqref{ene:density-UL}, we have the following estimates of the right-hand side of the above equation,
\begin{align*}
	& I_1 \leq C \norm{\nablah \rho}{L^{6}}  \norm{\dz v}{L^{3}}\norm{\nablah \dz v}{L^{2}} \\
	& ~~~~ \leq C \norm{\nablah \rho}{L^6} (\norm{\dz v}{L^2}^{1/2}\norm{\nabla \dz v}{L^2}^{3/2} + \norm{\dz v}{L^2} \norm{\nablah \dz v}{L^2})
	\\
	& ~~~~
	\leq \dfrac{1}{8} \norm{\nabla \dz v}{2}^2  + C( \norm{\nablah \rho}{L^6}^4 + \norm{\nablah \rho}{L^6}^2 )\norm{\dz v}{L^2}^2,\\
	& I_3 \leq C \norm{\nablah \rho}{L^{15}}^3 \norm{\dz v}{L^{10/3}} \norm{\nablah \dz v}{L^{2}} \\
	& ~~~~ ~~~~ \leq C \norm{\nablah \rho}{L^{15}}^3 (\norm{\dz v}{L^2}^{2/5} \norm{\nabla \dz v}{L^2}^{8/5} + \norm{\dz v}{L^2} \norm{\nablah \dz v}{L^2})
	\\
	& ~~~~
	\leq \dfrac{1}{8} \norm{\nabla \dz v}{L^2}^2 + C (\norm{\nablah \rho}{L^{15}}^{15} + \norm{\nablah \rho}{L^{15}}^6) \norm{\dz v}{L^2}^{2},\\
	& I_4 \leq 0,
\end{align*}
where we have applied the Gagliardo-Nirenberg inequality as follows,
\begin{gather*}
	\norm{\dz v}{L^3} \leq C \norm{\dz v}{L^2}^{1/2} \norm{\nabla \dz v}{L^2}^{1/2} + C \norm{\dz v}{L^2},  \\
	\norm{\dz v}{L^{10/3}} \leq C \norm{\dz v}{L^2}^{2/5} \norm{\nabla \dz v}{L^2}^{3/5} + C \norm{\dz v}{L^2},
\end{gather*}
noticing that $ \dz v = 0 $ at $ z =0, 1 $.
On the other hand, directly applying integration by parts implies,
\begin{align*}
	& I_2 = -2 \int (\dz v \cdot\nablah v + \dz w \dz v ) \cdot\dz v \idx - 2 \int (v \cdot\nablah \dz v + w \dz^2 v) \cdot\dz v \idx \\
	& ~~~~ = 2 \int \biggl( \dz \dvh v  (v \cdot\dz v) + (\dz v \cdot\nablah) \dz v \cdot v - (v \cdot\nablah) \dz v \cdot \dz v \biggr) \idx \\
	& ~~~~ ~~~~  - \int \dz w \abs{\dz v}{2} \idx.
\end{align*}
After substituting \eqref{GR-011}, we have
\begin{align*}
	& - \int \dz w \abs{\dz v}{2} \idx = \int \rho^{-1} \dvh(\rho \widetilde v)  \abs{\dz v}{2} \idx = - 2 \int \widetilde{v} \cdot \nablah \dz v \cdot \dz v \idx \\
	& ~~~~ + \int  ( \rho^{-1} \widetilde{v} \cdot \nablah \rho ) \abs{\dz v}{2} \idx.
\end{align*}
Therefore, we have after employing \eqref{ene:density-UL}
\begin{equation*}
	I_2 \leq C \int \abs{\nablah \rho}{} \abs{\widetilde{v}}{} \abs{\dz v}{2} \idx  + C \int (\abs{v}{} + \abs{\widetilde{v}}{}) \abs{\dz v}{} \abs{\nabla\dz v}{}\idx =: I_2' + I_2''.
\end{equation*}
As before, applying \eqref{GR-014}, the Sobolev embedding and Young's inequalities then yield,
\begin{align*}
	& I_2' \leq C \norm{\nablah \rho}{L^{15}} \norm{\widetilde{ v }}{L^{5}} \norm{\dz v}{L^{30/11}}^2 \leq C\norm{\nablah \rho}{L^{15}} \norm{v}{L^5} (\norm{\dz v}{L^2}^{6/5} \\
	& ~~~~ \times \norm{\nabla \dz v}{L^2}^{4/5} + \norm{\dz v}{L^2}^2)
	 \leq \dfrac{1}{8} \norm{\nablah \dz v}{L^2}^2 \\
	 & ~~~~ ~~~~ + C( \norm{\nablah \rho}{L^{15}}^{5/2} + \norm{\nablah \rho}{L^{15}}^{5/4} + \norm{v}{L^5}^5  )\norm{\dz v}{L^2}^2,\\
	& I_2'' \leq C( \norm{v}{L^{5}} + \norm{\widetilde{v}}{L^{5}} ) \norm{\dz v}{L^{10/3}} \norm{\nabla \dz v}{L^2} \leq C \norm{v}{L^5} (\norm{\dz v}{L^2}^{2/5} \\
	& ~~~~ \times \norm{\nabla\dz v}{L^2}^{8/5} + \norm{\dz v}{L^2} \norm{\nabla\dz v}{L^2}  ) \leq \dfrac{1}{8} \norm{\nabla \dz v}{L^2}^2 \\
	& ~~~~ ~~~~ + C (\norm{v}{L^5}^{5} + \norm{v}{L^5}^2) \norm{\dz v}{L^2}^2.
\end{align*}
After summing up these estimates and applying H\"older's and Young's inequalities, \eqref{GR-012} implies
\begin{equation}\label{GR-013}
	\begin{aligned}
		& \dfrac{d}{dt} \norm{\dz v}{L^2}^2 + (\dfrac{1}{2} + 2 \sqrt \varepsilon) \norm{\nablah \dz v}{L^2}^2 + \dfrac{3}{2} \norm{\partial_{zz}v}{2}^2 \\
		& ~~~~ \leq C (\norm{\nablah \rho}{L^{15}}^{15} + \norm{v}{L^5}^5 + 1) \norm{\dz v}{L^2}^2,
	\end{aligned}
\end{equation}
where from \eqref{GR-001} and \eqref{GR-009}, with $ p = 5 $, we have
\begin{equation}\label{GR-016}
	\int_0^T \biggl( \norm{\nablah \rho}{L^{15}}^{15} + \norm{v}{L^5}^5 + 1 \biggr) \,dt < C.
\end{equation}
Therefore, after applying Gr\"onwall's inequality to \eqref{GR-013}, we have
\begin{equation}\label{ene:H1-of-dz-v}
	\sup_{0\leq t\leq T} \norm{\dz v(t)}{L^2} + \norm{\nabla \dz v}{L^2(\Omega \times (0,T))} < C.
\end{equation}
Next, we establish the horizontal derivative estimates for $ v $. Take the inner product of \eqref{GR-010} with $ - 2 \deltah v $ and integrate the resultant in the spatial variable. After integration by parts, one has,
\begin{equation}\label{GR-015}
	\begin{aligned}
		& \dfrac{d}{dt} \int \abs{\nablah v}{2}\idx + \int \biggl( (1+2\sqrt{\varepsilon}) \abs{\deltah v}{2} + 2 \abs{\nablah \dz v}{2} + \abs{\nablah \dvh v }{2} \biggr)  \idx\\
		& ~~~~ = - \int \biggl( ( 1 + 2\sqrt{\varepsilon}) (\nablah \log \rho \cdot\nablah) v \cdot\deltah v + ( \deltah v \cdot \nablah ) v \cdot\nablah \log  \rho \biggr) \idx \\
		& ~~~~ ~~~~  + 2 \int  (v \cdot\nablah v + w \dz v) \cdot\deltah v \idx + 2 \int \rho^{-1} \nablah \rho^{\gamma} \cdot \deltah v \idx \\
		& ~~~~ ~~~~ - 2 \varepsilon \int ( \rho^{-1/2} \abs{\nablah \rho^{1/2}}{2} \nablah \rho^{1/2} \cdot\nablah ) v \cdot \deltah v \idx \\
		& ~~~~ ~~~~ + 2 \varepsilon \int \rho^{-p_0-1} v \cdot\deltah v \idx  - 2 \varepsilon \int \abs{v}{3} \abs{\nablah v}{2} + \dfrac{3}{4} \abs{v}{} \abs{\nablah \abs{v}{2}}{2} \idx 
		\\ & ~~~~ =: I_5 + I_6 + I_7 + I_8 + I_9 + I_{10}.
	\end{aligned}
\end{equation}
With \eqref{ene:density-UL}, we will have the following estimates of the right-hand side of the above terms,
\begin{align*}
	& I_5 \leq C \norm{\nablah \rho}{L^{6}} \norm{\nablah v}{L^{3}} \norm{\deltah v}{L^2} \leq C \norm{\nablah \rho}{L^6} (\norm{\nablah v}{L^2}^{1/2} \\
	& ~~~~ ~~~~ \times \norm{\nabla^2 v}{L^2}^{3/2} + \norm{\nablah v}{L^2} \norm{\deltah v}{L^2} ) \leq \dfrac{1}{8} \norm{\nabla^2 v}{L^2}^2 \\
	& ~~~~ ~~~~ + C ( \norm{\nablah \rho}{L^6}^4 + \norm{\nablah \rho}{L^6}^2 ) \norm{\nablah v}{L^2}^2,  \\
	& I_7 \leq C \norm{\nablah \rho}{L^2} \norm{\deltah v}{L^2} \leq \dfrac{1}{8} \norm{\deltah v}{L^2}^2 + C \norm{\nablah \rho}{L^2}^2, \\
	& I_8 \leq C \norm{\nablah \rho}{L^{15}}^3 \norm{\nablah v}{L^{10/3}} \norm{\deltah v}{L^2} \leq C \norm{\nablah \rho}{L^{15}}^3(\norm{\nablah v}{L^2}^{2/5} \\
	& ~~~~ ~~~~ \times \norm{\nabla^2 v}{L^2}^{8/5}  + \norm{\nablah v}{L^2}\norm{\deltah v}{L^2}) \leq \dfrac{1}{8} \norm{\nabla^2 v}{L^2}^2 \\
	& ~~~~ ~~~~ + C (\norm{\nablah \rho}{L^{15}}^{15} + \norm{\nablah \rho}{L^{15}}^6) \norm{\nablah v}{L^2}^2,  \\
	& I_9 \leq \dfrac{1}{8} \norm{\deltah v}{L^2}^2 + C \norm{v}{L^2}^2  \\
	& I_{10} \leq 0,
\end{align*}
where we have applied the following interpolation inequalities
\begin{gather*}
	\norm{\nablah v}{L^3} \leq C \norm{\nablah v}{L^2}^{1/2} \norm{\nabla^2 v}{L^2}^{1/2} + C\norm{\nablah v}{L^2}, \\
	\norm{\nablah v}{L^{10/3}} \leq C \norm{\nablah v}{L^2}^{2/5} \norm{\nabla^2 v}{L^2}^{3/5} + C \norm{\nablah v}{L^2}.
\end{gather*}
One the other hand, after substituting \eqref{GR-011},
\begin{align*}
	& I_6 = 2 \int \biggl\lbrack ( v \cdot\nablah) v \cdot\deltah v -  \rho^{-1} \int_0^z \biggl( \dvh(\rho \widetilde v)(\cdot,z',t) \biggr) \,dz'  ( \dz v \cdot\deltah v ) \biggr\rbrack \idx \\
	& ~~~~ = 2 \int ( v \cdot\nablah ) v \cdot \deltah v \idx - 2 \int \biggl\lbrack \biggl( \int_0^z \dvh \widetilde v \,dz' \biggr) ( \dz v\cdot\deltah v ) \biggr\rbrack \idx \\
	& ~~~~ ~~~~ - 2 \int \biggl\lbrack \biggl( \int_0^z ( \widetilde v \cdot\nablah \log\rho ) \,dz' \biggr) ( \dz v \cdot\deltah v) \biggr\rbrack \idx =: I_6' + I_6''+ I_6''',\\
	& I_6' \leq C \norm{v}{L^{5}}\norm{\nablah v}{L^{10/3}} \norm{\deltah v}{L^2} \leq \dfrac{1}{8}\norm{\nabla^2 v}{L^2}^2 \\
	& ~~~~ ~~~~ + C(\norm{v}{L^5}^5 + \norm{v}{L^5}^2) \norm{\nablah v}{L^2}^2.
\end{align*}
Meanwhile, we shall estimate $ I_6'', I_6''' $ by making use of the stratification structure. Indeed, we will apply the Minkowski and H\"older's inequalities as follows,
\begin{align*}
	& I_6''  = - 2 \int_0^1 \biggl\lbrack \int_{\Omega_h} \biggl\lbrack  \biggl( \int_0^z \dvh \widetilde{v}(\cdot,z',t) \,dz' \biggr) (\dz v\cdot\deltah v)(\cdot, z,t) \biggr\rbrack  \idxh  \biggr\rbrack \,dz \\
	& ~~~~ \leq C \biggl( \int_0^1 \normh{\nablah \widetilde v}{L^4}\,dz' \biggr) \times \biggl( \int_0^1 \normh{\dz v}{L^4} \normh{\deltah v}{L^2} \,dz \biggr) \\
	& ~~~~ \leq C \biggl( \int_0^1 \normh{\nablah \widetilde{v}}{L^2}^{1/2} \normh{\nablah^2 \widetilde{v}}{L^2}^{1/2} \,dz' \biggr) \times \biggl( \int_0^1 (\normh{\dz v}{L^2}^{1/2} \normh{\nablah \dz v}{L^2}^{1/2} + \normh{\dz v}{L^2}) \\
	& ~~~~ ~~~~ \cdot \normh{\deltah v}{L^2} \,dz \biggr) \leq C \norm{\nablah \widetilde{v}}{L^2}^{1/2} \norm{\nablah^2 \widetilde{v}}{L^2}^{1/2}(\norm{\dz v}{L^2}^{1/2} \norm{\nablah \dz v}{L^2}^{1/2} \\
	& ~~~~ ~~~~ \times \norm{\deltah v}{L^2} + \norm{\dz v}{L^2} \norm{\deltah v}{L^2} )  \leq \dfrac{1}{8}\norm{\nablah^2 v}{L^2}^2 \\
	& ~~~~ ~~~~ + C (\norm{\dz v}{L^2}^2 \norm{\nablah \dz v}{L^2}^2 + \norm{\dz v}{L^2}^4 )\norm{\nablah v}{L^2}^2,\\
	& I_6'''
	\leq C \biggl( \int_0^1 \normh{\widetilde{v}}{L^{6}} \normh{\nablah \rho}{L^{12}}\,dz'\biggr) \times \biggl( \int_0^1 \normh{\dz v}{L^{4}} \normh{\deltah v}{L^2} \,dz \biggr) \\
	& ~~~~ \leq C \biggl( \int_0^1 (\normh{\widetilde{v}}{L^5}^{5/6} \normh{\nablah \widetilde{v}}{L^2}^{1/6} + \normh{\widetilde v}{L^2}) \normh{\nablah \rho}{L^{12}} \,dz' \biggr) \times \biggl( \int_0^1 (\normh{\dz v}{L^2}^{1/2} \normh{\nablah \dz v}{L^2}^{1/2} \\
	& ~~~~ ~~~~ + \normh{\dz v}{L^2}) \normh{\deltah v}{L^2}  \,dz  \biggr)
	 \leq \dfrac{1}{8} \norm{\deltah v}{L^2}^2 + C\norm{v}{L^5}^{5} \norm{\nablah v}{L^2}^{} \\
	 & ~~~~ ~~~~ + C \norm{\nablah \rho}{L^{12}}^{12} + C \norm{\dz v}{L^2}^{2} \norm{\nablah \dz v}{L^2}^{2} + C \norm{v}{L^2}^6 + C \norm{\dz v}{L^2}^{4},
\end{align*}
where we have employed \eqref{GR-014} and the two-dimensional Sobolev embedding inequality. Therefore, by summing up all these inequalities and employing \eqref{GR-001} and \eqref{ene:H1-of-dz-v}, \eqref{GR-015} implies
\begin{align*}
	& \dfrac{d}{dt} \norm{\nablah v}{L^2}^2 + \dfrac{1}{8} \norm{\nablah^2 v}{L^2}^2 + \dfrac{9}{8} \norm{\nablah \dz v}{L^2}^2 \\
	& ~~~~ \leq C (\norm{\nablah \rho}{L^{15}}^{15} + \norm{v}{L^5}^5 + \norm{\nablah \dz v}{L^2}^2 + 1) \norm{\nablah v}{L^2}^2 \\
	& ~~~~ ~~~~ + C \norm{\nablah \rho}{L^{15}}^{15} + C \norm{v}{L^5}^5 + C \norm{\nabla\dz v}{L^2}^2 + C.
\end{align*}
Therefore, \eqref{GR-016}, \eqref{ene:H1-of-dz-v} and the Gr\"onwall's inequality imply
\begin{equation}\label{ene:H^1-of-v-1}
	\sup_{0\leq t\leq T} \norm{\nabla v(t)}{L^2} + \norm{\nabla^2 v}{L^2(\Omega\times(0,T))} < C.
\end{equation}
For simplicity, in the following, we denote, for every $ p > 1 $ the norm in space-time $ \norm{\cdot}{L^{p}(\Omega\times (0,T))} $ as $ \norm{\cdot}{L^p_{tx}} $.
The Sobolev embedding inequality and \eqref{ene:H^1-of-v-1} then yield,
\begin{equation}\label{GR-018}
\begin{aligned}
& \norm{v}{L^{10}_{tx}} + \norm{\nabla v}{L^{10/3}_{tx}} \leq (\sup_{0\leq t\leq T}\norm{\nabla v}{L^2})^{4/5} \norm{\nabla^2 v}{L^2_{tx}}^{1/5} \\
& ~~~~ + \norm{v}{L_{tx}^2} + (\sup_{0\leq t\leq T} \norm{\nabla v}{L^2})^{2/5} \norm{\nabla^2 v}{L^2_{tx}}^{3/5} < C.
\end{aligned}
\end{equation}
Then by taking $ p = 10 $ in \eqref{GR-009}, we have
\begin{equation}\label{GR-019}
\norm{\nablah \rho}{L^{30}_{tx}} \leq C \norm{\nablah \eta}{L^{30}_{tx}} \leq C_p + C_p \norm{v}{L^{10}_{tx}}^{2/3} < C.
\end{equation}
Here we have used \eqref{ene:density-UL}.
Now, we can obtain the integrability of $ g $ in \eqref{approximate:v-source}. Indeed, the interpolation and Young's inequalities imply,
\begin{align*}
& \norm{\nablah \rho \cdot\nablah v}{L^{3}_{tx}} + \norm{\nablah v \cdot\nablah  \rho}{L^{3}_{tx}} \leq C \norm{\nablah \rho}{L^{30}_{tx}} \norm{\nabla v}{L^{10/3}_{tx}}, \\
& \norm{v \cdot\nablah v}{L^{5/2}_{tx}} \leq C \norm{v}{L^{10}_{tx}}\norm{\nabla v}{L^{10/3}_{tx}}, \\
& \norm{\abs{\nablah \rho}{2} \nablah \rho \cdot\nablah v}{L^{5/2}_{tx}} \leq C \norm{\nablah \rho}{L^{30}_{tx}}^{3} \norm{\nabla v}{L^{10/3}},\\
& \norm{\abs{v}{3}v}{L^{5/2}_{tx}} \leq C \norm{v}{L^{10}_{tx}}^4.
\end{align*}
Moreover, from \eqref{ene:density-UL}, \eqref{GR-014}, \eqref{GR-011}, we have
\begin{align*}
& \norm{w \dz v}{L^{2}_{tx}} = \norm{\rho^{-1} \biggl( \int_0^z \dvh (\rho \widetilde v) \,dz' \biggr) \dz v}{L^{2}_{tx}} \\
& ~~~~ \leq C\norm{\biggl( \int_0^z \dvh \widetilde v \,dz' \biggr) \dz v}{L^{2}_{tx}}
 + C \norm{\biggl( \int_0^z \widetilde v \cdot\nablah \rho \,dz' \biggr) \dz v}{L^{2}_{tx}}.
\end{align*}
Observe that,
\begin{align*}
& \norm{\biggl( \int_0^z \dvh \widetilde v \,dz' \biggr) \dz v}{L^{2}_{tx}}^2 =\int_{0}^T \biggl\lbrack \int_0^1 \biggl( \int_{\Omega_h}\abs{ \bigl( \int_0^z \dvh \widetilde v \,dz' \bigr) \dz v}{2} \idxh \biggr) \,dz \biggr\rbrack \,dt  \\
& ~~~~ \leq C  \int_0^T \biggl\lbrack \int_0^1 \biggl( \normh{\int_0^z\dvh \widetilde v \,dz' }{L^{4}_{}}^{2} \normh{\dz v}{L^{4}_{}}^{2} \biggr)  \,dz\biggr\rbrack \,dt \\
& ~~~~ \leq C \int_0^T \biggl\lbrack  ( \int_0^1 \normh{\nablah^2 \widetilde v}{L^2}^{1/2} \normh{\nablah \widetilde v}{L^2}^{1/2}  \,dz')^2  \int_0^1 (\normh{\nablah \dz v}{L^2} \normh{\dz v}{L^2} \\
& ~~~~ ~~~~ + \normh{\dz v}{L^2}^2) \,dz \biggr\rbrack \,dt  \leq C \int_0^T \biggl\lbrack  \norm{\nablah^2 v}{L^2} \norm{\nablah v}{L^2} \times (\norm{\nablah \dz v}{L^2} \\
& ~~~~ ~~~~ \cdot \norm{\dz v}{L^2} + \norm{\dz v}{L^2}^2 ) \biggr\rbrack \,dt  \leq C\sup_{0\leq t\leq T} \norm{\nabla v(t)}{L^2}^2 \times \norm{\nabla^2 v}{L^2_{tx}}^2 \\
& ~~~~ ~~~~ + C\sup_{0\leq t\leq T} \norm{\nabla v(t)}{L^2}^4,\\
& \text{and that}\\
&  \norm{\biggl( \int_0^z \widetilde v \cdot\nablah \rho \,dz' \biggr) \dz v}{L^{2}_{tx}}^2 = \int_0^T \int_{0}^{1} \biggl( \int_{\Omega_h} \abs{ \bigl( \int_0^z \widetilde v \cdot\nablah \rho \,dz' \bigr) \dz v}{2} \idxh \biggr) \,dz \,dt\\
& ~~~~ \leq C \int_0^T \biggl\lbrack  \int_0^1 \normh{\int_0^z \widetilde v \cdot\nablah \rho \,dz' }{L^{4}}^2 \normh{\dz v}{L^{4}}^2\,dz \biggr\rbrack \,dt \\
& ~~~~ \leq C \int_0^T \biggl\lbrack (\int_0^1 \normh{\widetilde v}{L^{60/13}} \normh{\nablah \rho}{L^{30}} \,dz' )^2 \int_0^1 \biggl( \normh{\dz v}{L^2} \normh{\nablah \dz v}{L^{2}} \\
& ~~~~ ~~~~ + \normh{\dz v}{L^2}^2 \biggr) \,dz \biggr\rbrack \,dt \leq C \int_0^T \biggl\lbrack \norm{v}{L^{60/13}}^2 \norm{\nablah \rho}{L^{30}}^2  \times (\norm{\dz v}{L^2} \\
& ~~~~ ~~~~ \cdot \norm{\nablah \dz v }{L^2}+ \norm{\dz v}{L^2}^2) \biggr\rbrack \,dt \leq C \int_0^T \biggl\lbrack (\norm{v}{L^2}^{3/10}\norm{\nabla v}{L^2}^{17/10} \\
& ~~~~ ~~~~ + \norm{v}{L^2}^2 )\norm{\nablah \rho}{L^{30}}^2  (\norm{\dz v}{L^2}\norm{\nablah \dz v }{L^2}+ \norm{\dz v}{L^2}^2) \biggr\rbrack \,dt \\
& ~~~~ \leq C \norm{\nabla^2 v}{L^2_{tx}}^2 + C \norm{\nablah \rho}{L^{30}_{tx}}^{30} + C \sup_{0\leq t\leq T} ( \norm{v(t)}{L^2} + \norm{\nabla v(t)}{L^2} )^{90/13}\\
& ~~~~ + C.
\end{align*}
Therefore, we have $ \norm{g}{L^2_{tx}} < C $, thanks to \eqref{GR-001},
\eqref{ene:H^1-of-v}, \eqref{GR-018} and \eqref{GR-019}.
In particular, \eqref{GR-010} together with \eqref{ene:H^1-of-v-1} implies,
\begin{equation}\label{ene:H^1-of-v-2}
\norm{v_t}{L^2(\Omega\times(0,T))} + \norm{\nabla^2 v}{L^2(\Omega\times(0,T))} \leq C + C \norm{g}{L^2(\Omega\times(0,T))} < C.
\end{equation}
With \eqref{ene:H^1-of-v-1}, \eqref{GR-018} and \eqref{GR-019} at hand, this finishes the proof of the proposition.
\end{pf}

\subsection{The global existence of approximating solutions}\label{sec:aprx-exist}

Now we will establish the existence of the approximating solutions satisfying Definition \ref{def:aprxm-sols}, for any fixed $ \varepsilon > 0 $, via a modified Galerkin approximation scheme. This is equivalent to show the existence of strong solutions to the system consisting of \eqref{GR-003}, \eqref{GR-011} and the momentum equation with $ w $ given by \eqref{GR-011}. By omitting, for simplicity, the approximating parameter $ \varepsilon $ (taking, e.g., $ \varepsilon = 1 $), we arrive at the following parabolic system:
\begin{equation}\label{eq:GA}
	\begin{cases}
		2 \dt \eta - \dvh((1+\abs{\nablah \eta}{2})\nablah\eta) = - \dvh (\eta \overline v) - \overline v \cdot \nablah \eta + \eta^{-2p_0-1}& \text{in} ~ \Omega, \\
		 \dt (\eta^2 v ) + \dvh (\eta^2 v \otimes v) + \dz(\eta^2 w v) + \nablah \eta^{2\gamma} - \dvh(\eta^2\mathcal D (v)) \\
		~~~~ - \dz (\eta^2 \dz v)
		- \dvh(\eta^2\nablah v)  = \eta (\deltah \eta) v + \dvh (\eta \abs{\nablah \eta}{2} \nablah \eta \otimes v)\\
		~~~~  - \abs{\nablah\eta}{4}v - \eta^2 \abs{v}{3}v
	   & \text{in} ~ \Omega,\\
	   w = - \int_0^z \bigl( \dvh \widetilde v + 2 \widetilde v \cdot\nablah \log \eta \bigr) \,dz' & \text{in} ~ \Omega,  \\
	   \partial_z \eta = 0 & \text{in} ~ \Omega,
	\end{cases}
\end{equation}
with the boundary condition
\begin{equation}\label{GA-BC}
	\dz v\bigr|_{z=0,1} = 0,
\end{equation}
where the momentum equation \subeqref{eq:GA}{2} is obtained by making use of \subeqref{eq:approximating-CPE}{1} and \subeqref{eq:approximating-CPE}{2}. Notice, we are using the notation $ \eta = \rho^{1/2} $.

It is worth mentioning that a modified version of the
incompressible Navier-Stokes equations with similar nonlinear viscosity as in \subeqref{eq:GA}{1} have been studied in \cite[Page 194]{LadyzhanskayaBook}, and also \cite{Ladyzhanskaya, LadyzhanskayaSeregin}. Moreover, Smagorinsky \cite{smagorinsky1963general} has proposed this type of eddy viscosity in the context of turbulence modeling in meteorology.

We now describe our modified Galerkin approximation scheme.
Let $  e_i = e_i(x,y,z), i = 1,2\cdots $ be an orthonormal basis of $ L^2(\Omega) $ of eigenfunctions of the Laplacian operator $ - \Delta = - \deltah - \partial_{zz} $, with the Neumann boundary condition in the $ z $-variable and periodic in the horizontal variable, corresponding to the eigenvalue $ \lambda_i $, i.e. $ - \Delta e_i = \lambda_i e_{i}, \dz e_i|_{z=0,1} = 0 $ and $ \lambda_1 < \lambda_2 \leq \cdots $. Here $ \lambda_1 = 0, e_1 \equiv 1 $ and $ \int e_i \idx = 0 $ for $ i \geq 2 $. 
For any $ n \in \mathbb N^+ $,
define the space $ X_n : = [\spn \lbrace e_i, 1 \leq i \leq n \rbrace]^2 $.
For any fixed $ 0 < T < \infty $, given $  v_n \in C(0,T;X_n) $, we will construct the following map. As before, denote $ \overline v_n := \int_0^1 v_n\,dz' $ and $ \widetilde v_n := v_n - \overline v_n $. We set
\begin{equation*}
	\mathfrak S: v_n \rightsquigarrow (\eta_n,w_n),
\end{equation*}
where $ \eta_n $ is the unique solution to \subeqref{eq:GA}{1} and $ w_n $ is given by \subeqref{eq:GA}{3} with $(\eta, v) = (\eta_n, v_n) $. That is, $ \mathfrak S $ is the solution map of  equations \subeqref{eq:GA}{1} and \subeqref{eq:GA}{3} with given $ v_n \in C([0,T]; X_n) $.
Then we will apply the Banach fixed-point theorem to show that there exists a unique approximating solution in the space $ C(0,T;X_n) $ to
the projection into the finite-dimensional space $ X_n $,
at the $ n $-level of equation
 \subeqref{eq:GA}{2}. Then we will establish some estimates independent of $ n $. By employing a compactness theorem, we will extract a subsequence which converges to a solution to \eqref{eq:GA}, as $ n \rightarrow \infty $.

\subsubsection{On the equations \subeqref{eq:GA}{1} and \subeqref{eq:GA}{3}: the map $ \mathfrak S $}\label{sec:approx-sol-exist-themap}


We start by studying the map $ \mathfrak S $. It suffices to find the solution to \subeqref{eq:GA}{1}, since once $ \eta $ is obtained for a given $ v $, \subeqref{eq:GA}{3} follows by direct substitution.

Consider $ v \in C([0,T];C^\infty(\overline\Omega)) $ satisfying the boundary condition \eqref{CPE:bc}. Then $ \overline v \in C([0,T];C^\infty(\Omega_h)) $. We will find the strong solution to \subeqref{eq:GA}{1}. Notice, that equation \subeqref{eq:GA}{1} involves a singular term on the right-hand side. We shall consider the following regularization of \subeqref{eq:GA}{1}. Let $ \eta_\delta $ be the solution to the following regularized problem. For $ \delta > 0 $, consider the equation
\begin{equation}\label{eq:GA-002}
\begin{aligned}
	& 2 \dt \eta_\delta - \dvh((1+\abs{\nablah \eta_\delta}{2})\nablah\eta_\delta) = - \dvh (\eta_\delta \overline v) \\
	& ~~~~ ~~~~ - \overline v \cdot \nablah \eta_\delta + (\eta_\delta^{2} + \delta )^{-p_0-1/2},
\end{aligned}
\end{equation}
with initial data $ \eta_\delta|_{t=0} = \rho_{\varepsilon,0}^{1/2} $ in $ \Omega_h $.
Notice that, once $ \eta_\delta $ is obtained, formally, by taking $ \delta \rightarrow 0^+ $, $ \eta : = \lim_{\delta \rightarrow 0^+} \eta_\delta $ will be the solution to \subeqref{eq:GA}{1} by means of $ \delta $-independent estimates.
In order to solve equation \eqref{eq:GA-002}, we will perform another layer of the Galerkin's approximation here. In fact, we will prove the following:
\begin{prop}\label{lm:exist-regularization-density-eq}
	Consider $ v \in C(0,T;C^\infty(\overline{\Omega})) $ satisfying the boundary condition \eqref{CPE:bc}. There is a unique strong solution $ \eta_\delta $ to \eqref{eq:GA-002} with given initial data $ \rho_{\varepsilon,0}^{1/2} $. In particular, $ \eta_\delta $ admits the following estimates which are independent of $ \delta $.
	\begin{equation}\label{est:delta-independent-density}
		\begin{aligned}
			& C_{T,\overline v}^{-1} < \eta_\delta < C_{T,\overline v},
			 \sup_{0\leq t\leq T} \bigl\lbrace \normh{\eta_{\delta}(t)}{L^2} + \normh{\nablah\eta_{\delta}(t)}{L^2} + \normh{\nablah \eta_{\delta}(t)}{L^4} \bigr\rbrace\\
			 & ~~~~  + \normh{\nablah \eta_{\delta}}{L^2(\Omega_h\times(0,T))}
			  + \normh{\nablah \eta_{\delta}}{L^4(\Omega_h\times(0,T))} + \normh{\nablah^2 \eta_{\delta}}{L^2(\Omega_h\times(0,T))} \\
			  & ~~~~ + \normh{\abs{\nablah \eta_{\delta}}{} \abs{\nablah^2 \eta_{\delta}}{}}{L^2(\Omega_h\times(0,T))}
			  + \normh{\abs{\nablah \eta_{\delta}}{2} \abs{\nablah^2 \eta_{\delta}}{}}{L^2(\Omega_h\times(0,T))} \\
			  & ~~~~ + \normh{\dt \eta_\delta}{L^2(\Omega_h\times(0,T))}
				< C_{T,\overline v},
		\end{aligned}
	\end{equation}
	for some positive, finite constant $ C_{T,\overline v} $ which depends on $$ T, \normh{\overline v}{L^\infty(\Omega_h\times(0,T))}, \normh{\nablah \overline v}{L^\infty(\Omega_h\times(0,T))}, $$ and is independent of $ \delta $.
\end{prop}

\begin{pf}
For any $ m \in \mathbb N^+ $, let $ \bar X_m : = \spn\lbrace \bar e_i, 1 \leq i\leq m \rbrace $ be the eigen-space consisting of the first $ m $ orthonormal eigenfunctions of the two-dimensional Laplacian with the periodic boundary condition in $ L^2(\Omega_h) $. That is $ - \deltah \bar e_i = \bar \lambda_i \bar e_i $, where $ \bar\lambda_i $ are the eigenvalues with $ \bar\lambda_1 \leq \bar\lambda_2 \leq \cdots $. Here $ \normh{\bar e_i}{L^2} = 1 $, $ \bar \lambda_1 = 0, \bar e_1 \equiv 1 $ and $ \int_{\Omega_h} \bar e_i \idxh = 0 $ for $ i \geq 2 $.
Let $ \eta_{\delta, m}(t) := \sum_{i=1}^{m} a_i(t) \bar e_i \in \bar X_m $ and $ \eta_{\delta,m}|_{t=0} = \sum_{i=1}^m a_{i,0} \bar e_i = \bar P_{m} \rho_{\varepsilon, 0}^{1/2} \in \bar X_m $, being the projection of $ \rho_{\varepsilon,0}^{1/2} $ on the space $ \bar X_m $. Here $ \bar P_m $ is the $ L^2 $ projection operator onto $ \bar X_m $.

{\noindent \bf Step 1: Solving the regularization problem \eqref{eq:GA-002}\par}
\noindent Consider the following ODE system,
\begin{equation}\label{eq:GA-004}
\begin{aligned}
	& \dt \eta_{\delta,m} =  \dfrac{1}{2} \bar P_{m} \bigl( \dvh((1+\abs{\nablah \eta_{\delta,m}}{2})\nablah\eta_{\delta,m}) \\
	& ~~~~ ~~~~ - \dvh (\eta_{\delta,m} \overline v) - \overline v \cdot \nablah \eta_{\delta,m} + (\eta_{\delta,m}^{2} + \delta )^{-p_0-1/2} \bigr).
	\end{aligned}
\end{equation}
Notice, that the above equation is an implicit form of the equations of $ \lbrace a_i \rbrace_{i=1,2\cdots m} $ and the right-hand side is Lipschitz in $ \lbrace a_i \rbrace_{i=1,2\cdots m} $. Then by using the standard ODE theory, there exists a unique local solution in the time interval $ (0, T_{\delta,m}) $. In the meantime, multiply \eqref{eq:GA-004} with $ 2 \eta_{\delta,m} $ and integrate the resultant. One has, after applying integration by parts,
\begin{equation}\label{GA-006}
	\begin{aligned}
		& \dfrac{d}{dt} \int_{\Omega_h} \eta_{\delta,m}^2 \idxh + \int_{\Omega_h} \biggl( \abs{\nablah \eta_{\delta,m}}{2} + \abs{\nablah \eta_{\delta,m}}{4} \biggr) \idxh \\
		& ~~~~ = \int_{\Omega_h} \eta_{\delta,m} (\eta_{\delta,m}^2 + \delta)^{-p_0-1/2} \idxh \leq C_\delta + \int_{\Omega_h} \eta_{\delta,m}^2\idxh.
	\end{aligned}
\end{equation}
Therefore, the Gr\"onwall's inequality implies $ T_{\delta,m} = T $ and that
\begin{equation*}
	\sup_{0\leq t\leq T} \normh{\eta_{\delta,m}(t)}{L^2} + \normh{\nablah \eta_{\delta,m}}{L^2(\Omega_h\times(0,T))} + \normh{\nablah \eta_{\delta,m}}{L^4(\Omega_h\times(0,T))} < C_{T,\delta}.
\end{equation*}
On the other hand, multiply \eqref{eq:GA-004} with $ - 2 \bar P_m (\dvh ((1+\abs{\nablah \eta_{\delta,m}}{2})\nablah \eta_{\delta,m})) $ and integrate the resultant. After applying integration by parts and the Cauchy-Schwarz inequality,
\begin{equation}\label{GA-007}
	\begin{aligned}
		& \dfrac{d}{dt} \int_{\Omega_h} \biggl( \abs{\nablah \eta_{\delta,m}}{2} + \dfrac{1}{2}\abs{\nablah \eta_{\delta,m}}{4} \biggr) \idxh \\
		& ~~~~ + \int_{\Omega_h} \abs{\bar P_m (\dvh ((1+\abs{\nablah \eta_{\delta,m}}{2})\nablah \eta_{\delta,m}))}{2} \idxh \leq C_\delta \\
		& ~~~~ + C (\normh{\overline v}{L^\infty\spacetime}^2 + \normh{\nablah \overline v}{L^\infty\spacetime}^2)\int_{\Omega_h} \biggl( \eta_{\delta,m}^{2} \\
		& ~~~~ ~~~~ + \abs{\nablah \eta_{\delta,m}}{2} \biggr) \idxh.
	\end{aligned}
\end{equation}
Notice that, $ \deltah \eta_{\delta,m} \in \bar X_{m} $,
\begin{align*}
	& \int_{\Omega_h} \abs{\bar P_m (\dvh ((1+\abs{\nablah \eta_{\delta,m}}{2})\nablah \eta_{\delta,m}))}{2} \idxh = \int_{\Omega_h} \biggl( \abs{\deltah \eta_{\delta,m}}{2} \\
	& ~~~~ + 2 \deltah \eta_{\delta,m} \dvh (\abs{\nablah\eta_{\delta,m}}{2}\nablah \eta_{\delta,m}) + \abs{\bar P_{m} \dvh(\abs{\nablah \eta_{\delta,m}}{2}\nablah \eta_{\delta,m})}{2} \biggr) \idxh \\
	& = \int_{\Omega_h} \biggl( \abs{\deltah \eta_{\delta,m}}{2} + 2 \abs{\nablah \eta_{\delta,m}}{2} \abs{\nablah^2 \eta_{\delta,m}}{2} + \abs{\nablah \abs{\nablah \eta_{\delta,m}}{2}}{2} \\
	& ~~~~ ~~~~ + \abs{\bar P_{m} \dvh(\abs{\nablah \eta_{\delta,m}}{2}\nablah \eta_{\delta,m})}{2} \biggr)  \idxh.
\end{align*}
Then, by applying the Gr\"onwall's inequality, we have from \eqref{GA-007} the following,
\begin{equation}\label{17July2018-01}
	\begin{aligned}
		& \sup_{0\leq t\leq T} \lbrace \normh{\nablah\eta_{\delta,m}(t)}{L^2}^2 + \normh{\nablah \eta_{\delta,m}(t)}{L^4}^4 \rbrace + \normh{\nablah^2 \eta_{\delta,m}}{L^2(\Omega_h\times(0,T))}^2  \\
		& ~~ + \normh{\abs{\nablah \eta_{\delta,m}}{} \abs{\nablah^2 \eta_{\delta,m}}{}}{L^2(\Omega_h\times(0,T))}^2 + \normh{\bar P_{m} \dvh(\abs{\nablah \eta_{\delta,m}}{2}\nablah \eta_{\delta,m})}{L^2(\Omega_h\times(0,T))}^2 \\
		& ~~~~ < C_{T,\delta}(1 + \normh{\overline v}{L^\infty(\Omega_h\times(0,T))}^2 + \normh{\nablah \overline v}{L^\infty(\Omega_h\times(0,T))}^2).
	\end{aligned}
\end{equation}
Moreover from \eqref{eq:GA-004}
\begin{align*}
	& \normh{\dt \eta_{\delta,m}}{L^2(\Omega_h\times(0,T))} \leq C \normh{\nablah^2 \eta_{\delta,m}}{L^2(\Omega_h\times(0,T))} \\
	& ~~~~ + C\normh{\bar P_{m} \dvh(\abs{\nablah \eta_{\delta,m}}{2}\nablah \eta_{\delta,m})}{L^2(\Omega_h\times(0,T))} \\
	& ~~~~ + C \normh{\overline v}{L^\infty(\Omega_h\times(0,T))} \normh{\nablah \eta_{\delta,m}}{L^2(\Omega_h\times(0,T))} \\
	& ~~~~ + C \normh{\nablah \overline v}{L^\infty(\Omega_h\times(0,T))} \normh{\eta_{\delta,m}}{L^2(\Omega_h\times(0,T))} + C_\delta\\
	& ~~~~ \leq C_{T,\delta}(1 + \normh{\overline v}{L^\infty(\Omega_h\times(0,T))} + \normh{\nablah \overline v}{L^\infty(\Omega_h\times(0,T))}).
\end{align*}
Therefore,
applying the Aubin's compactness Theorem (see, e.g., \cite[Theorem 2.1]{Temam1984} and \cite{Simon1986,Chen2012}), the estimates above imply that there exists a subsequence, denoted also by $ \lbrace \eta_{\delta,m} \rbrace $, and that there is a $ \eta_\delta \in L^\infty(0,T;W^{1,4}(\Omega_h)) \cap L^2(0,T;H^2(\Omega_h))  $, $ \dt \eta_\delta \in L^2(\Omega_h\times(0,T)) $, such that, for $(p,q) \in  (1,\infty) \times [2,\infty)$,
\begin{equation}\label{GA-17July2018}
\begin{aligned}
	& \eta_{\delta,m} \rightarrow \eta_{\delta}, ~~ & \text{in}& ~ L^2(0,T; W^{1,p}(\Omega_h)) \cap C(0,T; L^q(\Omega_h)  , \\
	& \eta_{\delta,m} \buildrel\ast\over\rightharpoonup \eta_\delta, ~~  & \text{weak-$\ast$ in}& ~ L^\infty(0,T; W^{1,4}(\Omega_h)),\\
	& \eta_{\delta,m} \rightharpoonup \eta_\delta, ~~ & \text{weakly in}& ~ L^2(0,T;H^2(\Omega_h)),\\
	& \dt \eta_{\delta,m} \rightharpoonup \dt \eta_\delta, ~~ & \text{weakly in} & ~ L^2(\Omega_h\times(0,T)).
\end{aligned}
\end{equation}
Consequently, let $ \psi \in C_c^\infty(\overline\Omega_h \times [ 0, T) ) $ be the test function. After taking the $ L^2 $-inner product of \eqref{eq:GA-004} with $ \psi $ in the space-time domain $ \Omega_h \times (0,T) $ and applying integration by parts in the resultant, we have the identity
\begin{align*}
	& \int_0^T \int_{\Omega_h} \dt \eta_{\delta,m}\psi \idxh \,dt = - \dfrac{1}{2} \int_0^T \int_{\Omega_h} (1+ \abs{\nablah \eta_{\delta,m}}{2}) \nablah \eta_{\delta,m} \cdot \nablah \bar P_m \psi \idxh \,dt\\
	& ~~~~ + \int_0^T \int_{\Omega_h} \eta_{\delta,m} \overline v \cdot \nablah \bar P_m \psi \idxh \,dt + \int_0^T \int_{\Omega_h} \biggl( (- \overline v \cdot \nablah \eta_{\delta,m} \\
	& ~~~~ ~~~~ + (\eta_{\delta,m}^2 + \delta)^{-p_0-1/2} ) \bar P_{m} \psi \biggr)  \idxh \,dt.
\end{align*}
Then \eqref{GA-17July2018} implies that after taking $ m \rightarrow \infty $ in the above equation, we have,
\begin{align*}
	& \int_0^T \int_{\Omega_h} \dt \eta_{\delta}\psi \idxh \,dt = - \dfrac{1}{2} \int_0^T \int_{\Omega_h} (1+ \abs{\nablah \eta_{\delta}}{2}) \nablah \eta_{\delta} \cdot \nablah \psi \idxh \,dt\\
	& ~~~~ + \int_0^T \int_{\Omega_h} \eta_{\delta} \overline v \cdot \nablah \psi \idxh \,dt + \int_0^T \int_{\Omega_h} \biggl( (- \overline v \cdot \nablah \eta_{\delta} \\
	& ~~~~ ~~~~ + (\eta_{\delta}^2 + \delta)^{-p_0-1/2} )  \psi \biggr)  \idxh \,dt,
\end{align*}
where we have applied the fact that $ \bar P_m \psi \rightarrow \psi $ strongly in $ L^p(\Omega_h\times(0,T)) $ for every $ p \in (1, \infty) $.
Therefore $ \eta_\delta $ is a weak solution to \eqref{eq:GA-002}. Moreover, applying \eqref{GA-17July2018} again, for any $ \psi \in C_c^\infty(\overline\Omega_h\times (0,T)) $,
\begin{align*}
	& \int_0^T \int_{\Omega_h} \abs{\nablah \eta_{\delta}}{} \abs{\nablah^2 \eta_\delta}{} \psi  \idxh\,dt = \lim\limits_{m \rightarrow \infty} \int_0^T \int_{\Omega_h} \abs{\nablah \eta_{\delta,m}}{} \abs{\nablah^2 \eta_{\delta,m}}{} \psi  \idxh\,dt \\
	& ~~~~ \leq C \biggl( \sup_{m \geq 1} \normh{\abs{\nablah \eta_{\delta,m}}{} \abs{\nablah^2 \eta_{\delta,m}}{}}{L^2(\Omega_h \times (0,T))} \biggr) \times \normh{\psi}{L^2(\Omega_h\times(0,T))} , \\
	& \int_0^T \int_{\Omega_h}  \abs{\nablah \eta_\delta}{2} \nablah \eta_\delta \cdot \nablah \psi  \idxh\,dt \\
	&  = \lim\limits_{m\rightarrow\infty} \int_0^T \int_{\Omega_h} \abs{\nablah \eta_{\delta,m}}{2} \nablah \eta_{\delta,m} \cdot \nablah \bar P_m \psi  \idxh\,dt \\
	& = - \lim\limits_{m\rightarrow\infty} \int_0^T \int_{\Omega_h} \bar P_m \dvh (\abs{\nablah \eta_{\delta,m}}{2} \nablah \eta_{\delta,m})  \psi  \idxh\,dt\\
	& ~~~~ \leq C \sup_{m\geq1} \normh{\bar P_m \dvh (\abs{\nablah \eta_{\delta,m}}{2} \nablah \eta_{\delta,m})}{L^2(\Omega_h \times (0,T))} \times  \normh{\psi}{L^2(\Omega_h\times(0,T))}.
\end{align*}
Together with \eqref{17July2018-01}, these inequalities imply that
\begin{align*}
	& \normh{\abs{\nablah \eta_\delta}{}\abs{\nablah^2 \eta_\delta}{}}{L^2(\Omega_h\times(0,T))} + \normh{\dvh(\abs{\nablah \eta_{\delta}}{2}\nablah \eta_{\delta})}{L^2(\Omega_h\times(0,T))}^2 \\
	& ~~~~ \leq C_{T,\delta}(\normh{\overline v}{L^\infty(\Omega_h\times(0,T))}, \normh{\nablah \overline v}{L^\infty(\Omega_h\times(0,T))}).
\end{align*}
Consequently, similar arguments imply that
$ \eta_\delta $ is the strong solution  to \eqref{eq:GA-002} with $ \eta_\delta|_{t=0} = \rho_{\varepsilon,0}^{1/2} $ satisfying
\begin{equation}\label{GA-005}
	\begin{aligned}
		& \sup_{0\leq t\leq T} \lbrace \normh{\eta_{\delta}(t)}{L^2} + \normh{\nablah\eta_{\delta}(t)}{L^2} + \normh{\nablah \eta_{\delta}(t)}{L^4} \rbrace  + \normh{\nablah \eta_{\delta}}{L^2(\Omega_h\times(0,T))} \\
		& ~~~~ + \normh{\nablah \eta_{\delta}}{L^4(\Omega_h\times(0,T))} + \normh{\nablah^2 \eta_{\delta}}{L^2(\Omega_h\times(0,T))} + \normh{\abs{\nablah \eta_{\delta}}{} \abs{\nablah^2 \eta_{\delta}}{}}{L^2(\Omega_h\times(0,T))}\\
		& ~~~~ + \normh{\abs{\nablah \eta_{\delta}}{2} \abs{\nablah^2 \eta_{\delta}}{}}{L^2(\Omega_h\times(0,T))} + \normh{\dt \eta_\delta}{L^2(\Omega_h\times(0,T))} \\
		&  < C_{T,\delta}(\normh{\overline v}{L^\infty(\Omega_h\times(0,T))}, \normh{\nablah \overline v}{L^\infty(\Omega_h\times(0,T))}).
	\end{aligned}
\end{equation}
where we have applied the fact
\begin{align*}
	& \normh{\dvh(\abs{\nablah \eta_{\delta}}{2}\nablah \eta_{\delta})}{L^2(\Omega_h\times(0,T))}^2 = \int_0^T \int_{\Omega_h} \abs{\dvh(\abs{\nablah \eta_{\delta}}{2}\nablah \eta_{\delta})}{2}\idxh\,dt \\
	& ~~~~ = \int_0^T \int_{\Omega_h} \biggl( \abs{\nablah \eta_\delta}{4} \abs{\nablah^2 \eta_\delta}{2} + \abs{\nablah \eta_\delta \cdot \nablah\abs{\nablah \eta_\delta}{2}}{2} \\
	& ~~~~ ~~~~ + \abs{\nablah \eta_\delta}{2} \abs{\nablah \abs{\nablah \eta_{\delta}}{2}}{2} \biggr) \idxh\,dt
	\geq \normh{\abs{\nablah \eta_{\delta}}{2} \abs{\nablah^2 \eta_{\delta}}{}}{L^2(\Omega_h\times(0,T))}^2.
\end{align*}

{\noindent\bf Step 2: $ \delta $-independent estimates\par}
\noindent It suffices to show that $ \eta_\delta $, see \eqref{eq:GA-002}, admits uniform upper and lower bounds independent of $ \delta $. In fact, the above estimates will be independent of $ \delta $ once we have $ \eta_{\delta} > C $ a priorly, for some positive constant $ C $ independent of $ \delta $.

Multiply \eqref{eq:GA-002} with $ \eta_\delta^{2p_0 + 1} $ and with $ - \dvh((1+\abs{\nablah \eta_\delta}{2})\nablah\eta_\delta) $, respectively, and integrate the resultants in the space-time domain $ \Omega_h \times (0,T) $. Then after applying integration by parts and the Young's inequality, one obtains
\begin{align*}
	& \sup_{0\leq t\leq T} \bigl\lbrace \dfrac{1}{p_0+1}\int_{\Omega_h} \abs{\eta_\delta}{2p_0+2} \idxh \bigr\rbrace + (2p_0+1) \int_0^{T} \int_{\Omega_h} \biggl( \eta_\delta^{2p_0}  \abs{\nablah \eta_\delta}{2} \\
	& ~~~~ + \eta_\delta^{2p_0} \abs{\nablah \eta_\delta}{4} \biggr) \idxh \,dt
	 \leq  C \normh{\overline v}{L^\infty(\Omega_h \times(0,T))} \int_0^{T} \int_{\Omega_h} \eta_\delta^{2p_0+2} \idxh \,dt \\
	& ~~~~ + C\int_0^{T} \int_{\Omega_h} \underbrace{\eta_\delta^{2p_0+1} (\eta_\delta^2 + \delta)^{-p_0-1/2}}_{\leq C} \idxh\,dt + C, \\
	& \sup_{0\leq t\leq T} \bigl\lbrace \int_{\Omega_h} \abs{\nablah \eta_\delta}{2} + \dfrac{1}{2} \abs{\nablah \eta_\delta}{4} \idxh  \bigr\rbrace
	 + \int_0^{T} \int_{\Omega_h} \biggl( \dfrac{1}{2} \abs{\deltah\eta_\delta}{2} \\
	 & ~~~~ + \abs{\nablah\eta_\delta}{2}\abs{\nablah^2\eta_\delta}{2}
	+ \dfrac{1}{2} \abs{\nablah \eta_\delta }{4} \abs{\nablah^2 \eta_\delta}{2} + (2p_0+1)(\eta_{\delta}^2+\delta)^{-p_0-3/2}\\
	& ~~~~ ~~~~ \times (1+\abs{\nablah \eta_\delta}{2}) \abs{\nablah \eta_\delta}{2} \biggr) \idxh \,dt
	 \leq C \bigl(\normh{\overline v}{L^\infty(\Omega_h\times(0,T))}^2 \\
	 & ~~~~  + \normh{\nablah \overline v}{L^\infty(\Omega_h\times(0,T))}^2\bigr) \int_0^{T} \int_{\Omega_h} \biggl( 1 + \eta_{\delta}^{2p_0+2}
	  + \abs{\nablah \eta_{\delta}}{2} \biggr) \idxh\,dt + C,
\end{align*}
where $ C $ is independent of $ \delta $.
Therefore, after applying the Gr\"onwall's inequality to the above inequalities, and the following Gagliardo-Nirenberg inequality, we have
\begin{equation}\label{18July2018}
\begin{aligned}
	& \normh{\eta_\delta }{L^\infty(\Omega_h\times(0,T))} \leq C \sup_{0\leq t\leq T} \lbrace \normh{\eta_\delta}{L^{2p_0+2}} + \normh{\nablah \eta_\delta}{L^4} \rbrace \\
	& ~~~~ < C_{T}(\normh{\overline v}{L^\infty(\Omega_h\times(0,T))}, \normh{\nablah \overline v}{L^\infty(\Omega_h\times(0,T))}),
\end{aligned}
\end{equation}
where the right-hand side of the last inequality is
independent of $ \delta $. On the other hand, let $ \sigma_\delta := \eta_\delta^{-1} $. Then $ \sigma_\delta $ satisfies
\begin{equation*}
	\begin{aligned}
		& 2 \dt \sigma_\delta + 2 \overline v \cdot\nablah \sigma_\delta - \sigma_\delta \dvh \overline v + \sigma_\delta^2 (\sigma^{-2}_\delta +\delta)^{-p_0-1/2} + 2  \sigma_\delta^{-1} \abs{\nablah \sigma_\delta}{2} \\
		& ~~~~ + 2  \sigma^{-5}_\delta \abs{\nablah \sigma_\delta}{4}
		= \deltah \sigma_\delta + \dvh (\sigma^{-4}_\delta \abs{\nablah \sigma_\delta}{2} \nablah \sigma_\delta).
	\end{aligned}
\end{equation*}
Notice we have $ \sigma_\delta^{-1} = \eta_\delta < C_{T,\overline v} < \infty $ from \eqref{18July2018} and hence there is no singularity in the above equation. Here and in the following, $ C_{T,\overline v} $ denotes a constant depending on $ T, \normh{\overline v}{L^\infty(\Omega_h\times(0,T))}, \normh{\nablah \overline v}{L^\infty(\Omega_h\times(0,T))} $.
Multiply this equation with $ (3\sigma_{\delta}^5)_+ := 3 \sigma_\delta^5 \mathbbm 1_{\lbrace\sigma_\delta > 0\rbrace} $ and $ (\sigma_\delta - k)_+ := (\sigma_\delta - k ) \mathbbm 1_{\lbrace \sigma_\delta > k \rbrace} $, respectively, with some $ k \geq \normh{\rho_{\varepsilon,0}^{-1/2}}{L^\infty} = \normh{\sigma_\delta|_{t=0}}{L^\infty} > 0 $, and integrate the resultant. After applying integration by parts, it holds
\begin{align*}
	& \sup_{0\leq t\leq T} \int_{\Omega_h} \abs{(\sigma_{\delta})_+}{6} \idx + 9 \int_0^T \int_{\Omega_h} \biggl( (\sigma_{\delta}^4)_+ \abs{\nablah (\sigma_{\delta})_+}{2} + \abs{\nablah (\sigma_{\delta})_+}{4}\biggr) \idxh\,dt \\
	& ~~~~ \leq C \normh{\nablah \overline v}{L^\infty(\Omega_h\times(0,T))} \int_0^T \int_{\Omega_h} (\sigma_{\delta}^6)_+ \idxh \,dt + C , \\
	& \sup_{0<t<T} \int_{\Omega_h} \abs{(\sigma_\delta - k)_+}{2} \idxh + \dfrac{1}{2} \int_0^T \int_{\Omega_h}\biggl( \abs{\nablah (\sigma_\delta - k)_+}{2} \\
	& ~~~~ ~~~~ + \sigma_\delta^{-4} \abs{\nablah (\sigma_\delta - k)_+}{4} \biggr) \idxh \,dt\\
	& ~~~~ \leq C (\normh{\nablah \overline v}{L^\infty(\Omega_h\times(0,T))}^2 + 1) \int_0^T \int_{\Omega_h} \abs{(\sigma_\delta-k)_+}{2} \idxh\,dt \\
	& ~~~~ ~~~~ + C \int_0^T \int_{\Omega_h}\mathbbm{1}_{\lbrace \sigma_\delta > k \rbrace } \abs{\sigma_\delta}{2} \idxh \,dt .
\end{align*}
Then applying the Gr\"onwall's inequality implies
\begin{align*}
	& \sup_{0\leq t\leq T} \int_{\Omega_h} \abs{(\sigma_{\delta})_+}{6} \idx + \int_0^T \int_{\Omega_h} \biggl( (\sigma_{\delta}^4)_+ \abs{\nablah (\sigma_{\delta})_+}{2} + \abs{\nablah (\sigma_{\delta})_+}{4}\biggr) \idxh\,dt \\
	& ~~~~ ~~~~ \leq C_{T, \overline v} < \infty, \\
	& \sup_{0\leq t\leq T} \int \abs{(\sigma_\delta - k)_+}{2} \idxh +  \int_0^T \int_{\Omega_h} \biggl( \abs{\nablah (\sigma_\delta - k)_+}{2} \\
	& ~~~~ ~~~~ + \sigma_\delta^{-4} \abs{\nablah (\sigma_\delta - k)_+}{4}\biggr) \idxh \,dt
	 \leq C_{T,\overline v} \int_0^T \int_{\Omega_h}\mathbbm{1}_{\lbrace \sigma_\delta > k > 0 \rbrace } \abs{\sigma_\delta}{2} \idxh \,dt \\
	 & ~~~~ \leq C_{T,\overline v} \abs{\lbrace \sigma_\delta > k \rbrace }{2/3} (\int_0^T \int_{\Omega_h} \abs{(\sigma_{\delta})}{6}\idxh\,dt)^{1/3}
	\leq C_{T, \overline v}\abs{\lbrace \sigma_\delta > k \rbrace }{2/3}.
\end{align*}
Then using similar arguments as in \eqref{UL-010}, we have, after applying interpolation inequality, for $ l > k $,
\begin{align*}
	& \abs{\lbrace \sigma_\delta > l \rbrace}{}\leq (l-k)^{-4}\normh{(\sigma_\delta-k)_+}{L^4(\Omega_h\times(0,T))}^4  \leq  C (l-k)^{-4} \\
	& ~~~~ \times \bigl( \sup_{0\leq t\leq T} \int_{\Omega_h}\abs{(\sigma_\delta - k)_+}{2} \idxh
	 \int_0^T \int_{\Omega_h} \abs{\nablah (\sigma_\delta-k)_+}{2} \idxh \,dt \\
	 & ~~~~ + ( \sup_{0\leq t\leq T} \int_{\Omega_h}\abs{(\sigma_\delta - k)_+}{2} \idxh)^2\bigr)
	\leq C_{T,\overline v} (l-k)^{-4} \abs{\lbrace \sigma_\delta > k \rbrace }{4/3}.
\end{align*}
Then using the De Giorgi-type Lemma (i.e., Lemma \ref{lm:De-Giorgi}) as before, from \eqref{UL-010} to \eqref{UL-009} in section \ref{sec:ul-bounds}, yield that there is a positive constant $ C'_{T,\overline v} < \infty $ such that \begin{equation*}
	{\sigma_\delta} < C'_{T, \overline v}.
\end{equation*}
Therefore since $ \eta_\delta|_{t=0} > 0 $, it follows from the upper bounds of $ \abs{\eta_\delta}{}, {\sigma_\delta}{} = {\eta_\delta^{-1}}{} $ and the continuity of $ \eta_\delta $ in space-time that there is a positive constant $ 0 < C_{T,\overline v} < \infty $ depending on
$$ T, \normh{\overline v}{L^\infty(\Omega_h\times(0,T))}, \normh{\nablah \overline v}{L^\infty(\Omega_h\times(0,T))}, $$
such that
\begin{equation}
	1/C_{T,\overline v} < \eta_\delta < C_{T,\overline v}.
\end{equation}
Therefore, by performing similar arguments as in \eqref{GA-006}, \eqref{GA-007}, the estimate \eqref{GA-005} is independent of $ \delta $.
This finishes the proof.
\end{pf}
With Proposition \ref{lm:exist-regularization-density-eq} at hand, after letting $ \delta \rightarrow 0^+ $ in \eqref{eq:GA-002},
one will get the strong solution $ \eta $ to \subeqref{eq:GA}{1}.
\begin{prop}\label{lm:the-map-S-approx}
	Under the same assumption as in Proposition \ref{lm:exist-regularization-density-eq}, there is a unique strong solution $ \eta $ to \subeqref{eq:GA}{1} with the following estimates.
	\begin{equation}\label{est:the-map-S}
	\begin{aligned}
		&C^{-1}_{T,\overline v} < \eta < C_{T,\overline v}, ~ \sup_{0\leq t\leq T} \bigl\lbrace \normh{\eta_{}(t)}{L^2} + \normh{\nablah\eta_{}(t)}{L^2} + \normh{\nablah \eta_{}(t)}{L^4} \bigr\rbrace \\
		& ~~~~ + \normh{\nablah \eta_{}}{L^2(\Omega_h\times(0,T))}
		+ \normh{\nablah \eta_{}}{L^4(\Omega_h\times(0,T))}
		+ \normh{\nablah^2 \eta_{}}{L^2(\Omega_h\times(0,T))} \\
		& ~~~~ + \normh{\abs{\nablah \eta_{}}{} \abs{\nablah^2 \eta_{}}{}}{L^2(\Omega_h\times(0,T))}
		 + \normh{\abs{\nablah \eta_{}}{2} \abs{\nablah^2 \eta_{}}{}}{L^2(\Omega_h\times(0,T))} \\
		 & ~~~~
		+ \normh{\dt \eta}{L^2(\Omega_h\times(0,T))} < C_{T,\overline v},
	\end{aligned}
\end{equation}
for some positive constant $ C_{T,\overline v} =  C_{T,\overline v}(T, \normh{\overline v}{L^\infty(\Omega_h\times(0,T))}, \normh{\nablah \overline v}{L^\infty(\Omega_h\times(0,T))}) < \infty$.
Moreover,
\begin{equation}\label{18July2018-03}
\begin{aligned}
	& \normh{\eta_1 - \eta_2}{L^\infty(0,T;L^2(\Omega_h))} \leq C_{T,\overline v_1, \overline v_2} ( \normh{\eta_{1,0} - \eta_{2,0}}{L^2} \\
	& ~~~~ + T \cdot \normh{\overline v_1-\overline v_2}{L^\infty(\Omega_h\times (0,T))}
	 + T \cdot \normh{\nablah (\overline v_1-\overline v_2)}{L^\infty(\Omega_h\times (0,T))}),
	\end{aligned}
\end{equation}
where $ \eta_i $ satisfies \subeqref{eq:GA}{1}, with initial data $ \eta_i\bigr|_{t=0} = \eta_{i,0} $ and  $ v $ replaced by $ v_i \in C(0,T;C^\infty(\overline\Omega)) $, $ i = 1,2$, and $ C_{T,\overline v_1, \overline v_2} \in (0,\infty) $ depends on
\begin{gather*}
	T, \normh{\overline v_1}{L^\infty(\Omega_h\times (0,T))}, \normh{\nablah \overline v_1}{L^\infty(\Omega_h\times (0,T))}, \\
	\normh{\overline v_2}{L^\infty(\Omega_h\times (0,T))}, \normh{\nablah \overline v_2}{L^\infty(\Omega_h\times (0,T))}.
\end{gather*}
In particular, the map $ \mathfrak S: v \leadsto (\eta, w) $ 
for $ v \in C(0,T; C^\infty (\overline\Omega)) $ is well-defined.
\end{prop}

\begin{pf}
The existence of the strong solution $ \eta $ to \subeqref{eq:GA}{1} and the estimates \eqref{est:the-map-S} follow from similar arguments as from \eqref{GA-17July2018} to \eqref{GA-005} by replacing the sequence $ \lbrace \eta_{m,\delta} \rbrace_{m = 1,2,\cdots} $ with $ \lbrace \eta_\delta \rbrace_{\delta\in (0,\infty)} $ and taking the limit $ \delta \rightarrow 0^+ $. The uniqueness follows from \eqref{18July2018-03} by taking $ v_1 = v_2 $ and $ \eta_{1,0} = \eta_{2,0} $.

To show \eqref{18July2018-03}, denote by $ \eta_{12} := \eta_1 - \eta_2 $ where $ \eta_i $, $i=1,2$, are solutions to \subeqref{eq:GA}{1} with $ v = v_i $.
Then $ \eta_{12} $ satisfies
\begin{align*}
	& 2 \dt \eta_{12} - \deltah \eta_{12} - \bigl(\dvh ( \abs{\nablah \eta_1}{2} \nablah \eta_{1}) - \dvh ( \abs{\nablah \eta_2}{2} \nablah \eta_{2})\bigr) \\
	& ~~~~ = - \dvh (\eta_{12} \overline v_1) - \overline{v}_1 \cdot \nablah \eta_{12} - \dvh (\eta_{2}\overline{v}_{12} ) - \overline v_{12}\cdot \nablah \eta_2 \\
	& ~~~~ + ( \eta_1^{-2p_0-1} - \eta_2^{-2p_0 - 1}),
\end{align*}
where $ \overline v_{12} := \overline v_1 - \overline v_2 $.
Multiply the above equation with $ \eta_{12} $ and integrate the resultant over $ \Omega_h $. After applying integration by parts and the Cauchy-Schwarz inequality, one obtains
\begin{equation}\label{28July2018}
\begin{aligned}
	& \dfrac{d}{dt} \normh{\eta_{12}}{L^2}^2 + \normh{\nablah \eta_{12}}{L^2}^2 - \int_{\Omega_h} \biggl( \bigl(\dvh ( \abs{\nablah \eta_1}{2} \nablah \eta_{1})\\
	& ~~~~ - \dvh ( \abs{\nablah \eta_2}{2} \nablah \eta_{2})\bigr) \eta_{12} \biggr) \idxh  \\
	& ~~ = \int \bigl( \eta_2 \overline v_{12} \cdot \nablah \eta_{12} + \eta_2 \dvh( \eta_{12} \overline v_{12}) \bigr) \idxh  \\
	& ~~~~ + \int_{\Omega_h} ( \eta_1^{-2p_0-1} - \eta_2^{-2p_0 - 1}) \eta_{12}\idxh \leq \dfrac{1}{2} \normh{\nablah \eta_{12}}{L^2}^2 + C \normh{\eta_{12}}{L^2}^2 \\
	& ~~~~ + C \normh{\overline v_{12}}{L^\infty(\Omega_h\times (0,T))}^2 + C \normh{\nablah \overline v_{12}}{L^\infty(\Omega_h\times (0,T))}^2,
\end{aligned}
\end{equation}
where $ C $ depends on the upper and lower bounds of $ \eta_i $, i.e., $ C_{T,\overline v_i}^{-1} < \eta_i < C_{T,\overline v_i} $, $ i = 1,2 $. Notice that the operator $ \dvh (\abs{\nablah \eta}{2} \nablah \eta) $ in $ \eta $ admits the monotonicity (see, e.g., \cite[Page 194]{LadyzhanskayaBook}):
\begin{equation}\label{monotonicity-operator}
	- \int_{\Omega_h}  \bigl(\dvh ( \abs{\nablah \eta_1}{2} \nablah \eta_{1}) - \dvh ( \abs{\nablah \eta_2}{2} \nablah \eta_{2})\bigr) \bigl( \eta_1 -\eta_2\bigr) \idxh \geq 0.
\end{equation}
Then applying the Gr\"onwall's inequality to \eqref{28July2018} yields
\begin{align*}
	& \sup_{0\leq t\leq T} \normh{\eta_{12}(t)}{L^2}^2 \leq C \normh{\eta_{12}\bigr|_{t=0}}{L^2}^2 + C T \normh{\overline v_{12}}{L^\infty(\Omega_h\times (0,T))}^2 \\
	& ~~~~ + C T \normh{\nablah \overline v_{12}}{L^\infty(\Omega_h\times (0,T))}^2. 
\end{align*}
This finishes the proof of \eqref{18July2018-03}. Notice that \eqref{18July2018-03} implies the uniqueness and the Lipschitz continuity of $ \eta $ with respect to $ \overline v $.
\end{pf}

We present a proof of the monotonicity in \eqref{monotonicity-operator} here. In fact, after rearranging the left hand side of \eqref{monotonicity-operator} and applying integration by parts, one obtains
\begin{align*}
	& - \int_{\Omega_h}  \bigl(\dvh ( \abs{\nablah \eta_1}{2} \nablah \eta_{1}) - \dvh ( \abs{\nablah \eta_2}{2} \nablah \eta_{2})\bigr) \bigl( \eta_1 -\eta_2\bigr) \idxh\\
	& ~~~~ = \normh{\abs{\nablah \eta_1}{} \nablah \eta_{12}}{L^2}^2 + \normh{\nablah \eta_2 \cdot \nablah \eta_{12}}{L^2}^2  \\
	& ~~~~ ~~~~ + \int_{\Omega_h} (\nablah \eta_{1}\cdot \nablah \eta_{12}) (\nablah \eta_{2}  \cdot \nablah \eta_{12}) \idxh\\
	& ~~~~ \geq \dfrac{1}{2} \normh{\abs{\nablah \eta_1}{} \nablah \eta_{12}}{L^2}^2 + \dfrac{1}{2} \normh{\nablah \eta_2 \cdot \nablah \eta_{12}}{L^2}^2 \geq 0,
\end{align*}
where in the last inequality we have applied the Cauchy-Schwarz inequality.

%
%
%
%
%

\subsubsection{On the momentum equation of
the approximating system and the fixed point}\label{sec:approx-sol-exist-v}


Now we are in the place to establish the modified Galerkin approximation to \subeqref{eq:GA}{2}. This is done following a similar scheme as in \cite{Feireisl2001} with the map $ \mathfrak S $ defined in the last subsection. As mentioned before, let $ v_n = \sum_{i=1}^{n} a_{i}(t) e_i \in X_n $, for a fixed $ n \in \mathbb Z^+ $, with $ v_{n}|_{t=0} = P_n v_{\varepsilon,0} \in X_n $ is the projection of $ v_{\varepsilon,0} $ into the space $ X_n $. Here $ P_n $ is the $ L^2 $ projection operator onto the space $ X_n $. Let $ (\eta_n, w_n) := \mathfrak S (v_n) $, given by Proposition \ref{lm:the-map-S-approx}.

Consider the following ODE system for the coefficients $ a_j(t) $, $ j= 1,\cdots,n $,
\begin{equation}\label{GA-009}
	\begin{aligned}
		& \int \dt( \eta_n^2 v_{n} ) \cdot e_i \idx  = \int (-  \dvh (\eta_n^2 v_n \otimes v_n) - \dz(\eta_n^2 w_n v_n) - \nablah \eta_n^{2\gamma} \\
		& ~~~~ + \dvh(\eta_n^2\mathcal D (v_n))
		+ \dz (\eta_n^2 \dz v_n) + \dvh(\eta_n^2\nablah v_n) + \eta_n \deltah \eta_n v_n \\
		& ~~~~ + \dvh (\eta_n \abs{\nablah \eta_n}{2} \nablah \eta_n \otimes v_n)
		 - \abs{\nablah\eta_n}{4}v_n   - \eta_n^2 \abs{v_n}{3}v_n ) \cdot e_i \idx,
	\end{aligned}
\end{equation}
for $ i = 1, 2 \cdots n $, which is an implicit form of the equations for $ \lbrace a_j(t) \rbrace_{j=1,2\cdots n} $.

Define the map $ M[\eta_n]: X_n \mapsto X_n $ as follows. Let $ m_{ij}[\eta_n] := \int \eta_{n}^2 e_i \cdot e_j \idx $. For every $ \vec{u} = \sum_{i=1}^n a_i e_i \in X_n $, define $ M[\eta_n] \vec{u} := \sum_{i=1}^n (\sum_{j=1}^n a_j m_{ij}[\eta_n]) e_i \in X_n $. Thus
 $ \int M[\eta_n] \vec{a} \cdot \vec{b}\idx = \int \eta_n^2 \vec{a}\cdot\vec{b}\idx $, for every $ \vec{a},\vec{b} \in X_n $, which is positive definite, since  $ \eta_n $ is bounded from below by a strictly positive number. Thus $ M[\eta_n] $ is invertible. In particular,
 \begin{equation*}
 	(M[\eta_n])^{-1}(P_{n} (\eta_n^2 v_n)) = v_n.
 \end{equation*}

Then the following integral representation of $ v_n $ holds,
\begin{equation}\label{18July2018-02}
	\begin{aligned}
		& v_n = Q(v_n) := (M[\eta_n])^{-1} \biggl( P_n \biggl\lbrack  \int_0^t (-  \dvh (\eta_n^2 v_n \otimes v_n) - \dz(\eta_n^2 w_n v_n) \\
		& ~~~~ - \nablah \eta_n^{2\gamma}
		 + \dvh(\eta_n^2\mathcal D (v_n))
		  + \dz (\eta_n^2 \dz v_n) + \dvh(\eta_n^2\nablah v_n) \\
		& ~~~~ + \eta_n \deltah \eta_n v_n
		+ \dvh (\eta_n \abs{\nablah \eta_n}{2} \nablah \eta_n \otimes v_n)
		- \abs{\nablah\eta_n}{4}v_n \\
		& ~~~~ - \eta_n^2 \abs{v_n}{3}v_n) \,ds \biggr\rbrack
		 + P_n (\eta_{n}^2 v_{n})|_{t=0} \biggr).
	\end{aligned}
\end{equation}
Here 
the map $ Q $ is mapping $ C([0,T];X_n) $ into $ C([0,T];X_n) $.
Moreover, for $$ \eta \in \lbrace \eta\in L^2(\Omega_h), \inf_{(x,y) \in \Omega_h} \eta > c \rbrace ~~ \text{with some positive constant $ c \in (0,\infty) $}, $$
the following map
$$ \eta \mapsto (M[\eta])^{-1} ~~ \text{is a map of  $ L^2(\Omega_h) $ into $ \mathcal L (X_n, X_n) $}. $$
It is well-defined as mentioned before, and it is Lipschitz, thanks to the identity
\begin{equation*}
	(M[\eta_1])^{-1} - (M[\eta_2])^{-1} = (M[\eta_1])^{-1}(M[\eta_2]-M[\eta_1])(M[\eta_2])^{-1}, ~ \eta_1,\eta_2 \in L^2(\Omega_h).
\end{equation*}
Together with \eqref{18July2018-03}, this implies that $ Q $ is contracting in $ C([0,T_n];X_n) $ for some sufficiently small $ T_n \in (0,T] $, and consequently, after applying the Banach fixed-point theorem to the finite dimension system \eqref{18July2018-02} (or equivalently \eqref{GA-009}) in the Banach space $ C([0,T_n];X_n) $, one can obtain that there is a local solution $ v_n $ to \eqref{GA-009} in the time interval $ [0,T_n] \subset [0,T] $.

Notice that $ ( \rho_n = \eta_n^2,v_n,w_n) $ satisfies the equation
\begin{equation*}
	\dt \eta_n^2 + 	\dvh (\eta_n^2 v_n) + \dz (\eta_n^2 w_n) = \eta_n \dvh (( 1 + \abs{\nablah \eta_n}{2}) \nablah \eta_n) + \eta_n^{-2p_0}.
\end{equation*}
After multiplying \eqref{GA-009} with $ a_i $ and summing the resultant, for $ i = 1,2 \cdots n $, one has the following equation, which is the same as \eqref{BE-001},
\begin{align*}
	& \dfrac{d}{dt} \biggl( \dfrac{1}{2} \int \eta_n^2 v_n^2 \idx \biggr) + \int \biggl( \eta_n^2 \abs{\mathcal D(v_n)}{2} + \eta_n^2 \abs{\dz v_n}{2} + \eta_n^2 \abs{\nablah v_n}{2} \biggr) \idx \\
	& ~~~~ + \int \nablah \eta_n^{2\gamma} \cdot v_n \idx + \dfrac{1}{2} \int \biggl( \eta_n \deltah \eta_n \abs{v_n}{2} + \eta_n \dvh (\abs{\nablah \eta_n}{2} \nablah \eta_n) \abs{v_n}{2} \\
	& ~~~~ + 2 \eta_n \abs{\nablah \eta_n}{2} \nablah \eta_n \cdot \nablah v_n \cdot v_n - \eta_n^{-2p_0} \abs{v_n}{2} \biggr) \idx.
\end{align*}
Moreover, $ v_n \in C([0,T_n]; X_n) $ is smooth and $ \eta_n , w_n $ is sufficiently regular. Therefore, one can perform similar arguments as in section \ref{sec:basic-energy} and obtain the bound of the approximating basic energy
\begin{align*}
	& \sup_{0\leq t\leq T_n} \bigl\lbrace \int \eta_n^2 v_n^2 \idx + \int  \eta_n^{2\gamma} \idx + \int \eta_n^2 \idx + \int \eta_n^{-2p_0}  \idx \bigr\rbrace \\
	& ~~~~ + \int_0^{T_n} \int \biggl( \eta_n^2 \abs{\mathcal D(v_n)}{2} + \eta_n^2 \abs{\dz v_n}{2}
	+ \eta_n^2 \abs{\nablah v_n}{2} \biggr) \idx \,dt \\
	& ~~~~ + \int_0^{T_n} \int \biggl( \eta_n^{-2p_0} \abs{v_n}{2} + \eta_n^2 \abs{v_n}{5} + \abs{\nablah \eta_n}{4}\abs{v_n}{2} + \abs{ \nablah \eta_n}{2} \abs{v_n}{2} \\
	& ~~~~ + \abs{\nablah \eta_n}{4} + \abs{\eta_{n}}{-4p_0-2} \biggr) \idx\,dt < C_T,
\end{align*}
independent of $ n $. Therefore, it follows that $ T_n = T $.
Now we further discuss the estimates of $ \eta_n $. Indeed, multiply \subeqref{eq:GA}{1} for $ (\eta, v) = (\eta_n, v_n) $ with  $ - \dvh (( 1+ \abs{\nablah \eta_n}{2}) \nablah \eta_n) $ and integrate the resultant. It holds after applying integration by parts, as before we have
\begin{align*}
	&  \sup_{0\leq t\leq T} \biggl\lbrace \int_{\Omega_h} \biggl( \abs{\nablah \eta_n}{2} + \dfrac{1}{2} \abs{\nablah \eta_n}{4}\biggr) \idxh   \biggr\rbrace + \dfrac{1}{2} \int_0^T \int_{\Omega_h} \biggl( \abs{\deltah \eta_n}{2} \\
	& ~~~~ + \abs{\nablah \eta_n}{2}\abs{\nablah^2 \eta_n}{2}  + \abs{\nablah \eta_n}{4} \abs{\nablah^2 \eta_n}{2} \biggr) \idxh \,dt \leq C \int_0^T \int_{\Omega} \biggl( \abs{\eta_n}{-4p_0-2} \\
	& ~~~~ +\abs{\nablah \eta_n}{2}\abs{v_n}{2}
	+ \abs{\eta_n}{2} \abs{ \nablah v}{2} \biggr) \idx \,dt < C.
\end{align*}
Then one can follow the same arguments as in section \ref{sec:ul-bounds} to get the upper and lower bounds for $ \eta_n $ as in \eqref{ene:density-UL}
. Therefore the bound in \eqref{GR-001} holds for $ (\eta, v, w) $ replaced by $ (\eta_n, v_n, w_n) $ and \eqref{ene:H^1-of-density} holds for $ \eta_n $. With such estimates, from \eqref{GA-009} one will get the bound
$$ \norm{\dt v_n}{L^2(0,T;W^{-1,1}(\Omega))} < C.  $$
Therefore, there exists a subsequence, denoted also by $ \lbrace(\eta_n, v_n) \rbrace $, that after taking $ n \rightarrow \infty $, one will get the weak solution $ (\eta, v) $ to the equation \subeqref{eq:GA}{2} and the strong solution $ \eta $ to the equation \subeqref{eq:GA}{1}. Indeed, after applying the Aubin's compactness Theorem (see, e.g., \cite[Theorem 2.1]{Temam1984} and \cite{Simon1986,Chen2012}), as $ n \rightarrow \infty $, we have the following:
\begin{align*}
	& \eta_{n} \rightarrow \eta_{}, ~~ & \text{in}& ~ L^2(0,T; W^{1,p}(\Omega_h)) \cap C(0,T; L^q(\Omega_h)  , \\
	& & & \text{for $(p,q) \in  (1,\infty) \times [2,\infty)$, } \\
	& \eta_{n} \buildrel\ast\over\rightharpoonup \eta, ~~  & \text{weak-$\ast$ in}& ~ L^\infty(0,T; W^{1,4}(\Omega_h)),\\
	& \eta_{n} \rightharpoonup \eta, ~~ & \text{weakly in}& ~ L^2(0,T;H^2(\Omega_h)),\\
	& \dt \eta_{n} \rightharpoonup \dt \eta, ~~ & \text{weakly in} & ~ L^2(\Omega_h\times(0,T)), \\
	& v_n \rightarrow v,  ~~ & \text{in} & ~  L^2(0,T; L^p(\Omega)), ~~
	 \text{for $ p \in (1,6) $ },\\
	& v_n \buildrel\ast\over\rightharpoonup v, ~~ & \text{weak-$\ast$ in} & ~ L^\infty(0,T;L^2(\Omega)) ,\\
	& v_n \rightharpoonup v, ~~ & \text{weakly in} &  ~ L^5(\Omega\times(0,T)) \cap L^2(0,T;H^1(\Omega)),
	\\
	& \dt v_n \rightharpoonup \dt v, ~~ & \text{weakly in} & ~ L^2(0,T;W^{-1,1}(\Omega)), \\
	& w_n \rightharpoonup w, ~~ & \text{weakly in} & ~ L^2(\Omega\times(0,T))   ,\\
	& \dz w_n \rightharpoonup \dz w, ~~ & \text{weakly in} & ~ L^2(\Omega\times(0,T))   ,
\end{align*}
where
\begin{equation}\label{18July2018-04}
\begin{aligned}
	& \eta \in L^\infty(0,T;W^{1,4}(\Omega_h) \cap L^2(0,T;H^2(\Omega_h)) , ~ \dt \eta \in L^2(\Omega_h\times(0,T)) ,\\
	& v \in L^\infty(0,T;L^2(\Omega))\cap L^2(0,T;H^1(\Omega)) \cap L^5(\Omega\times (0,T)), \\
	& \dt v \in L^2(0,T;W^{-1,1}(\Omega)), ~ w, \dz w \in L^2(\Omega \times (0,T)).
\end{aligned}
\end{equation}
Then arguing as in Proposition \ref{lm:the-map-S-approx} yields that $ \eta $ is the strong solution to \subeqref{eq:GA}{1} and admits positive upper and lower bounds. Similar arguments yield that $ (\eta, v, w) $ satisfies \subeqref{eq:GA}{3} in $ L^2(\Omega\times(0,T)) $.
Also, after applying the convergence above in \eqref{GA-009}, one has the following weak formulation of \subeqref{eq:GA}{2}, for every $ \phi = (\phi_1,\phi_2)^\top \in C_c^\infty(\overline\Omega\times [0,T)) $,
\begin{align*}
	& \int_0^T \int \biggl\lbrack v \cdot \eta^2 \dt \phi - \eta^2 \mathcal D(v) : \nablah \phi - \eta^2 \nablah v : \nablah \phi - \eta^2 \dz v \cdot \dz \phi \\
	& ~~~~ + \eta^2 v \otimes v : \nablah \phi  + \eta^2 w v \cdot \dz \phi + \eta^{2\gamma} \dvh \phi + \eta \deltah \eta v \cdot \phi \\
	& ~~~~ - (\eta\abs{\nablah \eta}{2} \nablah \eta \otimes v) : \nablah \phi - \abs{\nablah \eta}{4} v \cdot \phi - \eta^2 \abs{v}{3} v \cdot \phi \biggr\rbrack \idx \,dt \\
	& = - \int \eta^2 v \cdot \phi \idx\biggr|_{t=0}.
\end{align*}
By taking a positive approximating sequence of $ \eta $ in $ C_c^\infty(\overline \Omega\times[0,T]) $, the above formulation holds true for $ \phi = \dfrac{\psi}{\eta^2} $, with $ \psi \in C_c^\infty(\overline\Omega\times[0,T)) $. Then after making use of \subeqref{eq:GA}{1} and \subeqref{eq:GA}{3}, one will arrive at:
\begin{equation}\label{GA-011}
	\begin{aligned}
	& \int_0^T \int \biggl( v \cdot \dt \psi - \dfrac{3}{2} \nablah v : \nablah \psi - \dz v\cdot \dz \psi - \dfrac{1}{2} \dvh v \dvh \psi \biggr) \idx \,dt \\
	& ~~~~ + \int_0^T \int g \cdot \psi \idx\,dt = - \int v \cdot \psi \idx \biggr|_{t=0},
	\end{aligned}
\end{equation}
where $ g $ is as in \eqref{GR-010} with $ \varepsilon = 1 $. This procedure is rigorous since $ \eta > 0 $ admits upper and lower bounds and enough regularity. In particular, $ (\eta, v) $  is the weak solution to the following equation which is \eqref{GR-010} with $ \varepsilon = 1 $,
\begin{equation}\label{GA-010}
	\dt v - \dfrac{3}{2} \deltah v - \partial_{zz} v - \dfrac{1}{2} \nablah\dvh v = g,
\end{equation}
with $ g $ defined as in \eqref{approximate:v-source}.

On the other hand, we extend $ (\eta, v, w) $ to $ (\eta, v_e, w_e) $ in $ \Omega_h \times 2 \mathbb T $ as follows.
 Let $ v_e $ be a function satisfying:
 \begin{center}
	$ v_e = v $ in $ \Omega_h \times (0,1) $, and
	$ v_e $ is even in $ z \in (-1,1) $ and periodic in $ z \in 2 \mathbb T$.
\end{center}
Then, since $ \eta $ is independent of the $ z $-variable, $ w_e $ given by \subeqref{eq:GA}{3} with $ v $ replaced by $ v_e $ is odd in $ z \in (-1,1) $ and periodic in $ z \in 2 \mathbb T $. 
Then it is easy to verify that $ (\eta,v_e,w_e) $ is a weak solution to \eqref{eq:GA} in $ ( \Omega_h \times 2\mathbb T ) \times [0,T) $, and that $(\eta,v_e,w_e)=(\eta,v,w)$ in $ \Omega \times [0, T)$ satisfying the same regularity of $ (\eta, v, w) $ as in \eqref{18July2018-04}.
In particular, \eqref{GA-011} are satisfied by $ (\eta,v_e,w_e) $ and $ (\eta,v_e)\bigr|_{t=0} = (\rho^{1/2}_{\varepsilon,0}, v_{\varepsilon,0}) $, $ v_e $ is a weak solution to \eqref{GA-010}.


Then after employing the standard difference quotient method in \eqref{GA-011}, following similar estimates as in section \ref{sec:global-approx-priori-v} with the differential operators replaced by corresponding difference quotients,
%
one can get the following estimates, similar to \eqref{ene:H^1-of-v},
\begin{align*}
	& \norm{\nabla v_e}{L^\infty(0,T;L^2(\Omega_h\times 2\mathbb T))} + \norm{\nabla^2 v_e}{L^2((\Omega_h\times 2\mathbb T) \times(0,T))} \\
	& ~~~~ + \norm{\dt v_e}{L^2((\Omega_h\times 2\mathbb T) \times(0,T))} < C.
\end{align*}
This, together with \subeqref{eq:GA}{3}, implies that $ (\eta, v_e,w_e) $ satisfies equations \eqref{GA-010} and \eqref{eq:GA} almost everywhere in $ (\Omega_h\times 2\mathbb T)\times(0,T) $. In particular, the trace theorem implies that $$ v_e \bigr|_{z \in \mathbb Z} \in L^2(0,T; H^{3/2}(\Omega_h)), $$
and in particular, thanks to the symmetry property, $ \dz v_e \bigr|_{z \in \mathbb Z} \equiv 0 $.
Since $ v = v_e $ in $ \Omega $, we then have
\begin{gather*}
	 v \in L^\infty(0,T;H^1(\Omega)) \cap L^2(0,T;H^2(\Omega)), ~ \dt v \in L^2(\Omega\times(0,T)), \\
	 v \bigr|_{z =0,1} \in L^2(0,T; H^{3/2}(\Omega_h)), ~ \text{and} ~ \dz v \bigr|_{z =0,1} \equiv 0.
\end{gather*}
 This verifies the boundary condition \eqref{GA-BC}.

In particular, now one can employ the global a priori estimates in section \ref{sec:aprx-propri-est} to show the existence of the the approximating solutions satisfying Definition \ref{def:aprxm-sols}. This finishes the proof of Proposition \ref{prop:aprx-sol-exist}.

\section{Compactness and existence of weak solutions}\label{sec:compactness}

%
%

In this section, we aim at taking $ \varepsilon \rightarrow 0^+ $ in a subsequence of approximating solutions satisfying Definition \ref{def:aprxm-sols}. This will eventually establish the global existence of weak solutions to the compressible primitive equations \eqref{CPE'}. We will need some  compactness arguments in order to do so. Recall, from Proposition \ref{prop:uniform-est}, we have the following uniform estimates of the approximating solutions,
\begin{equation}\label{CM-001}
	\begin{aligned}
		& \sup_{0 \leq  t \leq T} \bigl\lbrace \norm{\rho_\varepsilon^{1/2}v_\varepsilon}{L^2} + \norm{\rho_\varepsilon}{L^\gamma} + \norm{\rho_\varepsilon}{L^1} + \norm{\varepsilon^{1/p_{0}} \rho_\varepsilon^{-1}}{L^{p_0}} \\
		& ~~~~ + \norm{{\nablah \rho_\varepsilon^{1/2}}{}}{L^2} + \norm{\varepsilon^{1/4}  \abs{\nablah \rho_\varepsilon^{1/2}}{}}{L^4} + \norm{\rho_\varepsilon(e+v_\varepsilon^2)\log(e+v_\varepsilon^2)}{L^1}
		\bigr\rbrace \\
		& ~~~~ + \norm{\rho_\varepsilon^{1/2} \nabla v_\varepsilon}{L^2(\Omega\times(0,T))} + \norm{\rho_\varepsilon^{1/2} \dz w_\varepsilon}{L^2(\Omega\times(0,T))} \\
		& ~~~~ + \norm{\rho_\varepsilon^{\gamma/2-1}\nablah \rho_\varepsilon }{L^2(\Omega\times(0,T))} + \norm{\varepsilon^{1/2}\rho_\varepsilon^{-p_0/2}v_\varepsilon}{L^2(\Omega\times(0,T))} \\
		& ~~~~ + \norm{\varepsilon^{1/5}\rho_\varepsilon^{1/5} v_\varepsilon}{L^5(\Omega\times(0,T))} + \norm{\varepsilon^{1/2}\abs{\nablah\rho_\varepsilon^{1/2}}{2}v_\varepsilon}{L^2(\Omega\times(0,T))}\\
		& ~~~~ + \norm{\varepsilon^{1/2}\abs{\nablah\rho_\varepsilon^{1/2}}{}v_\varepsilon}{L^2(\Omega\times(0,T))} + \norm{\varepsilon \rho_\varepsilon^{-p_0-1/2}}{L^2\spacetime} \\
		& ~~~~ + \norm{\varepsilon^{1/2} \deltah \rho_\varepsilon^{1/2}}{L^2\spacetime} + \norm{\varepsilon^{1/2}\abs{\nablah \rho^{1/2}_\varepsilon}{} \abs{\nablah^2 \rho_\varepsilon^{1/2}}{}}{L^2\spacetime} \\
		& ~~~~ + \norm{\varepsilon\abs{\nablah \rho^{1/2}_\varepsilon}{2}\abs{\nablah^2 \rho_{\varepsilon}^{1/2}}{}}{L^2\spacetime} < C_T,
	\end{aligned}
\end{equation}
where $ C_T $ is independent of $ \varepsilon $.
Since $ \rho_\varepsilon $ is independent of the $ z $-variable, it follows from the Minkowski inequality that,
\begin{equation}\label{CM-001-1}
	\begin{aligned}
		& \norm{\rho_\varepsilon^{1/2}w_\varepsilon}{L^2\spacetime}= \norm{\int_0^z \rho_\varepsilon^{1/2}\dz w_\varepsilon \,dz'}{L^2\spacetime} \\
		& ~~~~ ~~~~ \leq \norm{\rho_\varepsilon^{1/2}\dz w_\varepsilon}{L^2\spacetime}.
	\end{aligned}
\end{equation}
In the following subsections, we will first discuss the compactness based on the uniform estimates \eqref{CM-001} and then establish the existence of weak solutions. Through the rest of this section, $ C_T  $ denotes a positive constant independent of $ \varepsilon $. Without further mentioned, all the convergences, below, hold in the sense of a subsequence of $ \varepsilon $, as $ \varepsilon \rightarrow 0^+ $. 

\subsection{Strong convergences}

\subsubsection*{Strong convergence of $ \rho_\varepsilon $}

By applying H\"older's inequality, one has, by virtue of \eqref{CM-001}, for every $ t \in [0, T ] $,
\begin{align*}
	& \norm{\nabla\rho_\varepsilon}{L^{2\gamma/(\gamma+1)}}  \leq \norm{\rho_\varepsilon^{1/2}}{L^{2\gamma}} \norm{\rho_\varepsilon^{-1/2}\abs{\nablah \rho_\varepsilon}{}}{L^2}\\
	& ~~~~ = 2 \norm{\rho_\varepsilon}{L^\gamma}^{1/2} \norm{\nablah \rho_\varepsilon^{1/2}}{L^2} \leq C_T,
\end{align*}
where $ C_T $ is independent of  $ \varepsilon $.
In particular, after applying the Gagliardo-Nirenberg inequality, one has,
\begin{equation}\label{CM-002}
	\norm{\rho_\varepsilon}{L^\infty(0,T;W^{1,2\gamma/(\gamma+1)}(\Omega))} < C_T.
\end{equation}
On the other hand,
from \subeqref{eq:approximating-CPE}{1}, we have
\begin{equation}\label{eq:appro-density}
	\begin{aligned}
		& \dt \rho_\varepsilon + \dvh (\rho_\varepsilon v_\varepsilon) + \dz (\rho_\varepsilon w_\varepsilon) = \varepsilon \rho_\varepsilon^{1/2} \deltah \rho_\varepsilon^{1/2}\\
		& ~~~~ + \varepsilon \dvh (\rho_\varepsilon^{1/2} \abs{\nablah \rho_\varepsilon^{1/2}}{2} \nablah \rho_\varepsilon^{1/2} )  - \varepsilon \abs{\nablah \rho_\varepsilon^{1/2}}{4} + \varepsilon \rho_\varepsilon^{-p_0}.
	 \end{aligned}
\end{equation}
Notice that after making use of \eqref{CM-001} and \eqref{CM-001-1},
\begin{align}
	& \sup_{0\leq t\leq T} \norm{\rho_\varepsilon v_\varepsilon}{L^1} \leq \sup_{0\leq t\leq T} \norm{\rho^{1/2}_\varepsilon v_\varepsilon}{L^2} \sup_{0\leq t\leq T} \norm{\rho_\varepsilon^{1/2}}{L^2} {\nonumber} \\
	& ~~~~ = \sup_{0\leq t\leq T} \norm{\rho^{1/2}_\varepsilon v_\varepsilon}{L^2} \sup_{0\leq t\leq T} \norm{\rho_\varepsilon}{L^1}^{1/2} < C_T,\\
	& \norm{\rho_\varepsilon w_\varepsilon}{L^1\spacetime} \leq T^{1/2} \sup_{0\leq t\leq T} \norm{\rho_\varepsilon^{1/2}}{L^2} \times  \norm{\rho_\varepsilon^{1/2}w_\varepsilon}{L^2\spacetime} {\nonumber}\\
	& ~~~~ \leq C T^{1/2} \sup_{0\leq t\leq T} \norm{\rho_\varepsilon}{L^1}^{1/2} \times  \norm{\rho^{1/2}_\varepsilon\dz w_\varepsilon}{L^2\spacetime} <  C_T, \\
	& \label{CM-101} \norm{\varepsilon \rho_\varepsilon^{1/2}\deltah \rho_\varepsilon^{1/2}}{L^1\spacetime} \leq \varepsilon^{1/2} T^{1/2} \sup_{0\leq t\leq T} \norm{\rho_\varepsilon^{1/2}}{L^2} \nonumber\\
	& ~~~~ ~~~~ \times  \norm{\varepsilon^{1/2}\deltah \rho_\varepsilon^{1/2}}{L^2\spacetime} < \varepsilon^{1/2} C_T, \\
	& \norm{\varepsilon \rho_\varepsilon^{1/2} \abs{\nablah\rho_\varepsilon^{1/2}}{2}\nablah\rho_\varepsilon^{1/2}}{L^1\spacetime} \leq \varepsilon \int_0^T \sup_{0\leq t\leq T} \norm{\rho_\varepsilon^{1/2}}{L^2} \nonumber \\
	& ~~~~ ~~~~ \times \norm{\abs{\nablah \rho_\varepsilon^{1/2}}{3}}{L^2} \,dt 
	< \varepsilon^{1/6} C_T,\\
	& \norm{\varepsilon \abs{\nablah \rho_\varepsilon^{1/2}}{4}}{L^1\spacetime} \leq \varepsilon \int_0^T \sup_{0\leq t\leq T} \norm{{\nablah \rho_\varepsilon^{1/2}}{}}{L^2} {\nonumber}\\
	& ~~~~ ~~~~ \times  \norm{\abs{\nablah \rho_\varepsilon^{1/2}}{3}}{L^2}^{} \,dt  < \varepsilon^{1/6} C_T, \\
	& \label{CM-102} \norm{\varepsilon \rho_\varepsilon^{-p_0}}{L^1 \spacetime} \leq \varepsilon^{1/(2p_0+1)} \norm{\varepsilon \rho_\varepsilon^{-p_0-1/2}}{L^2\spacetime}^{p_0/(p_0+1/2)} \nonumber \\
	& ~~~~ < \varepsilon^{1/(2p_0+1)} C_T,
\end{align}
where in the above estimates we have noticed the facts that $ \rho_\varepsilon $ is independent of the $ z $ variable and that $ \varepsilon \in (0,\varepsilon_0) $, with $ \varepsilon_0\in(0,1) $ which is small enough, and applied the H\"older and Minkowski inequalities and the two-dimensional Gagliardo-Nirenberg inequality,
\begin{align}
	& \int_0^T \norm{\abs{\nablah \rho_\varepsilon^{1/2}}{3}}{L^2} \,dt = \int_0^T \normh{\abs{\nablah \rho_\varepsilon^{1/2}}{3}}{L^2} \,dt {\nonumber} \\
	& ~~~~ \leq C \int_0^T \normh{
	\abs{\nablah \rho_\varepsilon^{1/2}}{3}}{L^{4/3}}^{2/3} \normh{\abs{\nablah \rho_\varepsilon^{1/2}}{2} \abs{\nablah^2 \rho_\varepsilon^{1/2}}{}}{L^2}^{1/3}  + \normh{
	\abs{\nablah \rho_\varepsilon^{1/2}}{3}}{L^{4/3}} \,dt \nonumber\\
	& ~~~~ = C \int_0^T \norm{\nablah \rho_\varepsilon^{1/2}}{L^4}^{2} \norm{\abs{\nablah \rho_\varepsilon^{1/2}}{2} \abs{\nablah^2\rho_\varepsilon^{1/2}}{}}{L^2}^{1/3}  + \norm{{\nablah \rho_\varepsilon^{1/2}}{}}{L^{4}}^{3} \,dt {\nonumber}\\
	& ~~~~ \leq C T^{5/6}  \sup_{0\leq t\leq T} \norm{\nablah \rho_\varepsilon^{1/2}}{L^4}^{2} \times \norm{\abs{\nablah \rho_\varepsilon^{1/2}}{2} \abs{\nablah^2\rho_\varepsilon^{1/2}}{}}{L^2\spacetime}^{1/3} {\nonumber}\\
	& ~~~~ ~~~~ + C T \sup_{0\leq t\leq T} \norm{\nablah \rho_\varepsilon^{1/2}}{L^{4}}^{3} \leq \varepsilon^{-1/2-1/3} C_T + \varepsilon^{-3/4} C_T 
	\leq \varepsilon^{-5/6}C_T.  \label{CM-004}
\end{align}
Therefore, from \eqref{eq:appro-density} we have
\begin{equation}\label{CM-003}
	\norm{\dt\rho_\varepsilon}{L^1(0,T;W^{-1,1}(\Omega))} \leq C_T.
\end{equation}
Thus together with \eqref{CM-002}, applying the Aubin's compactness Theorem (see, e.g., \cite[Theorem 2.1]{Temam1984} and \cite{Simon1986,Chen2012}) implies that there is a $ \rho \in L^p(\Omega\times(0,T)) $,
\begin{equation}\label{CM-density}
	\rho_\varepsilon \rightarrow \rho ~~~~~~ \text{in}~~ L^p\spacetime,
\end{equation}
for any $ 1 < p < 2\gamma $, as $ \varepsilon \rightarrow 0^+ $, up to a subsequence. Also we have
\begin{equation}\label{CM-density-2}
	\rho \in L^{\infty}(0,T;L^1(\Omega)\cap L^{\gamma}(\Omega)).
\end{equation}
In particular, as $ \varepsilon \rightarrow 0^+ $,
\begin{equation} \label{CM-density-ae}
	\rho_\varepsilon \rightarrow \rho, ~~~~ \text{almost everywhere in} ~ \Omega\times(0,T).
\end{equation}
On the other hand, from \eqref{CM-001}, for some positive constant independent of $\varepsilon $,
\begin{equation*}
	\sup_{0\leq t\leq T} \bigl\lbrace \norm{ \rho_\varepsilon^{1/2}}{L^2} +  \norm{\nablah \rho_\varepsilon^{1/2}}{L^2} \bigr\rbrace < C_T .
\end{equation*}
Therefore, one has, as $ \varepsilon \rightarrow 0^+ $,
\begin{equation}\label{CM-density-star}
	\rho_\varepsilon^{1/2} \stackrel{\ast}\rightharpoonup \rho^{1/2} ~~~~~~ \text{weak-$\ast$ in} ~~ L^\infty(0,T;W^{1,2}(\Omega)).
\end{equation}
Thus we have shown the following:
\begin{prop}
	As $ \varepsilon \rightarrow 0^+ $, taking a subsequence if necessary, the following convergences hold, for $ 1 < p < 2\gamma $,
	\begin{gather*}
		\rho_\varepsilon \rightarrow \rho ~~~~~~ \text{in}~~ L^p\spacetime, \tag{\ref{CM-density}}\\
		\tag{\ref{CM-density-ae}}
		\rho_\varepsilon \rightarrow \rho ~~~~ \text{almost everywhere in} ~ \Omega\times(0,T),  \\
		\tag{\ref{CM-density-star}}
		\rho_\varepsilon^{1/2} \stackrel{\ast}\rightharpoonup \rho^{1/2} ~~~~~~ \text{weak-$\ast$ in} ~~ L^\infty(0,T;W^{1,2}(\Omega)).
	\end{gather*}
\end{prop}

\subsubsection*{Strong convergence of $ \rho_\varepsilon v_\varepsilon $}
It follows from \eqref{CM-001} that,
	\begin{align*}
		& \sup_{0\leq t\leq T} \norm{\rho_\varepsilon v_\varepsilon}{L^1}\leq \sup_{0\leq t\leq T} ( \norm{\rho_\varepsilon^{1/2}}{L^2} \norm{\rho_\varepsilon^{1/2} v_\varepsilon}{L^2} ) \\
		& ~~~~ = \sup_{0\leq t\leq T} ( \norm{\rho_\varepsilon}{L^1}^{1/2} \norm{\rho_\varepsilon^{1/2} v_\varepsilon}{L^2}) < C_T,\\
		& \norm{\nabla (\rho_\varepsilon v_\varepsilon)}{L^2(0,T,L^1(\Omega))} \leq \norm{\nablah \rho_\varepsilon v_\varepsilon}{L^2(0,T,L^1(\Omega))} + \norm{\rho_\varepsilon \nabla v_\varepsilon}{L^2(0,T,L^1(\Omega))} \\
		& ~~~~ \leq T^{1/2} C \sup_{0\leq t\leq T} (\norm{\nablah \rho_\varepsilon^{1/2}}{L^2} \norm{\rho_\varepsilon^{1/2} v_\varepsilon}{L^2}) + C \sup_{0\leq t\leq T} \norm{\rho_\varepsilon^{1/2}}{L^2} \\
		& ~~~~ ~~~~ \times \norm{\rho_\varepsilon^{1/2}\nabla v_\varepsilon}{L^2\spacetime}
		 < C_T.
	\end{align*}
Hence,
\begin{equation}\label{CM-005}
	\norm{\rho_\varepsilon v_\varepsilon}{L^2(0,T;W^{1,1}(\Omega))} < C_T.
\end{equation}
To obtain the estimate of $ \dt(\rho_\varepsilon v_\varepsilon) $, from \eqref{eq:approximating-CPE}, we have
\begin{equation}\label{eq:approx-momentum}
	\begin{aligned}
		& \dt (\rho_\varepsilon v_\varepsilon ) + \dvh (\rho_\varepsilon v_\varepsilon \otimes v_\varepsilon) + \dz(\rho_\varepsilon w_\varepsilon v_\varepsilon) + \nablah \rho^\gamma_\varepsilon - \dvh(\rho_\varepsilon\mathcal D (v_\varepsilon)) \\
		& ~~~~ - \dz (\rho_\varepsilon \dz v_\varepsilon)
		= \varepsilon \rho_\varepsilon^{1/2}\deltah \rho_\varepsilon^{1/2} v_\varepsilon + \varepsilon \dvh (\rho^{1/2}_\varepsilon \abs{\nablah \rho^{1/2}_\varepsilon}{2} \nablah \rho_\varepsilon^{1/2} \otimes v_\varepsilon) \\
		& ~~~~ - \varepsilon \abs{\nablah\rho_\varepsilon^{1/2}}{4}v_\varepsilon
		+ \sqrt{\varepsilon} \dvh(\rho_\varepsilon\nablah v_\varepsilon) - \varepsilon \rho_\varepsilon \abs{v_\varepsilon}{3}v_\varepsilon.
	\end{aligned}
\end{equation}
Moreover, from \eqref{CM-001}, \eqref{CM-001-1}, \eqref{CM-004}, we have
\begin{align}
	& \norm{\rho_\varepsilon v_\varepsilon\otimes v_\varepsilon}{L^1\spacetime} \leq T \sup_{0\leq t\leq T} \norm{\rho_\varepsilon^{1/2}v_\varepsilon}{L^2}^2 < C_T,\\
	& \norm{\rho_\varepsilon v_\varepsilon w_\varepsilon}{L^1\spacetime} \leq T^{1/2} \sup_{0\leq t\leq T}\norm{\rho_\varepsilon^{1/2}v_\varepsilon}{L^2} {\nonumber}\\
	& ~~~~ ~~~~ \times \norm{\rho_\varepsilon^{1/2} \dz w_\varepsilon}{L^2\spacetime} < C_T, \\
	& \norm{\rho_\varepsilon^\gamma}{L^1\spacetime} \leq T \sup_{0\leq t\leq T} \norm{\rho_\varepsilon}{L^\gamma}^\gamma < C_T, \\
	& {\nonumber} \norm{\rho_\varepsilon \mathcal D(v_\varepsilon)}{L^1\spacetime} + \norm{\rho_\varepsilon \dz v_\varepsilon}{L^1\spacetime} \leq \norm{\rho_\varepsilon \nabla v_\varepsilon}{L^1\spacetime} \\
	& ~~~~\leq T^{1/2} \sup_{0\leq t\leq T} \norm{\rho_\varepsilon^{1/2}}{L^2} \times \norm{\rho_\varepsilon^{1/2} \nabla v_\varepsilon}{L^2\spacetime} < C_T,\\
	& \label{CM-201} \norm{\sqrt{\varepsilon} \rho_\varepsilon \nablah v_\varepsilon}{L^1\spacetime} < \sqrt{\varepsilon} C_T , \\
	& {\nonumber} \norm{\varepsilon \rho^{1/2}_\varepsilon\abs{\nablah\rho_\varepsilon^{1/2}}{2}\nablah \rho_\varepsilon^{1/2}\otimes v_\varepsilon}{L^1\spacetime} \leq \varepsilon C \sup_{0\leq t\leq T} \norm{\rho_\varepsilon^{1/2} v_\varepsilon}{L^2} {\nonumber}\\
	& ~~~~~ ~~~~ \times \int_0^T \norm{\abs{\nablah\rho_\varepsilon^{1/2}}{3}}{L^2}\,dt \leq \varepsilon^{1/6} C_T.
\end{align}
On the other hand, thanks to \eqref{CM-001}, we have
\begin{align}
	& \nonumber \norm{\varepsilon \rho_\varepsilon^{1/2} \deltah \rho_\varepsilon^{1/2} v_\varepsilon}{L^1\spacetime} \leq \varepsilon^{1/2} T^{1/2} \sup_{0\leq t\leq T} \norm{\rho_\varepsilon^{1/2}v_\varepsilon}{L^2} \\
	& ~~~~ ~~~~ \times \norm{\varepsilon^{1/2}\deltah\rho_\varepsilon^{1/2}}{L^2\spacetime} < \varepsilon^{1/2} C_T,\\
	& {\nonumber} \norm{\varepsilon \abs{\nablah \rho_\varepsilon^{1/2}}{4} v_\varepsilon}{L^1\spacetime} \leq\norm{\varepsilon^{1/2} \abs{\nablah \rho_\varepsilon^{1/2}}{2}}{L^2\spacetime}\\
	& ~~~~ ~~~~ \times \norm{\varepsilon^{1/2} \abs{\nablah\rho_\varepsilon^{1/2}}{2}v_\varepsilon}{L^2\spacetime} \leq \varepsilon^{1/2} \sup_{0\leq t\leq T} \norm{\nablah\rho_\varepsilon^{1/2}}{L^2}^{1/2} \nonumber\\
	& ~~~~ ~~~~ \times \bigl(\int_0^T \norm{\abs{\nablah \rho_\varepsilon^{1/2}}{3}}{L^2}\,dt\bigr)^{1/2} \times \norm{\varepsilon^{1/2} \abs{\nablah\rho_\varepsilon^{1/2}}{2}v_\varepsilon}{L^2\spacetime} {\nonumber} \\
	& ~~~~ \leq \varepsilon^{1/2-5/12}C_T = \varepsilon^{1/12}C_T,\\
	& {\nonumber}\norm{\varepsilon\rho_\varepsilon \abs{v_\varepsilon}{3}v_\varepsilon}{L^1\spacetime} \leq \varepsilon^{1/5} \norm{\varepsilon^{1/5}\rho_\varepsilon^{1/5} v_\varepsilon}{L^5\spacetime}^4 \norm{\rho_\varepsilon^{1/5}}{L^5\spacetime} \\
	& {\label{CM-202}} ~~~~ \leq  \varepsilon^{1/5} T^{1/5}\norm{\varepsilon^{1/5}\rho_\varepsilon^{1/5} v_\varepsilon}{L^5\spacetime}^4 \times \sup_{0\leq t\leq T} \norm{\rho_\varepsilon}{L^1\spacetime} \nonumber \\
	& ~~~~ < \varepsilon^{1/5} C_T.
\end{align}
Therefore, we have, for some positive constant $ C_T $ independent of $ \varepsilon $,
\begin{equation}\label{CM-006}
	\norm{\dt (\rho_\varepsilon v_\varepsilon)}{L^1(0,T;W^{-1,1}(\Omega))} < C_T.	
\end{equation}
Together with \eqref{CM-005}, after applying the Aubin's compactness Theorem (see, e.g., \cite[Theorem 2.1]{Temam1984} and \cite{Simon1986,Chen2012}), this shows that there is $ m \in L^2(0,T;L^{3/2}(\Omega)) $ such that, as $ \varepsilon \rightarrow 0^+ $, up to a subsequence,
\begin{equation}\label{CM-007}
	\rho_\varepsilon v_\varepsilon \rightarrow m ~~~~~~ \text{in} ~~ L^2(0,T;L^p(\Omega)),
\end{equation}
for any $ p \in [1,3/2) $.
Define the horizontal velocity $ v $ as follows,
\begin{equation}\label{CM-def-horizontal-velocity}
	v(x,y,z,t) := \begin{cases}
		m(x,y,z,t)/\rho(x,y,t) & \text{in} ~~ \lbrace \rho > 0\rbrace,\\
		0 & \text{in} ~~ \lbrace \rho = 0 \rbrace.
	\end{cases}
\end{equation}
Then
\begin{equation}\label{CM-008}
	\rho_\varepsilon v_\varepsilon \rightarrow \rho v ~~~~ \text{almost everywhere in}~ \Omega\times(0,T).
\end{equation}
Indeed, since \eqref{CM-007} implies $ \rho_\varepsilon v_\varepsilon \rightarrow m $ almost everywhere in $ \Omega\times (0,T) $, $ \rho_\varepsilon v_\varepsilon \rightarrow \rho v $ almost everywhere on the set $ \lbrace \rho > 0 \rbrace $ as the direct consequence of definition \eqref{CM-def-horizontal-velocity}.
On the other hand, over the set $ M_0 : = \lbrace \rho = 0 \rbrace \subset \Omega\times[0,T] $,
\begin{align*}
	& \norm{\rho_\varepsilon v_\varepsilon}{L^1(M_0)} \leq \norm{\rho_\varepsilon^{1/2}}{L^2(M_0)}\norm{\rho_\varepsilon^{1/2} v_\varepsilon}{L^2(M_0)} \\
	& ~~~~ \leq T^{1/2} \sup_{0\leq t\leq T} \norm{\rho_\varepsilon^{1/2} v_\varepsilon}{L^2} \times \norm{\rho_\varepsilon}{L^1(M_0)}^{1/2} < C_T \norm{\rho_\varepsilon}{L^1(M_0)}^{1/2} \rightarrow 0,
\end{align*}
as $ \varepsilon \rightarrow 0^+ $. Here we have applied the Lebesgue dominated convergence theorem, which follows from
the facts that $ \rho_\varepsilon \rightarrow 0 $ almost everywhere in $ M_0 $ from \eqref{CM-density-ae} and that $ \norm{\rho_\varepsilon}{L^1(M_0)} \leq T \sup_{0\leq t\leq T}\norm{\rho}{L^1} < C_T $.
Therefore $ \rho_\varepsilon v_\varepsilon \rightarrow 0 = \rho v $ almost everywhere in $ M_0 $.


Thus, we have shown the following:
\begin{prop}
As $ \varepsilon \rightarrow 0^+ $, taking a subsequence if necessary, the following convergences hold,
\begin{gather*}
	\tag{\ref{CM-008}}
	\rho_\varepsilon v_\varepsilon \rightarrow \rho v, ~~~~ \text{almost everywhere in} ~ \Omega\times(0,T),
\end{gather*}	
where $ v $ is defined in \eqref{CM-007} and \eqref{CM-def-horizontal-velocity}.
\end{prop}

\subsubsection*{Strong convergence of $ \sqrt{\rho_\varepsilon} v_\varepsilon $}
Notice that from the definition of $ v $ in \eqref{CM-def-horizontal-velocity} and \eqref{CM-density-ae}, for almost every $ (x,y,z,t) \in \lbrace \rho > 0 \rbrace \subset \Omega\times(0,T) $, we have
\begin{align*}
	& \rho \abs{v}{2} \log (e+ \abs{v}{2}) = \dfrac{\abs{m}{2} \log (e \rho^2 + \abs{m}{2}) - 2 \abs{m}{2} \log \rho}{\rho} \\
	& ~~~~ = \lim\limits_{\varepsilon\rightarrow 0^+} \dfrac{\abs{\rho_\varepsilon v_\varepsilon}{2} \log (e \rho_\varepsilon^2 + \abs{\rho_\varepsilon v_\varepsilon}{2}) - 2 \abs{\rho_\varepsilon v_\varepsilon}{2} \log \rho_\varepsilon}{\rho_\varepsilon} \\
	& ~~~~ = \lim\limits_{\varepsilon \rightarrow 0^+} \rho_\varepsilon \abs{v_\varepsilon}{2}\log(e+\abs{v_\varepsilon}{2}),
\end{align*}
where on the right-hand side we have applied the fact that $ \rho_\varepsilon > 0 ~\text{in} ~ \Omega\times(0,T) $ as a consequence of \eqref{ene:density-UL}, for any $ \varepsilon \in (0,\varepsilon_0) $, and $ \varepsilon_0\in(0,1) $ is sufficiently small. On the other hand, over the set $ \lbrace \rho = 0 \rbrace $, we have
\begin{equation*}
	\rho \abs{v}{2} \log (e+ \abs{v}{2}) = 0 \leq \rho_\varepsilon \abs{v_\varepsilon}{2}\log(e+\abs{v_\varepsilon}{2}),
\end{equation*}
for any $ \varepsilon \in (0,\varepsilon_0) $. Therefore, from \eqref{CM-001} and Fatou's Lemma, we have
\begin{equation}\label{CM-009}
	\begin{aligned}
		& \norm{\rho \abs{v}{2} \log (e+ \abs{v}{2})}{L^1\spacetime} \leq  \norm{\liminf\limits_{\varepsilon\rightarrow 0^+}\rho_\varepsilon \abs{v_\varepsilon}{2}\log(e+\abs{v_\varepsilon}{2})}{L^1\spacetime} \\
		& ~~~~ \leq \liminf\limits_{\varepsilon\rightarrow 0^+} \norm{\rho_\varepsilon \abs{v_\varepsilon}{2}\log(e+\abs{v_\varepsilon}{2})}{L^1\spacetime} < C_T.
	\end{aligned}
\end{equation}

Next, for any $ M > 0 $, we have,
\begin{equation}\label{CM-012}
	\begin{aligned}
		& \norm{\rho_\varepsilon^{1/2} v_\varepsilon - \rho^{1/2} v}{L^2\spacetime}^2 \leq C \int_0^T \int_\Omega \bigl| \rho_\varepsilon^{1/2} v_\varepsilon \mathbbm{1}_{\lbrace |{v_\varepsilon}|{}\leq M \rbrace } \\
		& ~~~~ ~~~~ - \rho^{1/2} v \mathbbm{1}_{\lbrace |v|\leq M \rbrace} \bigr|^{2} \idx \,dt  + C \int_0^T \int_\Omega \abs{\rho_\varepsilon^{1/2} v_\varepsilon \mathbbm{1}_{\lbrace |{v_\varepsilon}|{} > M \rbrace } }{2} \\
		& ~~~~ + \abs{\rho^{1/2} v \mathbbm{1}_{\lbrace |v| > M \rbrace}}{2}   \idx\,dt =: I_1 + I_2.
	\end{aligned}
\end{equation}
First, we show that $ I_2 \rightarrow 0 $, as $ M \rightarrow \infty $, uniformly in $ \varepsilon $. 
Indeed, from \eqref{CM-001}, \eqref{CM-009}, as $ M \rightarrow \infty $, one has
\begin{equation}\label{CM-010}
	\begin{aligned}
		& I_2 \leq \dfrac{1}{\log (e+M^2)} \bigl( \norm{\rho_\varepsilon \abs{v_\varepsilon}{2} \log (e+ \abs{v_\varepsilon}{2})}{L^1\spacetime} \\
		& ~~~~ + \norm{\rho \abs{v}{2} \log (e+\abs{v}{2}) }{L^1\spacetime} \bigr) \leq \dfrac{C_T}{\log (e+M^2)} \rightarrow 0.
	\end{aligned}
\end{equation}
What is left is to show that for any fixed $  M \in (0, \infty) $, $ I_1 \rightarrow 0 $, as $ \varepsilon \rightarrow 0^+ $. Consider any number $  \delta \in (0, \infty) $. Denote the set $ M_\delta : = \lbrace \rho > \delta \rbrace \subset \Omega\times (0,T) $. Obviously, $ M_\delta $ is of finite measure. From \eqref{CM-density-ae}, Egorov's Theorem implies that there is a subset $ Q_\delta \subset M_\delta $, with $ \abs{M_\delta \backslash Q_\delta}{} < \delta $, such that $ \rho_\varepsilon \rightarrow \rho $ uniformly in $ Q_\delta $, as $ \varepsilon \rightarrow 0^+ $. Hence, there is a $ \varepsilon_1 \in (0, \varepsilon_0) $ such that $ \rho_\varepsilon > \delta / 2 $ in $ Q_\delta $ for any $ \varepsilon \in (0, \varepsilon_1) $. Then we have, from \eqref{CM-density-ae} and \eqref{CM-008}, that, as $ \varepsilon \rightarrow 0^+ $,
\begin{equation*}
		v_\varepsilon = \dfrac{\rho_\varepsilon v_\varepsilon}{\rho_\varepsilon} \rightarrow \dfrac{\rho v}{\rho} = v, ~~~~ \text{almost everywhere} ~~ (x,y,z,t) \in Q_\delta.
\end{equation*}
Again, together with \eqref{CM-density-ae}, as $ \varepsilon \rightarrow 0^+ $,
\begin{equation*}
	\rho_\varepsilon^{1/2} v_\varepsilon \rightarrow \rho^{1/2} v, ~~~~ \text{almost everywhere} ~~ (x,y,z,t) \in Q_\delta.
\end{equation*}
Thus, \eqref{CM-001} and the dominated convergence theorem imply
\begin{equation*}
	\int_0^T \int_\Omega \mathbbm{1}_{Q_\delta} \abs{ \rho_\varepsilon^{1/2} v_\varepsilon \mathbbm{1}_{\lbrace |{v_\varepsilon}|{}\leq M \rbrace } - \rho^{1/2} v \mathbbm{1}_{\lbrace |v|\leq M \rbrace} }{2} \idx \,dt \rightarrow 0 , ~~ \text{as} ~ \varepsilon \rightarrow 0^+.
\end{equation*}
On the other hand, as $ \varepsilon \rightarrow 0^+ $,
\begin{equation*}
	\begin{aligned}
		& \int_0^T \int_\Omega \mathbbm{1}_{M_\delta \backslash Q_\delta} \abs{ \rho_\varepsilon^{1/2} v_\varepsilon \mathbbm{1}_{\lbrace |{v_\varepsilon}|{}\leq M \rbrace } - \rho^{1/2} v \mathbbm{1}_{\lbrace |v|\leq M \rbrace} }{2} \idx \,dt \\
		& ~~~~ \leq M^2 C \int_0^T (\norm{\rho_\varepsilon}{L^\gamma}+ \norm{\rho}{L^\gamma} )\abs{M_\delta\backslash Q_\delta}{1-1/\gamma} \,dt \leq \delta^{1-1/\gamma} M^2 C_T,\\
		& \int_0^T \int_\Omega \mathbbm{1}_{ M_\delta^c 
		}\abs{ \rho_\varepsilon^{1/2} v_\varepsilon \mathbbm{1}_{\lbrace |{v_\varepsilon}|{}\leq M \rbrace } - \rho^{1/2} v \mathbbm{1}_{\lbrace |v|\leq M \rbrace} }{2} \idx \,dt  \\
		& ~~~~ \leq  \int_0^T \int_\Omega \mathbbm{1}_{\lbrace \rho\leq \delta \rbrace}M^2 \abs{ \rho_\varepsilon - \rho }{} + C \mathbbm{1}_{\lbrace \rho \leq \delta \rbrace } M^2 \abs{\rho}{} \idx \,dt \\
		& ~~~~ \leq M^2 \norm{\rho_\varepsilon - \rho}{L^1\spacetime} +  \delta M^2 C_T \rightarrow \delta M^2 C_T.
	\end{aligned}
\end{equation*}
Consequently, for any $ \delta > 0 $, we have
\begin{equation}\label{CM-011}
	\begin{aligned}
		& \limsup\limits_{\varepsilon\rightarrow 0^+} I_1 = \limsup\limits_{\varepsilon\rightarrow 0^+} \int_0^T \int_\Omega (\mathbbm{1}_{Q_\delta} + \mathbbm{1}_{M_\delta\backslash Q_\delta} + \mathbbm{1}_{M_\delta^c})  \bigl|\rho_\varepsilon^{1/2} v_\varepsilon \mathbbm{1}_{\lbrace |{v_\varepsilon}|{}\leq M \rbrace } \\
		& ~~~~ ~~~~ - \rho^{1/2} v \mathbbm{1}_{\lbrace |v|\leq M \rbrace} \bigr|^{2} \idx \,dt \\
		& ~~~~  \leq \delta^{1-1/\gamma} M^2C_T + \delta M^2 C_T, ~~ \text{as} ~ \varepsilon \rightarrow 0^+.
	\end{aligned}
\end{equation}
Therefore, by taking $ \delta \rightarrow 0^+ $, we conclude that $ I_1 \rightarrow 0 $, as $ \varepsilon \rightarrow 0^+ $, for any fixed $ M \in (0,\infty) $. Together with \eqref{CM-012} and \eqref{CM-010}, we conclude with the following:
\begin{prop}
As $ \varepsilon \rightarrow 0^+ $, taking a subsequence if necessary,
\begin{equation}\label{CM-energy}
	\rho_\varepsilon^{1/2}v_\varepsilon \rightarrow \rho^{1/2} v ~~~~~~ \text{in} ~~ L^2 \spacetime,
\end{equation}
with
\begin{equation}\label{CM-energy-2}
	\rho^{1/2} v \in L^\infty(0,T;L^2(\Omega)).
\end{equation}
\end{prop}

\subsection{Existence of weak solutions}

Now we are at the stage of showing the existence of weak solutions to \eqref{CPE'}. In particular, this proves the main theorem of this work, Theorem \ref{main-theorem}.
Starting with the continuity equation, let $ \psi $ be a scalar function in $ C^\infty_c(\overline\Omega\times[0,T)) $. Multiply \eqref{eq:appro-density} with $ \psi $ and integrate the resultant over $ \Omega\times(0,T) $. After integration by parts, we have
\begin{equation}\label{CV-001}
	\begin{aligned}
		& \int_\Omega\rho_{\varepsilon,0} \psi|_{t=0} \idx +  \int_0^T \int_\Omega \biggl( \rho_\varepsilon \dt \psi + \rho_\varepsilon v_\varepsilon \cdot\nablah \psi + \rho_\varepsilon w_\varepsilon \dz \psi \biggr) \idx\,dt \\
		& ~~~~ = \varepsilon \int_0^T \int_\Omega \biggl( - \rho_\varepsilon^{1/2} \deltah \rho_\varepsilon^{1/2} \psi
		 + \rho_\varepsilon^{1/2}\abs{\nablah \rho_\varepsilon^{1/2}}{2} \nablah \rho_\varepsilon^{1/2}\cdot\nablah \psi \\
		 & ~~~~ ~~~~ + \abs{\nablah \rho_\varepsilon^{1/2}}{4} \psi - \rho^{-p_0}_\varepsilon\psi \biggr) \idx \,dt.
	\end{aligned}
\end{equation}
On the one hand, the right-hand side of \eqref{CV-001} is bounded by
\begin{align*}
	& (\varepsilon^{1/2} + \varepsilon^{1/6} + \varepsilon^{1/(2p_0+1)})C_T\bigl(\norm{\psi}{L^\infty\spacetime} + \norm{\nabla \psi}{L^\infty\spacetime}\bigr)  \rightarrow 0, \\
	& ~~~~ ~~~~ \text{as} ~ \varepsilon \rightarrow 0^+.
\end{align*}
The above follows from inequalities \eqref{CM-101}--\eqref{CM-102}.
On the other hand, \eqref{CM-001-1} yields that, as $ \varepsilon \rightarrow 0^+ $, there is a $ \kappa \in L^2\spacetime $ such that
\begin{equation}\label{CV-002}
	\rho_\varepsilon^{1/2} w_\varepsilon \rightharpoonup \kappa, ~~~~~~ \text{weakly in} ~~ L^2 \spacetime.
\end{equation}
We define the vertical velocity as
\begin{equation}\label{CV-def-vertical-w}
	w(x,y,z,t) := \begin{cases}
		\kappa(x,y,z,t) / \rho^{1/2}(x,y,t) & \text{in}  ~~ \lbrace \rho > 0 \rbrace, \\
		0 & \text{in} ~~ \lbrace \rho = 0 \rbrace.
	\end{cases}
\end{equation}
Then it is easy to verify from \eqref{CV-002} that, as $ \varepsilon \rightarrow 0^+ $,
\begin{equation*}
	(\rho_\varepsilon^{1/2} w_\varepsilon - \rho^{1/2} w)\rho^{1/2} = (\rho_\varepsilon^{1/2} w_\varepsilon - \kappa)\rho^{1/2}  \rightharpoonup 0 ~~~~ \text{weakly in} ~ L^q \spacetime,
\end{equation*}
with $ q \in (1, 2) $, since $ \rho = \rho(x,y,t) $ and $$ \rho^{1/2} \in L^\infty(0,T;W^{1,2}(\Omega_h)) \subset L^\infty(0,T;L^p(\Omega)) , $$ for every $ p > 1 $.
Together with \eqref{CM-001}, \eqref{CM-001-1}, \eqref{CM-density-ae} and the Lebesgue Dominated Convergence Theorem, this yields, as $ \varepsilon \rightarrow 0^+ $,
\begin{equation*}
	\begin{aligned}
		& \norm{\rho_\varepsilon^{1/2} - \rho^{1/2}}{L^2\spacetime} \rightarrow 0, ~~~~ \text{and}\\
		& \int_0^T \int_\Omega (\rho_\varepsilon w_\varepsilon - \rho w )\dz \psi \idx\,dt = \int_0^T \int_\Omega (\rho_\varepsilon^{1/2} - \rho^{1/2}) \rho_\varepsilon^{1/2} w_\varepsilon \dz \psi \idx \,dt \\
		& ~~~~  + \int_0^T \int_\Omega (\rho_\varepsilon^{1/2} w_\varepsilon - \rho^{1/2} w) \rho^{1/2} \dz \psi \idx \,dt \leq \norm{\rho_\varepsilon^{1/2} - \rho^{1/2}}{L^2\spacetime}\\
		& ~~~~ ~~~~ \times \norm{\rho_\varepsilon^{1/2}w_\varepsilon}{L^2\spacetime} \norm{\dz \psi}{L^\infty(\Omega\times(0,T))} \\
		& ~~~~ + \int_0^T \int_\Omega (\rho_\varepsilon^{1/2} w_\varepsilon - \rho^{1/2} w) \rho^{1/2} \dz \psi \idx \,dt \rightarrow 0.
	\end{aligned}
\end{equation*}
Then \eqref{CM-001}, \eqref{CM-density}, \eqref{CM-008} together with the above analysis imply that by taking $ \varepsilon \rightarrow 0^+ $, \eqref{CV-001} converges to
\begin{equation}\label{weak-sol-001}
	\int_\Omega \rho_0 \psi|_{t=0} \idx + \int_0^T \int_\Omega \biggl( \rho \dt \psi + \rho v \cdot\nablah \psi + \rho w \dz \psi \biggr) \idx\,dt = 0,
\end{equation}
which shows \subeqref{CPE'}{1}.

Finally, to establish \subeqref{CPE'}{2}, let $ \phi $ be a vector-valued function in $ C_c^\infty(\overline\Omega\times[0,T); \mathbb R^2 ) $. Take inner product of \eqref{eq:approx-momentum} with $ \phi $ and integrate the resultant over $ \Omega \times (0,T) $. After integration by parts, we have
\begin{equation}\label{CV-003}
	\begin{aligned}
		& \int_\Omega \rho_{\varepsilon,0} v_{\varepsilon,0} \cdot\phi|_{t=0} \idx +  \int_0^T \int_\Omega \biggl( \rho_\varepsilon v_\varepsilon \cdot \dt \phi + \rho_\varepsilon v_\varepsilon \otimes v_\varepsilon : \nablah \phi \\
		& ~~~~ ~~~~ + \rho_\varepsilon w_\varepsilon v_\varepsilon \cdot \dz \phi
		+ \rho_\varepsilon^\gamma \dvh \phi \biggr)  \idx\,dt - \int_0^T \int_\Omega \biggl( \rho_\varepsilon \mathcal D(v_\varepsilon) : \nablah \phi \\
		& ~~~~ ~~~~ + \rho_\varepsilon \dz v_\varepsilon \cdot \dz \phi \biggr) \idx\,dt
		  = \varepsilon \int_0^T\int_\Omega \biggl( - \rho_\varepsilon^{1/2}\deltah \rho_\varepsilon^{1/2} v_\varepsilon \cdot \phi \\
		 & ~~~~ ~~~~
		 + \rho^{1/2} \abs{\nablah \rho_\varepsilon^{1/2}}{2} \nablah \rho_\varepsilon^{1/2} \otimes v_\varepsilon : \nablah \phi
		 + \abs{\nablah \rho_\varepsilon^{1/2}}{4}v_\varepsilon\cdot\phi \\
		  & ~~~~ ~~~~
		  + \rho_\varepsilon\abs{v_\varepsilon}{3}v_\varepsilon\cdot\phi  \biggr)  \idx\,dt
		+
			\varepsilon^{1/2} \int_0^T\int_\Omega
			\rho_\varepsilon \nablah v_\varepsilon : \nablah \phi \idx\,dt,
	\end{aligned}
\end{equation}
where the right-hand side is bounded by
\begin{align*}
	& (\varepsilon^{1/2} + \varepsilon^{1/6} + \varepsilon^{1/12} + \varepsilon^{1/5})C\bigl(\norm{\psi}{L^\infty\spacetime} + \norm{\nabla \psi}{L^\infty\spacetime}\bigr)  \rightarrow 0, \\
	& ~~~~ ~~~~ \text{as} ~~ \varepsilon \rightarrow 0^+,
\end{align*}
as a consequence of inequalities \eqref{CM-201}--\eqref{CM-202}.

On the other hand, as $ \varepsilon \rightarrow 0^+ $, \eqref{CM-001}, \eqref{CM-001-1}, \eqref{CM-energy}, \eqref{CV-002} and \eqref{CV-def-vertical-w} imply
\begin{align*}
	& \int_0^T \int_{\Omega} (\rho_\varepsilon w_\varepsilon v_\varepsilon - \rho w v) \cdot \dz \phi \idx\,dt = \int_0^T \int_\Omega \rho_\varepsilon^{1/2} w_\varepsilon (\rho_{\varepsilon}^{1/2}v_\varepsilon - \rho^{1/2} v)  \cdot \dz \phi \idx\,dt \\
	& ~~~~ + \int_0^T \int_\Omega (\rho_\varepsilon^{1/2} w_\varepsilon - \rho^{1/2} w)\rho^{1/2} v \cdot \dz \phi \idx\,dt = \int_0^T \int_\Omega \rho_\varepsilon^{1/2} w_\varepsilon (\rho_{\varepsilon}^{1/2}v_\varepsilon \\
	& ~~~~ ~~~~ - \rho^{1/2} v)  \cdot \dz \phi \idx\,dt + \int_0^T \int_\Omega (\rho_\varepsilon^{1/2} w_\varepsilon - \kappa)\rho^{1/2} v \cdot \dz \phi \idx\,dt  \rightarrow 0.
\end{align*}
In the meantime, $ \norm{\rho^{1/2}_\varepsilon \dz v_\varepsilon}{L^2(\Omega\times(0,T))} < C_T $, by \eqref{CM-001}, which implies that there is a $ \iota \in L^2(\Omega\times(0,T)) $ such that
\begin{equation}\label{CV:dz-horizontal-v}
	\rho^{1/2}_\varepsilon \dz v_\varepsilon \rightharpoonup \iota, ~~~~~~ \text{weakly in} ~~ L^2(\Omega\times(0,T)).
\end{equation}
As before, define
\begin{equation}\label{CV:def-dz-horizontal-v}
	\dz v(x,y,z,t) := \begin{cases}
	\iota(x,y,z,t) / \rho^{1/2}(x,y,t) & \text{in}  ~~ \lbrace \rho > 0 \rbrace, \\
	0 & \text{in} ~~ \lbrace \rho = 0 \rbrace.
	\end{cases}
\end{equation}
Then, as $ \varepsilon \rightarrow 0^+ $,
\begin{align*}
	& \int_0^T \int_\Omega (\rho_\varepsilon \dz v_\varepsilon - \rho \dz v ) \cdot \dz \phi \idx \,dt = \int_0^T \int_\Omega \rho_\varepsilon^{1/2}\dz v_\varepsilon (\rho_\varepsilon^{1/2} - \rho^{1/2}) \cdot \dz \phi \idx \,dt \\
	& ~~~~ + \int_0^T \int_\Omega (\rho_\varepsilon^{1/2} \dz v_\varepsilon - \rho^{1/2} \dz v) \rho^{1/2} \cdot \dz \phi \idx\,dt = \int_0^T \int_\Omega \rho_\varepsilon^{1/2}\dz v_\varepsilon (\rho_\varepsilon^{1/2} \\
	& ~~~~ - \rho^{1/2}) \cdot \dz \phi \idx \,dt  + \int_0^T \int_\Omega (\rho_\varepsilon^{1/2} \dz v_\varepsilon - \iota) \rho^{1/2} \cdot \dz \phi \idx\,dt \rightarrow 0.
\end{align*}
This argument equally holds when replacing $ \dz $ with $ \partial_h $, and hence we have, as $ \varepsilon \rightarrow 0^+ $,
\begin{equation}\label{CV:def-horizontal-derivative-v}
\begin{gathered}
	 \rho^{1/2}_\varepsilon \nablah v_\varepsilon \rightharpoonup \Xi, ~~~~~~ \text{weakly in} ~~ L^2(\Omega\times(0,T)), \\
	 \nablah v(x,y,z,t) := \begin{cases}
	\Xi(x,y,z,t) / \rho^{1/2}(x,y,t) & \text{in}  ~~ \lbrace \rho > 0 \rbrace, \\
	0 & \text{in} ~~ \lbrace \rho = 0 \rbrace,
	\end{cases} \\
	 \int_0^T \int_\Omega ( \rho_\varepsilon \mathcal D(v_\varepsilon) - \rho \mathcal D(v) ) : \nablah \phi \idx \,dt \rightarrow 0.
\end{gathered}
\end{equation}
We remark that in \eqref{CV:def-dz-horizontal-v} \eqref{CV:def-horizontal-derivative-v}, by taking the convergence in the subsequence of that in \eqref{CM-007}, the velocity $ v $ here is equal to that in \eqref{CM-def-horizontal-velocity}.

Therefore, the left-hand side of \eqref{CV-003} will converge according to \eqref{CM-001}, \eqref{CM-density}, \eqref{CM-density-star}, \eqref{CM-008} and \eqref{CM-energy}
. Thus, by taking $ \varepsilon \rightarrow 0^+ $, \eqref{CV-003} converges to
\begin{equation}\label{weak-sol-002}
	\begin{aligned}
		& \int_\Omega m_0 \cdot \phi|_{t=0} \idx + \int_0^T \int_\Omega \biggl( \rho v \cdot \dt \phi + \rho v \otimes v : \nablah \phi + \rho w v \cdot \dz \phi \\
		& ~~~~ ~~~~ + \rho^\gamma \dvh \phi \biggr) \idx\,dt - \int_0^T \int_\Omega \biggl( \rho \mathcal D(v) : \nablah \phi + \rho \dz v \cdot \dz \phi \biggr) \idx\,dt =0. \\
	\end{aligned}	
\end{equation}
Together with \eqref{weak-sol-001}, this shows the global existence of weak solutions to \eqref{CPE'}.

What is left is to establish the energy inequality. In fact, from \eqref{BE-001} and \eqref{CM-001}, one has
\begin{align*}
		& \dfrac{d}{dt} \bigl\lbrace \dfrac 1 2 \int \rho_\varepsilon \abs{v_\varepsilon}{2} \idx + \dfrac{1}{\gamma - 1} \int \rho_\varepsilon^\gamma \idx \bigr\rbrace + \int \rho_\varepsilon \abs{\mathcal{D}(v_\varepsilon)}{2} \idx + \int \rho_\varepsilon \abs{\dz v_\varepsilon}{2} \idx \\
		& ~~~~  + (\sqrt{\varepsilon} - \varepsilon) \int \rho_\varepsilon \abs{\nablah v_\varepsilon}{2} \idx \leq \dfrac{\varepsilon \gamma}{\gamma - 1} \int \rho_\varepsilon^{\gamma - p_0 - 1 }\idx \\
		&  \leq \dfrac{\varepsilon\gamma}{\gamma-1} \bigl(\int \rho_\varepsilon^{-p_0} \idx\bigr)^{\frac{p_0 - \gamma + 1}{p_0}}\leq C_\gamma \varepsilon^{\frac{\gamma-1}{p_0}} \sup_{0<t<T}\norm{\varepsilon^{1/{p_0}} \rho_\varepsilon^{-1}}{L^{p_0}}^{p_0-\gamma+1} \leq C_{\gamma, T} \varepsilon^{\frac{\gamma-1}{p_0}}.
\end{align*}
Then by taking $ \varepsilon \rightarrow 0^+ $ in the above inequality, we arrive at
\begin{equation}
	\begin{aligned}
		& \dfrac{d}{dt} \bigl\lbrace \dfrac 1 2 \int \rho \abs{v}{2} \idx + \dfrac{1}{\gamma - 1} \int \rho^\gamma \idx \bigr\rbrace + \int \rho \abs{\mathcal{D}(v)}{2} \idx \\
		& ~~~~ + \int \rho \abs{\dz v}{2} \idx \leq 0,
	\end{aligned}
\end{equation}
in $ \mathcal D'(0,T) $. The entropy inequality follows from the Bresch-Desjardins entropy inequality and the Mellet-Vasseur estimate in Proposition \ref{prop:uniform-est} after taking $ \varepsilon \rightarrow 0^+ $. The regularity follows from \eqref{CM-001},  \eqref{CM-density},
\eqref{CM-density-2}, \eqref{CM-density-star}, \eqref{CM-energy-2}, \eqref{CV-002}, \eqref{CV-def-vertical-w}, \eqref{CV:dz-horizontal-v}, \eqref{CV:def-dz-horizontal-v}, \eqref{CV:def-horizontal-derivative-v}.
This finishes the proof of Theorem \ref{main-theorem}.

\subsection*{Acknowledgement}

The authors would like to thank the \'{E}cole Polytechnique for its kind hospitality, where this work was completed, and the \'{E}cole Polytechnique Foundation for its partial financial support through the 2017-2018 ``Gaspard Monge Visiting Professor" Program. This work is supported in part by
 the ONR grant N00014-15-1-2333. The work of E.S.T. was also supported in part by the Einstein Stiftung/Foundation - Berlin, through the Einstein Visiting Fellow Program and by the John Simon Guggenheim Memorial Foundation.


\bibliographystyle{plain}


%
%
%

\end{document}